\documentclass[sglanonrev]{informs4_arxiv} 
\usepackage{eqndefns-left} 
\RequirePackage{tgtermes}
\RequirePackage{newtxtext}
\RequirePackage{newtxmath}
\RequirePackage{bm}
\RequirePackage{endnotes}

\OneAndAHalfSpacedXII 

\usepackage{algorithm}
\usepackage[noend]{algpseudocode}
\usepackage{tikz}

\usepackage{natbib}
 \bibpunct[, ]{(}{)}{,}{a}{}{,}%
 \def\bibfont{\small}%
\usepackage{booktabs} 
\usepackage{multirow}
\usepackage{hyperref}
\usepackage{subcaption}
\usepackage{amsmath}
\usetikzlibrary{arrows.meta, positioning}
\usepackage{fix-cm}
\usepackage{caption}
\usepackage{ifthen}
\usepackage{enumitem}
\usepackage{makecell}

\newcommand\rev[1]{\color{black}{#1}} 
\newenvironment{revblock}{\color{black}}{}  
\newcommand{\p}{\phantom{0}}
\newcommand{\pp}{\phantom{00}}

\newcommand{\Yes}{\checkmark}
\newcommand{\No}{}
\newcommand{\specialcell}[2][c]{%
  \begin{tabular}[#1]{@{}c@{}}#2\end{tabular}}
\newcommand{\pst}{\phantom{\text{s.t. }}}
\newcommand{\setspacing}{ 
  \setlength{\abovedisplayskip}{7.1pt}%
  \setlength{\belowdisplayskip}{4.6pt}%
  \setlength{\abovedisplayshortskip}{7.1pt}%
  \setlength{\belowdisplayshortskip}{4.6pt}%
} 

\newcommand{\epi}{\mathrm{epi}}

\EquationsNumberedThrough    

\TheoremsNumberedThrough     
\ECRepeatTheorems  %

\MANUSCRIPTNO{MOOR-0001-2024.00}

\begin{document}


\RUNAUTHOR{Legault et al.} 

\RUNTITLE{Bus Fleet Electrification Planning Through Logic-Based Benders Decomposition and Restriction Heuristics}

\TITLE{Bus Fleet Electrification Planning Through Logic-Based Benders Decomposition and Restriction Heuristics}

\ARTICLEAUTHORS{%

\AUTHOR{Robin Legault}
\AFF{Operations Research Center, Massachusetts Institute of Technology, \EMAIL{legault@mit.edu}}

\AUTHOR{Filipe Cabral}
\AFF{Georgia Institute of Technology, \EMAIL{fcabral1290@gmail.com}}

\AUTHOR{Xu Andy Sun}
\AFF{Sloan School of Management, Massachusetts Institute of Technology, \EMAIL{sunx@mit.edu}}
} 

\ABSTRACT{%
To meet sustainability goals and regulatory requirements, transit agencies must plan partial and full transitions to electric bus fleets. This paper presents a comprehensive and computationally efficient multi-period optimization framework integrating the key decisions required to support such electrification initiatives. Our model is formulated as a two-stage integer program with integer subproblems. These two levels optimize, respectively, yearly fleet sizing and charging infrastructure investments, and hourly vehicle scheduling and charging operations. We develop an exact logic-based Benders decomposition algorithm enhanced by several acceleration techniques, including preprocessing, master problem strengthening, and efficient cut separation methods applied to different relaxations of the operational problem. These accelerations achieve speedups of up to three orders of magnitude and provide new theoretical and practical insights into Benders cut selection. We also propose a heuristic tailored for long-term, citywide electrification planning. This approach imposes and progressively relaxes additional scheduling constraints derived from auxiliary problems. It delivers high-quality solutions with optimality gaps below 1\% for instances an order of magnitude larger than those considered in prior work. We illustrate our model using real data from the Chicago public bus system, providing managerial insights into optimal investment and operational policies.
}%



\maketitle

\EquationErrorMsgdimen=500pt
\EquationErrorBoxdimen=500pt
\setlength{\maxdseq}{500pt} 

\section{Introduction}
\label{sec:Intro}
Over the last decade, major transit agencies \begin{revblock}worldwide\end{revblock} have established ambitious electrification targets for their bus fleets. In 2017, 35 cities across six continents pledged to procure only zero-emission buses starting from 2025 \citep{C40_2023}\begin{revblock}, and several bus systems, particularly in China and Europe, are already largely or fully electric \citep{IEA2025}\end{revblock}. In the United States, cities such as Boston, Chicago, and New York have committed to operating fully electrified fleets by 2040, and all public transit agencies in California are mandated to replace conventional buses with zero-emission models by this time \begin{revblock}\citep{WSP2024}\end{revblock}. While full electrification is the stated goal, practical considerations such as funding availability, operational complexity, and infrastructure readiness also lead agencies to initiate partial electrification projects targeting selected routes or depots \citep[e.g.,][]{CTA_Electrification}. Despite these initiatives, \begin{revblock}most of the transition still lies ahead. As of July 2025, 7,261 full-size battery electric buses (BEBs) had been funded, ordered, delivered, or deployed nationwide \citep{CALSTART2026}, about one tenth of the country's roughly 70,000 transit buses \citep{APTA2026}. The scale of the intended electrification of public bus systems makes effective planning tools essential to achieving this transition in a timely and cost-effective manner.\end{revblock}

An electrification plan involves interdependent decisions across different time scales. At the strategic level, transit agencies must determine fleet composition, infrastructure deployment, and the phasing of these investments over multiple planning periods. At the operational level, they must ensure that their fleet can provide the required service while respecting range limitations and charging dynamics of BEBs. Recent years have seen significant progress in optimization techniques for the operational scheduling of electric bus fleets, which construct deployable schedules for a fixed fleet and charging infrastructure under precise operational constraints. \begin{revblock}This literature includes the electric vehicle scheduling problem (EVSP) \citep[e.g.,][]{parmentier2023electric, de2024electric}, which assigns timetabled trips and charging times to vehicles, and the charge scheduling problem \citep{klein2023electric}, which assumes a fixed trip-to-vehicle assignment and optimizes charging under flexible departure times and detailed battery dynamics\end{revblock}. In contrast to operational planning, strategic planning for bus fleet electrification has received limited attention from the optimization community \citep{perumal2022electric, zhou2024charging}. Existing approaches either rely on highly simplified operational models or provide very limited scalability. To bridge this gap, we propose an optimization framework for bus fleet electrification planning that balances operational realism with computational tractability, along with exact and heuristic algorithms that can solve instances of practical size to proven or near optimality. We summarize our contributions as follows.

\subsection*{Contributions}
\begin{enumerate}[align=left, listparindent=\parindent, parsep=0pt]
    \item \textbf{Modeling.} We introduce a comprehensive model for bus fleet electrification planning (Section~\ref{sec:Model}). At the strategic level, we model the procurement of BEBs, retirement of conventional buses, and placement of chargers over multiple investment periods. At the operational level, we model the hourly operations of the fleet through a flow-based formulation that tracks the state of charge at the individual vehicle level. This model provides more realism and flexibility than those classically used in strategic planning. Using historical service schedules and geospatial data from eight major US transit agencies, we construct a set of realistic benchmark instances for the electrification of bus networks.
    \item \textbf{Algorithm design and computation.}
    \begin{itemize}[wide=0pt]
        \item \textbf{Exact method.} To solve small instances representative of partial electrification projects, we develop a logic-based Benders decomposition accelerated with preprocessing, master problem strengthening, and custom cuts obtained from different relaxations of the operational model (Sections~\ref{sec:Algorithms:Preprocessing}--\ref{sec:Algorithms:Benders}). Moreover, we establish the equivalence of two recently proposed Benders cut selection techniques (Section \ref{sec:Algorithms:Benders:linear}). Our accelerations achieve speedups of three orders of magnitude compared to a baseline logic-based Benders decomposition algorithm (Section~\ref{sec:Experiments:aceclerationLBBD}), yielding an exact method that significantly outperforms Gurobi (Section~\ref{sec:Experiments:small}).
        \item \textbf{Heuristic method.} For large-scale instances, we propose an easily implementable heuristic that solves a sequence of problems in which restrictions on the schedule of BEBs are imposed and progressively relaxed (Section~\ref{sec:Algorithms:Heuristic}). This algorithm achieves optimality gaps of around 1\% for 10-year-long electrification planning instances defined on citywide networks with more than 100 candidate charging locations, 100 routes, and 1000 buses (Section~\ref{sec:Experiments:large}).
    \end{itemize}
    \item \textbf{Case study.} To illustrate the scalability of our model and the type of insights it can provide for transit agencies, we present a case study of the Chicago bus network (Section~\ref{sec:Experiment:chicago}). We analyze the optimal sequencing of investments, the spatial deployment of chargers, and trends in bus utilization throughout the electrification process. Under our stated assumptions, our findings suggest prioritizing early electrification of high-usage routes supported by fast chargers, while completing electrification in less dense areas with BEBs relying on depot charging.
\end{enumerate}

The remainder of the paper consists of a literature review (Section~\ref{sec:Literature}) and a conclusion (Section~\ref{sec:Conclusion}).

\section{Related works}
\label{sec:Literature}

Bus fleet electrification planning encompasses decisions on fleet composition, infrastructure investments, and their phasing over time. Table \ref{tab:literature_bus} reviews relevant publications based on the strategic and operational decisions they model and the scale of instances considered.

Early works on bus system electrification focused on strategic decisions without explicitly modeling operations. \citet{xylia2017locating} optimize the location of fast charging stations across a large network, but assume that charging demand is exogenously given at each location. \citet{pelletier2019electric} introduce a multi-period model for fleet electrification that determines the mix of electric and conventional buses over time. They account for different charging technologies, but assume that buses follow predefined schedules with known energy consumption patterns.

The first model to jointly optimize BEB fleet composition, charging infrastructure, and vehicle scheduling was proposed by \citet{rogge2018electric}. Their formulation constructs a heterogeneous fleet, selects a number of depot chargers, and designs schedules to cover a set of timetabled trips. \citet{wang2022integrated} and \citet{he2023joint} extend this framework to optimize the deployment of fast on-route chargers for opportunity charging during layovers at terminal stations. \begin{revblock} \citet{nath2024impact} co-optimize charging-station locations, vehicle scheduling, and charging scheduling for a homogeneous fleet. Their work highlights the value of integrating charging scheduling into strategic planning, reporting average operational and total cost savings of 14.1\% and 4.1\%, respectively, compared to a sequential approach in which charging scheduling is optimized after the other decisions have been fixed.\end{revblock} \citet{liu2021optimizing} optimize similar decisions, but model the hourly flow of vehicles between each route and charging station to construct bus schedules without relying on predefined timetabled trips. 

\begin{table}[t]
\caption{Selected publications on bus fleet electrification planning}
\centering
\makebox[\textwidth]{%
\resizebox{1.0\textwidth}{!}{%
\renewcommand{\arraystretch}{0.75}
\begin{tabular}{lcccccccrrr} 
    \toprule
    & \multicolumn{3}{c}{Strategic decisions} & \multicolumn{4}{c}{Operational decisions} & \multicolumn{3}{c}{Largest instance} \\
    \cmidrule(lr){2-4} \cmidrule(lr){5-8} \cmidrule(lr){9-11}
    & \multirow{2}{*}{\specialcell{\rev{Heterogenous} \\ \rev{fleet composition}}} & \multirow{2}{*}{\specialcell{Charger \\ placement}} & \multirow{2}{*}{\specialcell{Multi-period \\ investments}} & \multirow{2}{*}{\specialcell{Depot\\ charging}} & \multirow{2}{*}{\specialcell{On-route \\ charging}} & \multirow{2}{*}{\specialcell{Charging \\ scheduling}} & \multirow{2}{*}{\specialcell{Vehicle \\ scheduling}} & \multirow{2}{*}{\specialcell{Routes \\ \phantom{..} }} & \multirow{2}{*}{\specialcell{Buses \\ \phantom{..} }} & \multirow{2}{*}{\specialcell{Charging \\ locations }}  \\ \\
    \midrule
    Xylia et al.~(\citeyear{xylia2017locating})        & \No & \Yes & \No  & \No & \Yes & \No & \No & 143 & - & 403 \\ 
    Rogge et al.~(\citeyear{rogge2018electric})        & \Yes & \No & \No  & \Yes & \No & \Yes & \Yes & 3 & 14 & 1  \\
    Pelletier et al.~(\citeyear{pelletier2019electric})    & \Yes & \No & \Yes  & \Yes & \Yes & \No & \No & - & 100 & - \\ 
    Liu et al.~(\citeyear{liu2021optimizing})        & \No  & \Yes & \No  & \Yes  & \Yes & \Yes & \Yes & 17 & 252 & 15 \\
    Dirks et al.~(\citeyear{dirks2022integration})     & \Yes & \Yes & \Yes & \Yes & \Yes & \Yes & \No & 60 & 357 & 182 \\
    Hu et al.~(\citeyear{hu2022joint})              & \Yes & \Yes & \No  & \No & \Yes & \Yes & \No & 3 & 16 & 111 \\
    Wang et al.~(\citeyear{wang2022integrated})       & \No & \Yes & \No  & \Yes & \Yes & \Yes & \Yes & 4 & 28 & 31 \\
    He et al.~(\citeyear{he2023joint})              & \Yes & \Yes & \No  & \Yes & \Yes & \Yes & \Yes & 3 & 6 & 3 \\
    He et al.~(\citeyear{he2023time})               & \Yes & \Yes & \Yes & \Yes & \Yes & \Yes & \No  & 36 & 170 & 29 \\
    Gairola et al.~(\citeyear{gairola2023optimization}) & \No & \Yes & \No & \Yes & \Yes & \Yes & \No & 18 & 285 & 21 \\
    \rev{Nath et al.~(\citeyear{nath2024impact})}                       & \No & \rev{\Yes} & \No & \No  & \rev{\Yes}  & \rev{\Yes}  & \rev{\Yes}  & \rev{34} & \rev{183}  & \rev{58} \\
    This work                       & \Yes & \Yes & \Yes & \Yes & \Yes & \Yes & \Yes & 140 & 1817 & 200 \\
    \bottomrule
\end{tabular}
}}
\label{tab:literature_bus}
\end{table}

Other studies focus on charger placement while assuming vehicle service patterns are given. \citet{hu2022joint} and \citet{gairola2023optimization} design charging infrastructure and charging schedules to meet the fleet's energy needs. The models of \citet{he2023time} and \citet{dirks2022integration} address strategic planning as a multi-period process, thereby recognizing that fleet electrification is typically phased over time rather than accomplished in a single step. For computational tractability, these studies assign each new BEB to a service schedule currently performed by a conventional bus, substituting a scheduling problem for a simple assignment decision for each vehicle. While this modeling choice greatly reduces problem complexity, it comes at the cost of ignoring opportunities to optimize vehicle schedules in light of BEB range limitations and charging dynamics.

Existing literature on bus fleet electrification planning tends to prioritize new modeling approaches and case-study insights, with algorithmic development often treated as secondary. As a result, most studies either omit vehicle scheduling or are limited to small systems. \begin{revblock}Among the works that integrate strategic decisions and vehicle scheduling, we are not aware of any\end{revblock} that solve instances to proven optimality. \citet{rogge2018electric} and \citet{he2023joint} use heuristic genetic algorithms, \citet{liu2021optimizing} solve a surrogate relaxation without providing global optimality bounds, \begin{revblock}\citet{nath2024impact} develop an iterated local search heuristic in which charging-scheduling subproblems are addressed via surrogate linear programs\end{revblock}, and \citet{wang2022integrated} solve their model in extensive form, obtaining an optimality gap of 4\% after 25 hours of computation on a four-route instance.

To address these realism and scalability limitations, we develop a multi-period framework that optimizes fleet management and charger placement while integrating vehicle scheduling through a flow-based formulation defined on a time-expanded graph of state of charge. In contrast to the flow-based operational model of \citet{liu2021optimizing}, \begin{revblock}which aggregates the state of charge across all buses assigned to the same route\end{revblock}, our formulation tracks the state of charge of individual vehicles, capturing range limitations and charging dynamics more faithfully. \begin{revblock}Similar flow-based formulations have been used in planning applications for other electric mobility systems. \citet{zhang2019vehicle} introduced three-dimensional networks tracking location, time, and battery level to optimize charging operations in one-way electric carsharing, and \citet{santos2023space} employed this construct for the strategic design of a regional shared automated electric vehicle system, solving the resulting model in extensive form with a general-purpose solver. We instead exploit the structure of our formulation through an exact logic-based Benders decomposition and a tailored restriction heuristic.\end{revblock} This allows us to solve instances representative of realistic small-scale electrification projects to proven optimality, and to obtain high-quality solutions for systems an order of magnitude larger than those considered in previous studies on bus fleet electrification planning.

\section{Problem formulation}
\label{sec:Model}
This section presents our formulation of the bus fleet electrification problem (BFEP). Section \ref{sec:Model:assumptions} delimits the scope of this work and discusses our modeling assumptions. The mathematical formulations of the strategic and operational problems are presented in Sections \ref{sec:Model:strategic} and \ref{sec:model:operations}. Key properties of the BFEP are discussed in Section \ref{sec:model:properties}. A notation table is provided in \ref{ec:notation}.

\subsection{Scope and modeling assumptions}
\label{sec:Model:assumptions}
Our objective is to develop a comprehensive optimization framework for bus fleet electrification planning that integrates fleet composition and infrastructure deployment over a multi-year horizon. To inform these strategic decisions, we incorporate an operational model that captures the essential trade-offs between investment and operational expenses. To balance modeling realism, tractability, and clarity of exposition, we make simplifying assumptions that we outline next.

\begin{itemize}
    \item \textbf{Electrification targets, investment budget, and transition timeline.} We divide the planning horizon into periods (e.g., years), in which the fleet composition and charging infrastructure can be updated subject to budget constraints. Electrification targets are enforced through lower bounds on the number of acquired BEBs and upper bounds on the number of conventional buses still in operation. The latter can be set to zero in the final period to enforce full electrification by a target date, in line with regulatory mandates \citep[e.g.,][]{CARB2018} and agency commitments \citep[e.g.,][]{MBTA_Electrification}. This formulation translates policy goals and budget limitations into capacity constraints while leaving the timing and composition of investments endogenous.

    \item \textbf{BEB types and charging policies.} 
    BEBs available for purchase are partitioned into two categories: (i) \emph{depot BEBs}, which charge exclusively at depot locations using plug-in chargers; and (ii) \emph{on-route BEBs}, which rely on fast chargers installed at route terminals for opportunity charging during layovers. \begin{revblock}In addition, we assume that depot BEBs never interrupt charging before reaching their full state of charge. We adopt these assumptions primarily for presentation and modeling clarity. In \ref{ec:partial-charging}, we show how our framework can be extended to relax the distinction between depot and on-route BEBs and allow partial charging.\end{revblock}

    \item \textbf{Assignment of vehicles to routes.} In each investment period, we assign vehicles to a single bus route and do not allow interlining. This assumption is standard in planning models \citep[e.g.,][]{liu2021optimizing, he2023time} and is a reasonable approximation of real transit operations, where routes often operate with dedicated fleets \citep{zahedi2025dynamic}. Depot chargers and fast chargers at bus terminals are shared across routes.

    \item \textbf{Representative day, time discretization, and cyclicity.} As in previous multi-period electrification models \citep[e.g.,][]{dirks2022integration, he2023time}, we describe operations in each investment period using a representative day, divided into hourly intervals. Unlike existing strategic models that impose the same daily schedule for each individual bus, we only assume cyclicity at the fleet level: what repeats daily in our model is how many buses perform each service, idling, and charging activity in each possible state of charge. Our model thus allows individual buses to follow staggered multi-day rotations, consistent with recent multi-day BEB scheduling studies \citep{vende2023matheuristics}. \rev{\ref{ec:Example_operational_model} illustrates that relaxing bus-level cyclicity to fleet-level cyclicity can reduce the number of buses needed to provide a given service level.}

    \item \begin{revblock}\textbf{Deterministic parameters.} Throughout the paper, we treat all parameters, including technology and energy-cost trajectories, as deterministic. When the model is solved, these parameters should be set using the planner's best available forecasts. \ref{ec:rolling_horizon} illustrates how deploying the model in a rolling horizon keeps future electrification decisions revisable as new information, such as changes in market conditions, is observed during the transition.\end{revblock}
\end{itemize}

These assumptions position the BFEP as a decision-support tool for determining the pace of electrification, fleet composition, charging infrastructure deployment, and representative operational costs under policy and budget constraints. Next, we formulate the BFEP as a two-stage program. The first stage determines the fleet composition and charger deployment in each investment period, and the second stage models the operations of the fleet.

\subsection{Strategic problem}
\label{sec:Model:strategic}

In each investment period $p \in \mathcal{P} = \{1,\dots,P\}$, the state $x_p = \left(\chi_p,\{\eta_{pr}\}_{r \in \mathcal{R}}\right) \in \mathbb{Z}_{+}^n$ of the system is described by charger variables $\chi_p = (\bar{\chi}^p, \widetilde{\chi}^p)$ and vehicle assignment variables $\eta_{pr} = (\bar{\eta}^p_r , \widetilde{\eta}^p_r, \widehat{\eta}^p_r)$ for each bus route $r \in \mathcal{R}$. The number of depot BEBs of each type $b \in \mathcal{B}$, of on-route BEBs, and of conventional buses assigned to route $r \in \mathcal{R}$ are respectively controlled by the decision variables $\bar{\eta}_{rb}^{p}$, $\widetilde{\eta}_{r}^{p}$, and $\widehat{\eta}_{r}^{p}$. The number of chargers installed at depot $i \in \mathcal{I}$ is given by $\bar{\chi}^p_{i}$, and $\widetilde{\chi}^p_{j}$ denotes the number of on-route chargers at terminal $j \in \mathcal{J}$. The BFEP is formulated as follows:
\begingroup
\setspacing
\allowdisplaybreaks
\begin{subequations}\label{model:BFEP}
\begin{align}\label{model:BFEP:invest_obj}
&\min_{x_p, \ \forall p \in \mathcal{P} } \quad \sum_{p \in \mathcal{P}} \gamma^p
\left(I_p(x_p-x_{p-1})
+ O_p(x_p) \right) \hspace{-3cm} \\
\label{model:BFEP:invest_budget_constr}
\text{s.t.} \quad &
I_p(x_p-x_{p-1}) \leq I_p^{\text{UB}},
&& \forall p \in \mathcal{P}, \\[.3ex]
\label{model:BFEP:invest_charg_mono}
& \bar{\chi}_{i}^p \geq \bar{\chi}_{i}^{p - 1},  \quad
\widetilde{\chi}_j^p \geq \widetilde{\chi}_j^{p-1},
&& \forall p \in \mathcal{P},\ i \in \mathcal{I}, \ j \in\mathcal{J}, \\
\label{model:BFEP:invest_fleet_mono}
& \sum_{r \in \mathcal{R}}\bar{\eta}_{rb}^p \geq \sum_{r \in \mathcal{R}}\bar{\eta}_{rb}^{p - 1}, \quad
\sum_{r \in \mathcal{R}}\widetilde{\eta}_{r}^p \geq
\sum_{r \in \mathcal{R}}\widetilde{\eta}_{r}^{p-1}, \quad
&& \forall p \in \mathcal{P},\  b \in \mathcal{B}, \\
\label{model:BFEP:invest_fleet_mono_conv}
& \sum_{r \in \mathcal{R}}\widehat{\eta}_{r}^p \leq \sum_{r \in \mathcal{R}}\widehat{\eta}_r^{p-1},\quad
&& \forall p \in \mathcal{P}, \\
\label{model:BFEP:aquisition_and_retirement}
& \sum_{r \in \mathcal{R}}\sum_{b \in \mathcal{B}}\bar{\eta}_{rb}^p + \sum_{r \in \mathcal{R}}\widetilde{\eta}_r^p \geq \eta_p^{\text{LB}}, \quad \sum_{r \in \mathcal{R}}\widehat{\eta}_r^p \leq \widehat{\eta}_p^{\text{UB}},
&& \forall p \in \mathcal{P}, \\
\label{model:BFEP:invest_charger_bnds}
& \bar{\chi}_{i}^p \leq \bar{\chi}^{\text{UB}}_{i},
\quad \widetilde{\chi}_j^p \leq \widetilde{\chi}^{\text{UB}}_{j},
&& \forall p \in \mathcal{P},\ i \in \mathcal{I}, \ j \in\mathcal{J}, \\
\label{model:BFEP:charger_variables}
& \bar{\chi}^p \in \mathbb{Z}^{\mathcal{I}}_+,\ \widetilde{\chi}^p \in \mathbb{Z}^{\mathcal{J}}_{+}, 
&& \forall p \in \mathcal{P},\\
\label{model:BFEP:fleet_variables}
& \bar{\eta}^p_r \in \mathbb{Z}^{\mathcal{B}}_+, \ \widetilde{\eta}^p_r \in \mathbb{Z}_+, \
\widehat{\eta}^p_r \in \mathbb{Z}_+,
&& \forall p \in \mathcal{P}, \ r \in \mathcal{R}. 
\end{align}
\end{subequations}
\endgroup

A discount factor $\gamma \in (0,1]$ and the initial state $x_0$ are given. The objective \eqref{model:BFEP:invest_obj} is to minimize the time-discounted investment and operational costs over the planning horizon. The investment costs $I_p(x_p-x_{p-1})$ are assumed to depend linearly on the number of each type of assets acquired and sold in period $p$. The function $O_p(x_p) := H_p(x_p) + \mathcal{Q}_p(x_p)$ denotes the cost of maintaining and operating the system in state $x_p$ during period $p$. It is composed of a nondecreasing linear function $H_p(x_p)$ modeling fixed maintenance costs and the optimal value $\mathcal{Q}_p(x_p)$ of the fleet scheduling problem. Constraints \eqref{model:BFEP:invest_budget_constr} impose an investment budget $I_p^{\text{UB}}$ for each period $p \in \mathcal{P}$. Constraints \eqref{model:BFEP:invest_charg_mono}--\eqref{model:BFEP:invest_fleet_mono} ensure that the chargers and BEBs acquired at any period remain in the system during the following periods, whereas \eqref{model:BFEP:invest_fleet_mono_conv} forces the number of conventional buses to be nonincreasing. Constraints \eqref{model:BFEP:aquisition_and_retirement} impose acquisition and retirement targets for BEBs and conventional buses, respectively. Constraints \eqref{model:BFEP:invest_charger_bnds} specify the charger hosting capacity of each location. Constraints \eqref{model:BFEP:charger_variables}--\eqref{model:BFEP:fleet_variables} give the domain of the strategic variables. We denote by $\mathcal{X}$ the set of feasible solutions to problem \eqref{model:BFEP}, and by $\mathcal{X}_p$ the projection of $\mathcal{X}$ onto the space of the strategic variables $x_p$ of period $p \in \mathcal{P}$.

\subsection{Operational problem}\label{sec:model:operations}
In each period $p \in \mathcal{P}$, the operational problem is to construct a minimum-cost schedule that satisfies hourly service level requirements. Time is discretized into a set of intervals $\mathcal{T} = \{0,1,\dots,T-1\}$, viewed as the residue classes modulo~$T$ to encode the daily cyclicity of the fleet's operations.

The operations of on-route BEBs and conventional buses servicing each route $r \in \mathcal{R}$ are modeled as simple assignments, \begin{revblock}without explicitly tracking their state of charge or fuel level\end{revblock}. While in service, an on-route BEB is assigned to a terminal $j \in \mathcal{J}(r) \subseteq \mathcal{J}$, where it can charge during its layover times, with $\mathcal{R}(j) = \{r \in \mathcal{R} : j \in \mathcal{J}(r)\}$ denoting the set of routes connected to $j$. This assignment is time-dependent, allowing charging at different terminals across intervals. \begin{revblock}We model fast-charger capacity through a throughput parameter $\rho$, defined as the maximum number of on-route BEBs in operation that each fast charger can support. This approximation captures opportunity charging at terminals by assuming that layover charging opportunities can be sufficiently sequenced over the course of service.\end{revblock} The number of on-route BEBs from route $r$ assigned to terminal $j$ during interval $t \in \mathcal{T}$ is denoted by $\widetilde{w}^{pt}_{rj}$, and $\widehat{w}^{pt}_{r}$ denotes the number of conventional buses in service at time $t$. 

Depot BEBs differ in that their schedules are represented as circulations on cyclic time-expanded graphs of state of charge (see \ref{ec:Example_operational_model} for an illustrative example). Specifically, for each route $r \in \mathcal{R}$ and each type of bus $b \in \mathcal{B}$, we maintain a graph with nodes $(t,s) \in \mathcal{T} \times \{0,\dots,s_b\}$, where $s=0$ and $s=s_b$ respectively represent the fully depleted and fully charged states. In each interval, a depot BEB can service a route, idle, or initiate a charging trip to one of the depots $i \in \mathcal{I}$ equipped with plug-in chargers. The flow variables $w^{pt}_{rbs}, v^{pt}_{rbs}$, and $\{z^{pt}_{rbis}\}_{i \in \mathcal{I}}$ respectively represent the number of vehicles initiating each possible service, idling, and charging operations from node $(t,s)$. Working for one interval reduces the state of charge by one unit, idling does not affect the battery level, and $\kappa_{rbis}$ intervals are needed to perform a round-trip from route $r$ to depot $i$ and fully recharge from state $s$. We denote by $\mathcal{S}_b^w = \{1,2,\dots,s_b\}$, $\mathcal{S}_b = \{0,1,\dots,s_b\}$ and $\mathcal{S}_b^z = \{0,1,\dots,s_b-1\}$ the states of charge from which service, idling and charging operations can respectively be initiated. 

\begin{revblock}Under these dynamics, the operational problem determines a minimum-cost schedule that meets the required service level, with $d_r^{pt}$ denoting the minimum number of buses in service on route $r \in \mathcal{R}$ during interval $t \in \mathcal{T}$.\end{revblock} Denoting by $y_{pr} = (w^p_r, v^p_r, z^p_r, \widetilde{w}^p_r, \widehat{w}^p_r) \in \mathbb{Z}_+^{m_{pr}}$ the decision variables that pertain to each route $r \in \mathcal{R}$, \begin{revblock}and by $c^y_{pr} \geq 0$ the corresponding vector of operational costs, including energy costs, driver wages, and variable maintenance costs,\end{revblock} the problem can be expressed as: 
\begingroup 
\setspacing
\allowdisplaybreaks
\begin{subequations}\label{model:operations_detailed}
\begin{align}\label{model:operations_detailed:invest_obj}
&\mathcal{Q}_p(x_p) := \min_{y_{pr}, \ \forall r \in \mathcal{R} } \quad  \sum_{r \in \mathcal{R}} c^{y\top}_{pr} y_{pr} \hspace{-7cm} \\
\label{model:operations_detailed:depot_capacity} \text{s.t.} \quad &
\sum_{b\in \mathcal{B}}\sum_{r\in \mathcal{R}}\sum_{s\in\mathcal{S}_b^z}\sum_{\substack{l = 0}}^{\kappa_{rbis}-1} z_{rbis}^{p(t - l)} \leq \bar{\chi}_{i}^{p}
&& \forall \ i \in \mathcal{I}, \ t \in \mathcal{T}, \\
\label{model:operations_detailed:terminal_capacity} & \sum_{r\in \mathcal{R}(j)} \widetilde{w}_{rj}^{pt} \leq \rho \widetilde{\chi}_{j}^p,
&& \forall j \in\mathcal{J}, \ t \in \mathcal{T}, \\
\label{model:operations_detailed:service_level} & \sum_{b\in \mathcal{B}}\sum_{s\in\mathcal{S}^w_b} w_{rbs}^{p t} +
\sum_{j\in \mathcal{J}(r)} \widetilde{w}_{rj}^{pt} + \widehat{w}_{r}^{pt} \geq d_{r}^{pt},
&& \forall r \in\mathcal{R}, \ t \in \mathcal{T}, \\
\label{model:operations_detailed:depot_flow_full} & v_{rbs}^{p t} + w_{rbs}^{p t} = v_{rbs}^{p (t-1)} + \sum_{i\in \mathcal{I}}\sum_{s'\in \mathcal{S}^z_b} z_{rbis'}^{p(t-\kappa_{rbis'})},
&& \forall r \in \mathcal{R}, \ b \in\mathcal{B}, \ t \in \mathcal{T}, \ s = s_b, \\
\label{model:operations_detailed:depot_flow_intermediate} & v_{rbs}^{p t} + w_{rbs}^{p t} + \sum_{i\in \mathcal{I}} z_{rbis}^{p t} = v_{rbs}^{p (t-1)} + w_{rb(s+1)}^{p (t-1)},
&& \hspace{-0.7cm}\forall r \in \mathcal{R}, \ b \in\mathcal{B}, \ t \in \mathcal{T}, \ s \in \mathcal{S}_b{\setminus}\{0,s_b\},  \\
\label{model:operations_detailed:depot_flow_empty} & v_{rbs}^{p t} + \sum_{i\in \mathcal{I}} z_{rbis}^{p t} = v_{rbs}^{p (t-1)} + w_{rb(s+1)}^{p (t-1)},
&& \forall r \in \mathcal{R}, \ b \in\mathcal{B}, \ t \in \mathcal{T}, \ s = 0,  \\
\label{model:operations_detailed:fleet_depot} & \sum_{s \in \mathcal{S}_b} v_{rbs}^{p 0} + \sum_{s \in \mathcal{S}_b^w} w_{rbs}^{p 0} + \sum_{i\in \mathcal{I}} \sum_{s \in \mathcal{S}_b^z}\sum_{\substack{l = 0}}^{\kappa_{rbis}-1} z_{rbis}^{p(- l)} \leq \bar{\eta}_{rb}^{p},
&& \forall r \in\mathcal{R}, \ b \in \mathcal{B},  \\
\label{model:operations_detailed:fleet_conv_route} & \sum_{j\in\mathcal{J}(r)} \widetilde{w}_{rj}^{pt} \leq \widetilde{\eta}_{r}^{p}, \quad \widehat{w}_{r}^{pt} \leq \widehat{\eta}_{r}^{p},
&& \forall r \in\mathcal{R}, \ t \in \mathcal{T},  \\
\label{model:operations_detailed:var_depot} & v^{p}_r \in \mathbb{Z}_+^{\mathcal{T}\times\mathcal{B}\times\mathcal{S}_b},
\ w^{p}_r\in \mathbb{Z}_+^{\mathcal{T}\times\mathcal{B}\times\mathcal{S}^w_b}, \ z^{p}_r \in \mathbb{Z}_+^{\mathcal{T}\times\mathcal{B}\times\mathcal{I}\times\mathcal{S}^z_b},
&& \forall r \in\mathcal{R},  \\
\label{model:operations_detailed:var_conv_route} & \widetilde{w}^{p}_r \in \mathbb{Z}_+^{\mathcal{T}\times\mathcal{J}(r)},
\ \widehat{w}^{p}_r \in \mathbb{Z}_+^{\mathcal{T}},
&& \forall r \in\mathcal{R}.
\end{align}
\end{subequations}
\endgroup

Constraints \eqref{model:operations_detailed:depot_capacity} ensure that the number of BEBs simultaneously in charge at a depot does not exceed the number of installed chargers. The charger usage at depot $i \in \mathcal{I}$ during interval $t \in \mathcal{T}$ is given by the number of BEBs of any type $b \in \mathcal{B}$ that initiated a charging operation from any route $r \in \mathcal{R}$ and any state of charge $s \in \mathcal{S}_b^z$ in the last $\kappa_{rbis}$ intervals. Constraints \eqref{model:operations_detailed:terminal_capacity} ensure that on-route charger usage respects the installed capacity at each terminal $j \in \mathcal{J}$. Constraints \eqref{model:operations_detailed:service_level} \begin{revblock}enforce the service level requirement\end{revblock} on route $r \in \mathcal{R}$ in interval $t\in \mathcal{T}$. Constraints \eqref{model:operations_detailed:depot_flow_full}--\eqref{model:operations_detailed:depot_flow_empty} model the depot BEB charging dynamics for each route $r \in \mathcal{R}$ and bus type $b \in \mathcal{B}$. In each flow-balance equation, the left- and right-hand sides respectively correspond to the outflow and inflow at a node $(t,s) \in \mathcal{T} \times \mathcal{S}_b$. For the fully charged state $s=s_b$, the outflow is the number of vehicles in service or idling in state $s$ at time $t$, and the inflow is the number of vehicles that were idling in state $s$ in interval $t-1$, or initiated, $\kappa_{rbis'}$ intervals earlier, a charging trip to depot $i \in \mathcal{I}$ from state $s' \in \mathcal{S}^z_b$, and thus become available and fully charged at time $t$. A similar logic applies to the partly depleted states $s \in \mathcal{S}_b{\setminus}\{0,s_b\}$ and the fully depleted state $s=0$. Constraints \eqref{model:operations_detailed:fleet_depot} and \eqref{model:operations_detailed:fleet_conv_route} ensure that the allocated fleet size is respected. Constraints \eqref{model:operations_detailed:var_depot}--\eqref{model:operations_detailed:var_conv_route} give the domain of the operational variables. 

Throughout the paper, we will use the following compact form of problem \eqref{model:operations_detailed}:
\begingroup
\setspacing
\allowdisplaybreaks
\begin{subequations}
\label{model:operations}
\begin{align}
\label{model:operations:obj} \mathcal{Q}_p(x_p) := &\min_{y_{pr} \in \mathbb{Z}^{m_{pr}}_+, \ \forall r \in \mathcal{R} } \quad  \sum_{r \in \mathcal{R}} c^{y\top}_{pr} y_{pr} \\ \label{model:operations:linking_constraints}  &\text{s.t. } \sum_{r \in \mathcal{R}}A_{r}y_{pr} \leq B\chi_p,\\
\label{model:operations:route_constraints} & \pst D_{r}y_{pr} \leq e_{pr} + E\eta_{pr}, \quad \forall r \in \mathcal{R}.
\end{align}
\end{subequations}
\endgroup
Constraints \eqref{model:operations:linking_constraints} correspond to the charger availability constraints \eqref{model:operations_detailed:depot_capacity}--\eqref{model:operations_detailed:terminal_capacity}, and the route-separable constraints \eqref{model:operations:route_constraints} include the service level, battery dynamics, and fleet size constraints \eqref{model:operations_detailed:service_level}--\eqref{model:operations_detailed:fleet_conv_route}. If solution $x_p$ leads to an infeasible operational problem, then $\mathcal{Q}_p(x_p)= +\infty$.

\subsection{Properties of the problem}\label{sec:model:properties}
We conclude this section by highlighting three key structural properties of the BFEP.

\begin{remark}[NP-hardness]\label{remark1:NP_hardness}
Due to the set covering structure of charger placement, the BFEP is strongly NP-hard even in simple special cases. This can be established by reduction from the set covering problem \citep{garey1979computers}. \begin{revblock}By identifying each terminal with a candidate set covering its connected routes, any instance of the set covering problem can be encoded as a BFEP with one investment period, one time interval, only on-route BEBs, unit service requirements, and charger installations as the only costs.\end{revblock} A similar reduction gives the same result for the BFEP with one investment period and only depot BEBs. We omit \begin{revblock}both\end{revblock} proofs for brevity.
\end{remark}

\begin{remark}[Integer operational problem]\label{remark2:IP_operational_model}
Model \eqref{model:operations_detailed} resembles a minimum-cost circulation problem due to the flow conservation constraints \eqref{model:operations_detailed:depot_flow_full}--\eqref{model:operations_detailed:depot_flow_empty}. Yet, the charger capacity constraints \eqref{model:operations_detailed:depot_capacity} and service level requirements \eqref{model:operations_detailed:service_level} break total unimodularity \citep[see][]{schrijver2003combinatorial}. Instances with nonzero integrality gaps can be constructed to conclude that the problem does not have the integrality property and must be solved as an integer program (IP). \begin{revblock} Importantly, the integrality of the operational variables can significantly impact strategic decisions. In \ref{ec:Proofs}, we exhibit a small single-route instance in which problem \eqref{model:operations_detailed} requires a fleet of at least three depot BEBs for feasibility, whereas its LP relaxation is feasible with only two.\end{revblock}
\end{remark}

\begin{remark}[Monotonicity of the operational value function]\label{remark3:monotonicity}
    By construction of problem \eqref{model:operations}, the strategic variables $x_p = \left(\chi_p,\{\eta_{pr}\}_{r \in \mathcal{R}}\right)$ only appear with nonnegative coefficients on the right-hand side of constraints \eqref{model:operations:linking_constraints}--\eqref{model:operations:route_constraints}. The value function $\mathcal{Q}_p : \mathcal{X}_p \to [0,+\infty]$ is thus nonincreasing, i.e., if $x_p \preceq x'_p $, where $\preceq$ denotes componentwise inequality, then $\mathcal{Q}_p(x_p) \geq \mathcal{Q}_p(x'_p)$. In particular, infeasibility is downward closed: if $x_p \preceq x'_p$ and $\mathcal{Q}_p(x'_p)=+\infty$, then $\mathcal{Q}_p(x_p)=+\infty$.
\end{remark}
\section{Solution methodology}
\label{sec:Algorithms}
The BFEP can be expressed in extensive form by replacing in problem \eqref{model:BFEP} the value function $\mathcal{Q}_p(x_p)$ by its integer programming formulation \eqref{model:operations}. However, the resulting formulation rapidly becomes intractable. Leveraging the structure of the problem is thus essential to achieve efficient solution methods. Section \ref{sec:Algorithms:Preprocessing} presents preprocessing steps to reduce the size of the problem and strengthen its LP relaxation. In Section \ref{sec:Algorithms:Benders}, we present our logic-based Benders decomposition method along with several acceleration techniques. Section \ref{sec:Algorithms:Heuristic} introduces our restriction heuristic algorithm.

\subsection{Preprocessing}
\label{sec:Algorithms:Preprocessing}
This section presents valid inequalities (i.e., inequalities satisfied by any feasible solution) and dual reductions (i.e., inequalities that can remove feasible solutions while guaranteeing that at least one optimal solution remains) for the BFEP. In Section \ref{sec:Algorithms:Preprocessing:dominance}, we exploit the set covering structure of the on-route BEB charging dynamics to devise dominance relations between terminals and eliminate unnecessary locations. Section \ref{sec:Algorithms:Preprocessing:fleet_size} presents valid inequalities on the size of the depot BEB fleet. 

\subsubsection{Dominance relation between terminals}
\label{sec:Algorithms:Preprocessing:dominance}
Since we assume the installation cost of on-route chargers to be linear in the number of chargers and to be the same across terminals, the choice of the best location for fast chargers is driven by route covering. Definition \ref{def:bound_dominated_terminals} establishes a dominance relation between terminals $j \in \mathcal{J}$ based on their connectivity to bus routes and their hosting capacity $\widetilde{\chi}^{\text{UB}}_{j}$. We break ties based on an arbitrary indexing $\mathcal{J} = \{j_1,j_2,\dots,j_m\}$\begin{revblock}, so that a single representative is retained among fully equivalent terminals. Different indexings yield reduced problems that are identical up to a relabeling of such terminals.\end{revblock}
\begin{definition} \label{def:bound_dominated_terminals}
Terminal $j_b$ is dominated by $j_a$ if one of the following conditions holds:
\begin{enumerate}
    \item $\mathcal{R}(j_b) \subset \mathcal{R}(j_a)$,
    \item $\mathcal{R}(j_b) = \mathcal{R}(j_a)$ and $\widetilde{\chi}^{\text{UB}}_{j_b} < \widetilde{\chi}^{\text{UB}}_{j_a}$,
    \item $\mathcal{R}(j_b) = \mathcal{R}(j_a)$ and $\widetilde{\chi}^{\text{UB}}_{j_b} = \widetilde{\chi}^{\text{UB}}_{j_a}$ and $b < a$.
\end{enumerate}
\end{definition}

Proposition \ref{prop1:bound_dominated_terminals} provides a valid dual reduction on the number of on-route chargers installed at dominated terminals. The proof is given in \ref{ec:Proofs}.
\begin{proposition} \label{prop1:bound_dominated_terminals}
    Let $\mathcal{J}(j)$ be the set of terminals dominated by $j\in \mathcal{J}$. Imposing inequality~\eqref{eq:bound_dominated_terminals} for each $j \in \mathcal{J}, \ j' \in \mathcal{J}(j), \ p \in \mathcal{P}$ is a valid dual reduction for the BFEP.
    \begin{equation} \label{eq:bound_dominated_terminals}
        \widetilde{\chi}^{p}_{j'} \leq \max\bigg\{0, \Big\lceil \frac{1}{\rho} \max_{p' \in \mathcal{P}, \ t \in \mathcal{T}} \sum_{r \in \mathcal{R}(j) } d_r^{p't} \Big\rceil - \widetilde{\chi}^{\text{UB}}_{j} \bigg\}
    \end{equation}
\end{proposition}

The upper bounds of Proposition \ref{prop1:bound_dominated_terminals} can be computed outside the optimization process. In addition to strengthening the upper bounds $\widetilde{\chi}^{\text{UB}}_{j}$, this inexpensive preprocessing step generally eliminates several terminals, hence reducing the size of both the strategic and operational models. 

\subsubsection{Valid inequalities on fleet size}\label{sec:Algorithms:Preprocessing:fleet_size}
The service level requirements \eqref{model:operations_detailed:service_level} imply that the number of buses assigned to a route $r \in \mathcal{R}$ in period $p \in \mathcal{P}$ cannot be less than the peak demand $\max_{t \in \mathcal{T}}d^{pt}_r$. For conventional buses and on-route BEBs, which can be in service without interruption, this simple bound on fleet size cannot be tightened in general. In contrast, the LP relaxation of the operational model does not provide a tight lower bound on the number of depot BEBs needed to satisfy the constraints \eqref{model:operations_detailed:service_level}. We thus precompute lower bounds on fleet size for each route. 

Suppose at most $m$ conventional buses and on-route BEBs are assigned to route $r \in \mathcal{R}$ in period $p \in \mathcal{P}$. The minimum number of depot BEBs needed to satisfy the residual service level requirements can be computed by solving the following restricted operational problem:
\begingroup
\setspacing
\allowdisplaybreaks
\begin{subequations}
\label{model:min_fleet}
\begin{align}
\label{model:min_fleet:obj} \bar{\eta}^{\text{LB}}_{pr}(m) := &\min_{\widetilde{\eta}^p_r, 
\widehat{\eta}^p_r \in \mathbb{Z}_+, \ \bar{\eta}^p_{r} \in \mathbb{Z}^{\mathcal{B}}_+, \ y_{pr} \in \mathbb{Z}^{m_{pr}}_+}  \sum_{b \in \mathcal{B}} \bar{\eta}^p_{rb} \\ \label{model:min_fleet:linking_constraints}  &\text{s.t. }  D_{r}y_{pr} \leq e_{pr} + E\eta_{pr},\\
\label{model:min_fleet:route_constraints} & \pst \widetilde{\eta}^p_r + \widehat{\eta}^p_{r} \leq m.
\end{align}
\end{subequations}
\endgroup

Problem \eqref{model:min_fleet} is constructed by projecting problem \eqref{model:operations} onto the space of the operational variables $y_{pr}$ of route $r \in \mathcal{R}$, relaxing the charger capacity constraints \eqref{model:operations:linking_constraints}, and taking the fleet assignment parameters $\eta_{pr}$ as decision variables. Its optimal value $\bar{\eta}^{\text{LB}}_{pr}(m)$ provides a lower bound on the number of depot BEBs that must be assigned to route $r$ in period $p$ to complement a fleet of $\widehat{\eta}_{r}^{p} + \widetilde{\eta}^{p}_r \leq m$ conventional and on-route buses. These lower bounds can be enforced as Big-M constraints. However, to avoid introducing additional binary variables in the model, we instead approximate these Big-M constraints using the piecewise-linear lower envelope of the set of points $\mathcal{M}_{pr} := \{(m, \bar{\eta}^{\text{LB}}_{pr}(m))\}_{m=0}^{\max_{t \in \mathcal{T}}d^{pt}_r}$, i.e., the convex hull of the epigraph of the function $\bar{\eta}^{\text{LB}}_{pr}(m)$ evaluated for a number of buses $m$ ranging from 0 to the peak demand on route $r$ in period $p$. This yields constraints \eqref{ineq:flett_size_lin_constraints}, which we implement as linear inequalities on the fleet size variables $\eta_{pr}$.
\begin{equation} \label{ineq:flett_size_lin_constraints}
    \left(\widetilde{\eta}_r^p + \widehat{\eta}_r^p, \sum_{b \in \mathcal{B}} \bar{\eta}^p_{rb}\right) \in \operatorname{conv}\left(\text{epi}\left( \mathcal{M}_{pr} \right)\right)
\end{equation}

\subsection{Logic-based Benders decomposition algorithm}
\label{sec:Algorithms:Benders}
Our exact algorithm relies on logic-based Benders decomposition (LBBD) \citep{hooker2003logic}. Like classical Benders decomposition, LBBD projects a problem onto the subspace defined by a subset of its decision variables, the \emph{master variables}, and replaces the subproblems by their \emph{value functions} of the master variables. These value functions are represented by auxiliary variables that are bounded by valid inequalities. However, instead of relying on linear programming strong duality to generate cutting planes, LBBD uses problem-specific procedures to separate tight cuts and obtain an exact Benders reformulation. This is necessary in the case of the BFEP, where the subproblems are IPs (Remark \ref{remark2:IP_operational_model}), which leads to nonconvex value functions.

We use the strategic variables $x \in \mathcal{X}$ of model \eqref{model:BFEP} as master variables, and represent the value function $\mathcal{Q}_p : \mathcal{X}_p \to [0,+\infty]$ of the operational problem \eqref{model:operations_detailed} of each period $p \in \mathcal{P}$ by an auxiliary variable $\Theta_p \geq 0$. 
For any point $x'_p \in \mathcal{X}_p$, we define the lower-orthant indicator function of $x_p$ as $\mathbf{1}\{x_p \preceq x'_p\}=1$ if $x_p \preceq x'_p$, and $=0$ otherwise. The logic-based master problem is formulated as:
\begingroup
\setspacing
\allowdisplaybreaks
\begin{subequations}
\label{LBBD_MP}
\begin{align}
    \label{LBBD_MP:obj} & \min_{\substack{x \in \mathcal{X}, \ \Theta \in \mathbb{R}_+^{\mathcal{P}} }} \ f(x) + \sum_{p \in \mathcal{P}} \gamma^p\Theta_p,\\
    \label{LBBD_MP:opt_cut} & \text{s.t. } \Theta_p \geq \mathcal{Q}_p(x'_p)\mathbf{1}\{x_p \preceq x'_p\},  \quad & \forall p \in \mathcal{P}, \ \forall x'_p \in \mathcal{X}_p  : \mathcal{Q}_p(x'_p) < + \infty, \\
    \label{LBBD_MP:feas_cut} & \pst \mathbf{1}\{x_p \preceq x'_p\} = 0, \quad & \forall p \in \mathcal{P}, \ \forall x'_p \in \mathcal{X}_p  : \mathcal{Q}_p(x'_p) = + \infty,
\end{align}
\end{subequations}
\endgroup
where the first-stage objective function $f(x) := \sum_{p \in \mathcal{P}}\gamma^p\left(I_p(x_p - x_{p-1}) + H_p(x_p) \right)$ aggregates the discounted investment and fixed maintenance costs of model \eqref{model:BFEP}. We implement the monotone cuts \eqref{LBBD_MP:opt_cut}--\eqref{LBBD_MP:feas_cut} as disjunctive constraints using shared binary variables (see \ref{ec:sec:Monotone_cuts_implementation}). Proposition \ref{prop2:LBBD} shows that the logic-based master problem \eqref{LBBD_MP} is an equivalent reformulation of the BFEP \eqref{model:BFEP}.
\begin{proposition}\label{prop2:LBBD}
For each $p\in\mathcal{P}$, the set $LBBD_p:=\{(x_p,\Theta_p)\in\mathcal{X}_p\times\mathbb{R}_+ :$ \eqref{LBBD_MP:opt_cut}--\eqref{LBBD_MP:feas_cut} $\}$ is equal to the epigraph $\epi(\mathcal{Q}_p)=\{(x_p, \Theta_p)\in\mathcal{X}_p\times\mathbb{R}_+: \mathcal{Q}_p(x_p)\leq \Theta_p\}$ of the value function $\mathcal{Q}_p$. 
\end{proposition}
\begin{proof}{Proof.} $(\epi(\mathcal{Q}_p)\subseteq LBBD_p):$ Take any $(x_p, \Theta_p)\in\epi(\mathcal{Q}_p)$, so $\Theta_p\ge \mathcal{Q}_p(x_p)$ and $\mathcal{Q}_p(x_p)<+\infty$. Take any $x'_p\in\mathcal{X}_p$ such that $\mathcal{Q}_p(x'_p)<+\infty$. If $x_p \preceq x'_p$, then by Remark~\ref{remark3:monotonicity} we have $\mathcal{Q}_p(x_p)\ge \mathcal{Q}_p(x'_p)$ and hence $\Theta_p\ge \mathcal{Q}_p(x'_p)$, so \eqref{LBBD_MP:opt_cut} holds. Otherwise $x_p\not\preceq x'_p$ and \eqref{LBBD_MP:opt_cut} holds trivially since $\mathbf{1}\{x_p\preceq x'_p\}=0$. Now, take any $x'_p\in\mathcal{X}_p$ such that $\mathcal{Q}_p(x'_p)=+\infty$. If $x_p\preceq x'_p$, then infeasibility being downward closed (Remark~\ref{remark3:monotonicity}) would imply $\mathcal{Q}_p(x_p)=+\infty$, a contradiction. Hence $x_p\not\preceq x'_p$ and \eqref{LBBD_MP:feas_cut} holds. We conclude that $(x_p, \Theta_p)\in LBBD_p$.

$(LBBD_p \subseteq \epi(\mathcal{Q}_p)):$ Take any $(x_p, \Theta_p)\in LBBD_p$. If $\mathcal{Q}_p(x_p)=+\infty$, then \eqref{LBBD_MP:feas_cut} with $x'_p=x_p$ is violated since $\mathbf{1}\{x_p\preceq x_p\}=1$. Thus $\mathcal{Q}_p(x_p)<+\infty$, and \eqref{LBBD_MP:opt_cut} with $x'_p=x_p$ yields $\Theta_p\ge \mathcal{Q}_p(x_p)$. We conclude that $(x_p, \Theta_p)\in \epi(\mathcal{Q}_p)$. \hfill \halmos 
\end{proof}

The LBBD formulation \eqref{LBBD_MP} can be solved by constraint generation. Constraints \eqref{LBBD_MP:opt_cut} and \eqref{LBBD_MP:feas_cut} are initially relaxed. At each iteration $l = 1,2,\dots$, the current relaxed master problem is solved, and its optimal solution $(x^l, \Theta^l) \in \mathcal{X} \times \mathbb{R}_+^{\mathcal{P}}$ provides a lower bound $f(x^l) + \sum_{p \in \mathcal{P}} \gamma^p \Theta^l_p$ on the optimal value of the master problem. For each period $p \in \mathcal{P}$, the operational problem \eqref{model:operations} is then solved for $x_p = x_p^l$. The feasibility cut $\mathbf{1}\{x_p \preceq x^l_p\} = 0$ is added to the relaxed master problem if $\mathcal{Q}_p(x_p^l) = + \infty$, and the optimality cut $\Theta_p \geq \mathcal{Q}_p(x^l_p)\mathbf{1}\{x_p \preceq x^l_p\}$ is otherwise generated. Replacing each variable $\Theta^l$ by $\mathcal{Q}_p(x_p^l)$ in the objective of the relaxed master then yields the upper bound $f(x^l) + \sum_{p \in \mathcal{P}} \gamma^p\mathcal{Q}_p(x^l_p)$. This process is repeated until the optimality gap reaches the desired tolerance. If $\mathcal{X}$ has finite cardinality, this algorithm converges in a finite number of iterations \citep{hooker2003logic}.

In practice, this textbook implementation performs poorly for the BFEP. This is due in part to the fact that minimizing the investment and operational costs are antagonistic objectives (Remark \ref{remark3:monotonicity}). As the initial relaxed master problem neglects operational costs, early iterations are thus characterized by significant underinvestment, producing monotone cuts that are inactive in regions of $\mathcal{X}_p$ with reasonable investments. In the remainder of this section, we develop a range of techniques that aim to provide a strong and tractable approximation of the value functions early in the solving process.

Sections \ref{sec:Algorithms:Benders:Disaggregation} and \ref{sec:Algorithms:Benders:MPS} introduce relaxations of model \eqref{model:operations} that are respectively used to disaggregate the operational costs by route and to strengthen the relaxed master problem formulation. Section \ref{sec:Algorithms:Benders:linear} reviews Benders cut selection techniques and their application to the LP relaxations of our operational problems. Section \ref{sec:Algorithms:Benders:monotone} presents problem-specific monotone feasibility cuts. The outline of our LBBD algorithm is presented in Section \ref{sec:Benders:outline}.

\subsubsection{Disaggregation of the operational problem} 
\label{sec:Algorithms:Benders:Disaggregation}
For a fixed strategic solution $x \in \mathcal{X}$, the BFEP decomposes into a set of $P$ independent operational problems. This is exploited in the master problem \eqref{LBBD_MP}, where a variable $\Theta_p$ bounds the operational costs of each period $p \in \mathcal{P}$. Here, we further decompose the operational problems by route. Since the operational variables $\{y_{pr}\}_{r \in \mathcal{R}}$ of each period are readily partitioned by route, applying the linking constraints \eqref{model:operations:linking_constraints} to individual routes yields a separable relaxation of problem \eqref{model:operations}. The single-route relaxation for route $r \in \mathcal{R}$ is:
\begingroup
\setspacing
\allowdisplaybreaks
\begin{subequations}
\label{model_sr:operations}
\begin{align}
\label{model_sr:operations:obj} \widetilde{\mathcal{Q}}_{pr}(x_p) := &\min_{y_{pr} \in \mathbb{Z}^{m_{pr}}_+ } \quad c^{y\top}_{pr} y_{pr} \\ \label{model_sr:operations:linking_constraints} & \text{s.t. } A_{r}y_{pr} \leq B\chi_p,\\
\label{model_sr:operations:route_constraints} & \pst D_{r}y_{pr} \leq e_{pr} + E\eta_{pr}.
\end{align}
\end{subequations}
\endgroup

By the nonnegativity of the variables $\{y_{pr}\}_{r \in \mathcal{R}}$ and the matrices $\{A_{r}\}_{r \in \mathcal{R}}$, constraints \eqref{model_sr:operations:linking_constraints} are implied by \eqref{model:operations:linking_constraints}. For any feasible solution $\{\bar{y}_{pr}\}_{r \in \mathcal{R}}$ to the operational problem \eqref{model:operations}, $\bar{y}_{pr}$ is thus feasible for \eqref{model_sr:operations} for each route $r\in \mathcal{R}$. Since the objective coefficients of each variable are identical in both formulations, solving the single-route subproblems separately and summing their objectives provides a lower bound on the original operational problem, i.e., $\sum_{r \in \mathcal{R}} \widetilde{\mathcal{Q}}_{pr}(x_p) \leq \mathcal{Q}_p(x_p)$. To exploit the single-route subproblems in our LBBD, we define an auxiliary variable $\theta_{pr}$ for each period $p \in \mathcal{P}$ and each route $r \in \mathcal{R}$, and add the following constraints to the master problem:
\begin{equation}
\label{theta_eq}
    \Theta_p = \sum_{r \in \mathcal{R}}\theta_{pr}, \quad \forall p \in \mathcal{P}.
\end{equation}

The LP relaxation of the single-route subproblems \eqref{model_sr:operations} is used to generate Benders cuts. These single-route cuts have the advantage of being sparse, as they only involve the strategic decisions impacting the operations on a specific route. In practice, they can provide a good approximation of the value functions $\mathcal{Q}_p$ in a few iterations, and significantly improve the performance of our LBBD.

\subsubsection{Master problem strengthening}
\label{sec:Algorithms:Benders:MPS}
A common drawback of Benders decomposition is that the initial relaxed master usually provides a weak relaxation of the problem. This limitation can be mitigated by using an alternative master formulation that includes explicit information from the subproblems. For two-stage stochastic programs, partial decomposition, which consists of retaining a subset of second-stage scenarios in the master problem, can significantly improve the overall performance of Benders decomposition \citep{crainic2021partial, legault2025model}. Alternatively, for two-stage problems with integer subproblems, the initial master problem can be strengthened by retaining all the second-stage variables and relaxing their integrality \citep{gendron2016branch}. However, the resulting master formulation, sometimes described as \emph{semi-relaxed} \citep{liu2024generalized}, may include a prohibitively large number of variables and become impractical to solve. 

Taking inspiration from these techniques, we propose a strengthened master problem formulation that retains a compact LP relaxation of the operational problems. We define, for each period $p \in \mathcal{P}$, and each route $r \in \mathcal{R}$, the collection $\omega^{p}_r := ( \{\bar{\omega}^p_{rb}\}_{b \in \mathcal{B}}, \widetilde{\omega}^p_r, \widehat{\omega}^p_r) = (\{\frac{1}{T}\sum_{t \in \mathcal{T}}\sum_{s \in \mathcal{S}^{w}_b}w^{pt}_{rbs}\}_{b \in \mathcal{B}}, \ \frac{1}{T}\sum_{t \in \mathcal{T}}\sum_{j \in \mathcal{J}(r)}\widetilde{w}^{pt}_{rj}, \ \frac{1}{T}\sum_{t \in \mathcal{T}}\widehat{w}^{pt}_{r} ) \in \mathbb{R}^{|\mathcal{B}| + 2}_+$ of nonnegative continuous variables representing the average service level provided by each type of bus in model \eqref{model:operations_detailed}. From there, averaging the service level constraints \eqref{model:operations_detailed:service_level} of each period and each route over the intervals $t \in \mathcal{T}$ yields a surrogate relaxation that can be expressed using these auxiliary variables:
\begingroup
\setspacing
\allowdisplaybreaks
\begin{align}
    \label{eq:pb:demand_total}
    &\sum_{b \in \mathcal{B}}\bar{\omega}^p_{rb} + \widetilde{\omega}^p_{r} + \widehat{\omega}^p_{r} \geq \frac{1}{T}\sum_{t \in \mathcal{T}}d_{r}^{p t}, \quad & \  \forall r \in \mathcal{R}, \ p \in \mathcal{P}.
\end{align}
\endgroup

For constraints \eqref{eq:pb:demand_total} to serve as nontrivial feasibility cuts in the master problem, the auxiliary variables are then upper bounded by linear functions of the strategic variables. We obtain a first set of bounds by averaging the fleet size constraints \eqref{model:operations_detailed:fleet_depot}--\eqref{model:operations_detailed:fleet_conv_route} over the operational time intervals:
\begingroup
\setspacing
\allowdisplaybreaks
\begin{subequations}\label{UB:omegas:fleet}
\begin{align}
    \label{UB:omega:fleet}
    &\bar{\omega}^p_{rb} \leq \alpha_{rb} \bar{\eta}^p_{rb} , \quad & \forall p \in \mathcal{P}, \ r \in \mathcal{R}, \ b \in \mathcal{B},\\
    \label{UB:omega_tilde:fleet}
    &\widetilde{\omega}^p_{r} \leq \widetilde{\eta}^p_{r} , \quad  & \forall p \in \mathcal{P}, \ r \in \mathcal{R},\\
    \label{UB:omega_hat:fleet}
    &\widehat{\omega}^p_{r} \leq \widehat{\eta}^p_{r} , \quad  & \forall p \in \mathcal{P}, \ r \in \mathcal{R}.
\end{align}
\end{subequations}
\endgroup
where $\alpha_{rb} := \max_{i \in \mathcal{I}, s \in \mathcal{S}^z_b}\frac{s_b-s}{s_b-s+\kappa_{rbis}}$ is an upper bound on the fraction of the time intervals where depot BEBs of type $b \in \mathcal{B}$ can be in service on route $r \in \mathcal{R}$ while respecting the conservation constraints \eqref{model:operations_detailed:depot_flow_full}--\eqref{model:operations_detailed:depot_flow_empty}. This bound is attained by alternating indefinitely between working and recharging from the state $s \in \mathcal{S}^z_b$ that maximizes the ratio of regained battery level $s_b-s$ to charging time $\min_{i \in \mathcal{I}}\kappa_{rbis}$ achieved by charging at the nearest depot.

The average service level can also be bounded by the installed charging infrastructure. To do so, we average the on-route charger capacity constraints \eqref{model:operations_detailed:terminal_capacity} of period $p \in \mathcal{P}$ over the time intervals $t \in \mathcal{T}$. From there, summing over the set of terminals $j \in \mathcal{J}(r)$ connected to a route $r \in \mathcal{R}$, and over all the terminals $j \in \mathcal{J}$, gives inequalities \eqref{bound_omega_tilde1:2} and \eqref{bound_omega_tilde2:2}, respectively. 
\begingroup
\setspacing
\allowdisplaybreaks
\begin{subequations}\label{bound_omega_tilde}
\begin{align}
    \label{bound_omega_tilde1:2} 
    & \widetilde{\omega}_{r}^{p} \leq \sum_{j \in \mathcal{J}(r)} \rho \widetilde{\chi}_{j}^p, & \quad \forall p \in \mathcal{P}, \ r \in \mathcal{R}, \\
    \label{bound_omega_tilde2:2} 
    &\sum_{r \in \mathcal{R}} \widetilde{\omega}_{r}^{p} \leq \sum_{j \in \mathcal{J}} \rho \widetilde{\chi}_{j}^p, & \quad \forall p \in \mathcal{P}.
\end{align}
\end{subequations}
\endgroup

Constraints \eqref{bound_omega_tilde1:2} require that the chargers installed at terminals connected to a route $r\in \mathcal{R}$ suffice to satisfy the average energy requirements of the on-route BEBs assigned to this route, while \eqref{bound_omega_tilde2:2} ensures that the total energy requirements of the fleet can be satisfied over the complete network. 

Similar bounds can be devised for depot BEBs. Let $\beta_{i} := \max_{r \in \mathcal{R}, b \in \mathcal{B}, s \in \mathcal{S}^z_b}\frac{s_b-s}{\kappa_{rbis}}$ denote the maximum number of battery states that can be regained per time interval by charging at depot $i \in \mathcal{I}$. By constraints \eqref{model:operations_detailed:depot_capacity}, which limit to $\bar{\chi}^p_{i}$ the number of vehicles that charge simultaneously at depot $i \in \mathcal{I}$ during period $p \in \mathcal{P}$, it follows that no more than $\sum_{i \in \mathcal{I}} \beta_{i} \bar{\chi}^p_{i}$ units of charge can be regained by the fleet of depot BEBs in each time interval. From there, since the flow conservation constraints \eqref{model:operations_detailed:depot_flow_full}--\eqref{model:operations_detailed:depot_flow_empty} imply that the total flow $\sum_{t \in \mathcal{T}}\sum_{r \in \mathcal{R}}\sum_{b \in \mathcal{B}}\sum_{s \in \mathcal{S}^w_b}w^{pt}_{rbs}$ on service arcs equals the number of units of charge regained by the fleet over the operational horizon, the following inequalities hold:
\begingroup
\setspacing
\allowdisplaybreaks
\begin{align}
\label{bound_omega} 
& \sum_{r \in \mathcal{R}}\sum_{b \in \mathcal{B}}\bar{\omega}^{p}_{rb} \leq  \sum_{i \in \mathcal{I}} \beta_{i} \bar{\chi}^p_{i} , & \quad \forall p \in \mathcal{P}.
\end{align}
\endgroup

Upper bounds on the average service level of on-route BEBs and conventional buses can also be devised from the service level constraints. Indeed, since, for each $(p,r,t)\in \mathcal{P}\times\mathcal{R}\times\mathcal{T}$, the variables $\widetilde{w}_{r}^{pt}$ and $\widehat{w}_{r}^{pt}$ are assumed to have positive objective coefficients in model \eqref{model:operations}, and the service level constraint \eqref{model:operations_detailed:service_level} is the only one that impose a lower bound on their value, the level of service of conventional buses and on-route BEBs never exceeds the demand $d^{pt}_r$ in an optimal solution. From there, given that $m:=\widetilde{\eta}^p_{r}+\widehat{\eta}^p_{r}$ conventional buses and on-route BEBs are assigned to route $r$ in period $p$, their average service level $\widetilde{\omega}^p_r + \widehat{\omega}^p_r$ can be upper bounded by $\sigma_{pr}(m) := \frac{1}{T}\sum_{t \in \mathcal{T}} \min\{d^{pt}_r, m\}$. Analogously to the lower bounds \eqref{ineq:flett_size_lin_constraints}, these upper bounds can be imposed by constructing the piecewise-linear upper envelope of the set of points $\mathcal{M}'_{pr} := \left\{\left(m, \sigma_{pr}(m) \right) \right\}_{m=0}^{\max_{t \in \mathcal{T}}d^{pt}_r}$, giving:
\begin{equation} \label{ineq:conv_hypo}
    \left(\widetilde{\eta}_r^p + \widehat{\eta}_r^p, \widetilde{\omega}^p_r + \widehat{\omega}^p_r\right) \in \operatorname{conv}\left(\text{hypo}\left( \mathcal{M}'_{pr} \right)\right).
\end{equation}

For $(p,r)\in \mathcal{P}\times\mathcal{R}$, let $c^{\widehat{\omega}}_{pr}/T := \min_{t \in \mathcal{T}}c^{\widehat{w}}_{ptr}$ and $c^{\widetilde{\omega}}_{pr}/T := \min_{t \in \mathcal{T}, j \in \mathcal{J}(r)}c^{\widetilde{w}}_{ptrj}$ be the smallest per-interval operating cost of conventional buses and an on-route BEBs. Similarly, we write the smallest per-interval operating cost of depot BEBs of type $b \in \mathcal{B}$ as $c^{\bar{\omega}}_{prb}/T :=  \min_{t \in \mathcal{T}, s \in \mathcal{S}^w_b}c^{w}_{ptrbs} + \min_{t \in \mathcal{T}, i \in \mathcal{I}, s \in \mathcal{S}^z_b} \frac{c^{z}_{ptrbis}}{s_b-s}$, where the second term is the minimal charging costs needed to replenish one unit of charge. The operational costs on route $r$ during period $p$ admit the following lower bound:
\begin{equation}
    \label{ec:model:operations:objective_approximate}
    \theta_{pr} \geq \sum_{b \in \mathcal{B}}c^{\bar{\omega}}_{prb} \bar{\omega}^{p}_{rb} + c^{\widetilde{\omega}}_{pr} \widetilde{\omega}^{p}_{r} + c^{\widehat{\omega}}_{pr} \widehat{\omega}^{p}_{r}.
\end{equation}

We approximate the operational model by including in the master problem the auxiliary variables $\{\omega_r^p\}_{p,r}\in \mathcal{P \times \mathcal{R}}$, along with constraints \eqref{eq:pb:demand_total}--\eqref{ec:model:operations:objective_approximate}. Although this compact relaxation of the operational problem is weaker than its full LP relaxation, it only marginally increases the computational cost of solving the master problem, making it a tractable alternative to the semi-relaxed formulation.

\subsubsection{Benders cut selection}
\label{sec:Algorithms:Benders:linear}

Relying exclusively on the non-convex constraints \eqref{LBBD_MP:opt_cut} and \eqref{LBBD_MP:feas_cut} to approximate the operational value functions is inefficient for multiple reasons: monotone cuts are active on a small region of the first-stage domain, their separation requires solving integer subproblems, and their generation adds binary variables to the master problem. To mitigate these limitations, we generate Benders cuts from the LP relaxations of the subproblems before resorting to monotone cuts. \begin{revblock}Although these relaxations can be weak (Remark~\ref{remark2:IP_operational_model}), they are sufficiently tight in practice for Benders cuts to guide the master iterates toward near-optimal first-stage solutions, leaving to monotone cuts the roles of certifying feasibility and closing residual gaps.\end{revblock}

\begin{revblock}
The selection of strong Benders cuts is a key component of our algorithm. We thus review and computationally evaluate the main selection methods from the literature. Notably, we establish the following equivalence, which has not been previously noted in the literature.\end{revblock}
\begin{proposition}\label{prop3:closest_deepest}
    The closest cut selection method proposed by \cite{seo2022closest} and the Conforti–Wolsey deepest cut selection method proposed by \cite{hosseini2025deepest} are equivalent.
\end{proposition}
\begin{revblock}
To keep the exposition focused, we  defer the review and the proof of Proposition \ref{prop3:closest_deepest} to \ref{ec:benders_theory}, and here only present the closest cuts of \citet{seo2022closest}, which we retain in our final algorithm. 

Let $\pi = (\lambda_p, \{\mu_{pr}\}_{r \in \mathcal{R}})$ denote the dual vector of constraints \eqref{model:operations:linking_constraints}--\eqref{model:operations:route_constraints} in the LP relaxation of the operational problem \eqref{model:operations} of period $p \in \mathcal{P}$, and let $\mathcal{D}_p(x_p; \pi) := \lambda_p^{\top} B\chi_p + \sum_{r \in \mathcal{R}} \mu_{pr}^{\top}\!\left(e_{pr} + E\eta_{pr}\right)$ denote the corresponding dual objective function. Given a master solution $(x'_p, \Theta'_p)$ and a guiding point $(x^o_p, \Theta^o_p)$, where $\Theta^o_p$ is an upper bound on the optimal value of the LP relaxation of model \eqref{model:operations} at $x_p = x^o_p$, so that the guiding point violates no Benders cut, a closest cut is selected by solving:
\begingroup
\setspacing
\allowdisplaybreaks
\begin{subequations}
\label{model:closest_cut}
\begin{align}
\label{model:closest_cut:obj} \max_{(\pi, \pi_0) \leq 0} \quad & \pi_0\Theta'_p + \mathcal{D}_p(x'_p; \pi) \\
\label{model:closest_cut:dual} \text{s.t. } \quad & \pi_0 c^{y}_{pr} + A_{r}^{\top} \lambda_p + D_{r}^{\top} \mu_{pr} \leq 0, \quad \forall r \in \mathcal{R}, \\
\label{model:closest_cut:norm} & \pi_0(\Theta'_p-\Theta^o_p) + \mathcal{D}_p(x'_p; \pi) -  \mathcal{D}_p(x^o_p; \pi) = 1.
\end{align}
\end{subequations}
\endgroup
The objective \eqref{model:closest_cut:obj} maximizes the violation of the selected cut at the master solution. Constraints \eqref{model:closest_cut:dual} restrict the multipliers to those that yield Benders cuts valid for the master problem \eqref{LBBD_MP}. Constraint \eqref{model:closest_cut:norm} normalizes the dual vectors relative to the guiding point, so that the hyperplane of the selected cut crosses the segment joining $(x^o_p, \Theta^o_p)$ and $(x'_p, \Theta'_p)$ at the point closest to the guiding point. The optimal solution $(\bar{\pi}, \bar{\pi}_0)$ of problem \eqref{model:closest_cut} yields the Benders cut:
\begin{equation}\label{eq:benders_cut}
    \bar{\pi}_0 \Theta_p + \mathcal{D}_p(x_p; \bar{\pi}) \leq 0,
\end{equation}
which can be expressed as the feasibility cut $\mathcal{D}_p(x_p; \bar{\pi}) \leq 0$ if $\bar{\pi}_0 =0$, and as the optimality cut $\Theta_p \geq \mathcal{D}_p(x_p; -\bar{\pi}/\bar{\pi}_0)$ otherwise.\end{revblock}

\paragraph{Implementation.} 
We generate the Benders cuts \eqref{eq:benders_cut} from the operational problems \eqref{model:operations} and their single-route relaxations \eqref{model_sr:operations} by solving the closest cut selection problem \eqref{model:closest_cut}. The single-route cuts of period $p \in \mathcal{P}$ for route $r \in \mathcal{R}$ are obtained by replacing the set of routes $\mathcal{R}$ by the singleton $\{r\}$ and the auxiliary variable $\Theta_p$ by $\theta_{pr}$. \begin{revblock}An empirical comparison of the Benders cut selection methods reviewed in \ref{ec:benders_theory} is presented in \ref{ec:Benders_cuts}. \end{revblock}

\subsubsection{Monotone cuts}
\label{sec:Algorithms:Benders:monotone}
If a master problem iteration returns a solution $(x', \Theta') \in \mathcal{X} \times \mathbb{R}^{\mathcal{P}}_+$ for which $(x'_p, \Theta'_p)$ violates a Benders cut \eqref{eq:benders_cut} in at least one period $p \in \mathcal{P}$, one can directly proceed to the next iteration without considering integer subproblems. Monotone cuts are only needed when classical cuts do not suffice to eliminate a solution from the feasible domain of the restricted master problem. In the integer phase of our algorithm, we solve for each period $p \in \mathcal{P}$ the integer subproblem \eqref{model:operations} associated with solution $x'_p$. If this problem is feasible and its optimal value is strictly underestimated at the current master solution, i.e., $\Theta'_p < \mathcal{Q}_p(x'_p) < \infty$, then we generate the optimality cut $\Theta_p \geq \mathcal{Q}_p(x'_p)\mathbf{1}\{x_p \preceq x'_p\}$. However, if the integer problem is infeasible, instead of directly generating a generic feasibility cut $\mathbf{1}\{x_p \preceq x'_p\} = 0$, we execute a sequence of tests aimed at identifying a subset of components of $x'_p$ that cause the current solution to be infeasible. 

\paragraph{Single-route fleet cuts.} For each route $r \in \mathcal{R}$ and period $p \in \mathcal{P}$, we verify whether the fleet assignments $\eta'_{pr}$ imply the infeasibility of the operational problem by solving the single-route feasibility problem $y_{pr} \in \mathbb{Z}^{m_{pr}}_+, D_ry_{pr}\leq e_{pr} + E\eta'_{pr}$, which relaxes the charger capacity constraints \eqref{model:operations:linking_constraints}. If the feasible set $\mathcal{Y}_{pr}(\eta'_{pr}):= \{y_{pr} \in \mathbb{Z}^{m_{pr}}_+: D_ry_{pr}\leq e_{pr} + E\eta'_{pr}\}$ is empty, then we conclude that the fleet of depot BEBs is insufficient to satisfy the service level requirements unless the fleet size of conventional buses and on-route BEBs is increased. We obtain the following feasibility cut:
\begin{equation}
\label{feas_cut_sr}
    \mathbf{1}\{(\bar{\eta}_{r}^p, \widetilde{\eta}^p_r{+}\widehat{\eta}^p_r) \preceq (\bar{\eta}'^p_{r}, \widetilde{\eta}'^p_r{+}\widehat{\eta}'^p_r)\} =0.
\end{equation}

\paragraph{Aggregated depot cuts.} If the single-route fleet assignment tests are inconclusive, we solve a multi-route feasibility problem that considers the on-route charging infrastructure and the total number of installed depot chargers. We replace the original set $\mathcal{I}$ of depots by a unique aggregate depot $i^*$. For each tuple $r \in \mathcal{R}$, $b \in \mathcal{B}$, $s \in \mathcal{S}^z_b$, we set the charging time parameter of the aggregate depot to $\kappa_{rbi^*s} = \min_{i \in \mathcal{I}} \kappa_{rbis}$. This relaxes the depot charger capacity constraints by considering optimistic charging times. The feasibility test consists of solving the operational problem of period $p$ for the current fleet $\eta'_p$ and on-route chargers $\widetilde{\chi}'^p$, with the singleton $\{i^*\}$ replacing the original set $\mathcal{I}$, and the objective of minimizing the number of depot chargers, which is treated as a variable. If $\bar{\chi}'^{\text{LB}}_p := \min\{\bar{\chi}_{i^*}^p \in \mathbb{Z}_+ :  $\eqref{model:operations:linking_constraints}--\eqref{model:operations:route_constraints}$, y_{pr} \in \mathbb{Z}^{m_{pr}}_+ \ \forall r \in \mathcal{R}\}$ exceeds the number of depot chargers $\sum_{i \in \mathcal{I}}\bar{\chi}'^p_{i}$ from the current solution, we conclude that the total number of depot chargers or at least one component of the fleet sizing or on-route chargers installation decision vectors must be increased to recover feasibility. This results in the following feasibility cut:
\begin{equation}
\label{feas_cut_mr}
    \mathbf{1}\bigg\{\Big(\eta_p, \widetilde{\chi}^p, \sum_{i \in \mathcal{I}}\bar{\chi}^p_{i}\Big) \preceq \Big(\eta'_p, \widetilde{\chi}'^p, \bar{\chi}'^{\text{LB}}_p{-}1\Big)\bigg\} =0.
\end{equation}

\paragraph{Generic monotone feasibility cuts.} Finally, if these feasibility tests fail to eliminate $x'_p$ from the feasible domain of the relaxed master problem, we resort to a generic monotone cut, which states that at least one component of the strategic solution $x_p$ must be strictly higher than in solution $x'_p$:
\begin{equation}
\label{feas_cut_gen}
    \mathbf{1}\{x_p \preceq x'_p\} = 0.
\end{equation}

It is immediate that $\eqref{feas_cut_sr} \implies \eqref{feas_cut_mr} \implies \eqref{feas_cut_gen}$, with the first implication holding for any route $r \in \mathcal{R}$. In addition to eliminating a larger subset of the master problem's domain, the sparser feasibility cuts require fewer componentwise comparisons in their indicator function. 

\subsubsection{Outline of the algorithm} \label{sec:Benders:outline}
In Algorithm \ref{alg:Benders}, we describe our LBBD algorithm, which is implemented in a classical Benders decomposition fashion. The relaxed master problem, equipped with the disaggregated auxiliary variables of Section \ref{sec:Algorithms:Benders:Disaggregation} and the approximate operational model of Section \ref{sec:Algorithms:Benders:MPS}, is initialized in step \ref{alg:Benders:init_MP}. Each iteration of the algorithm then comprises three phases: (i) solving the current master problem; (ii) separating Benders cuts from linear subproblems; and (iii) separating monotone cuts from integer subproblems.

In the first phase (steps \ref{alg:Benders:solve_MP}--\ref{alg:Benders:update_LB}), the current master problem is solved, and the global lower bound is updated. In the second phase (steps \ref{alg:Benders:loop_linear_cuts}--\ref{alg:Benders:loop_linear_cuts:end}), violated multi-route (steps \ref{alg:Benders:closest_cut_problem}--\ref{alg:Benders:closest_cut_generation}) and single-route (step \ref{alg:Benders:single_route}) Benders cuts are separated using the closest cut selection method from Section \ref{sec:Algorithms:Benders:linear}. In the third phase (steps \ref{alg:Benders:enter_IP_phase}--\ref{alg:Benders:enter_IP_phase:end}), the integer operational problems are solved (step \ref{alg:Benders:solve_IP_subproblem}). Monotone optimality (step \ref{alg:Benders:integer_opt_cut}) and feasibility (step \ref{alg:Benders:integer_feas_cut}) cuts are generated as described in Section \ref{sec:Algorithms:Benders:monotone}, and the incumbent solution is updated (step \ref{alg:Benders:update_incumbent}). These steps are repeated until the relative tolerance rel\_tol (e.g. 0.01\%) is reached. As long as the upper bound UB$^R$ on the current relaxed master problem and the global upper bound UB respect (UB-UB$^R$)/UB $>$ rel\_tol, which means that the incumbent cannot be shown to attain the desired optimality gap from the current relaxed master, we exit step \ref{alg:Benders:solve_MP} early when more than two seconds have been spent on the current master problem and the lower bound LB$^R$ on the current master problem exceeds the global lower bound LB.

\begin{algorithm}[t]
\caption{Logic-based Benders decomposition algorithm}
\label{alg:Benders}
\begin{algorithmic}[1]
\State \label{alg:Benders:init_bounds} Initialize $\text{LB} \gets 0 ;\quad \text{UB} \gets +\infty ; $
\State \label{alg:Benders:init_MP} Initialize the relaxed master problem MP, \begin{revblock}with first-stage variables $x \in \mathcal{X}$, auxiliary variables $\Theta{\in}\mathbb{R}^{\mathcal{P}}_+, \theta{\in}\mathbb{R}^{\mathcal{P} \times \mathcal{R}}_+$ and $\omega{\in}\mathbb{R}^{\mathcal{P} \times \mathcal{R} \times { (|\mathcal{B}|+2) } }_+$, objective $\min f(x){+}\hspace{-0.1cm}\sum\limits_{p \in \mathcal{P}} \gamma^p\Theta_p$, and constraints \eqref{theta_eq}--\eqref{ec:model:operations:objective_approximate}.\end{revblock}
\While{$(\text{UB}-\text{LB})/\text{UB} > \text{rel\_tol} $}
  \State \label{alg:Benders:solve_MP} Solve MP, with bounds LB$^R$ and UB$^R$, and save the components $(x', \Theta', \theta')$ of the solution 
  \State \label{alg:Benders:update_LB} Update $\text{LB} \gets$ LB$^R$
  \For{$p \in \mathcal{P}$} \label{alg:Benders:loop_linear_cuts}
    \State \label{alg:Benders:closest_cut_problem} Solve the cut separation problem \eqref{model:closest_cut} and get the solution $(\bar{\pi}, \bar{\pi}_0)$
    \If{$(x'_p, \Theta_p')$ violates the Benders cut \eqref{eq:benders_cut} associated with $(\bar{\pi}, \bar{\pi}_0)$}
    \State \label{alg:Benders:closest_cut_generation} Generate the Benders cut \eqref{eq:benders_cut} 
    \For{$r \in \mathcal{R}$}
        \State \label{alg:Benders:single_route} Repeat steps \ref{alg:Benders:closest_cut_problem} to \ref{alg:Benders:closest_cut_generation}, with $\{r\}$ replacing $\mathcal{R}$ and $\theta_{pr}$ replacing $\Theta_p$
    \EndFor
  \EndIf
  \EndFor \label{alg:Benders:loop_linear_cuts:end}
  \If{the LP relaxation of \eqref{model:operations} at $x'_p$ is feasible for each $p \in \mathcal{P}$} \label{alg:Benders:enter_IP_phase}
  \For{$p \in \mathcal{P}$}
  \State Solve the operational problem \eqref{model:operations}, with optimal value $\mathcal{Q}_p(x'_p)$ \label{alg:Benders:solve_IP_subproblem}
  \If{$\Theta'_p < \mathcal{Q}_p(x'_p) < \infty$}
  \State Generate the optimality cut $\Theta_p \geq \mathcal{Q}_p(x'_p)\mathbf{1}\{x_p \preceq x'_p\}$  \label{alg:Benders:integer_opt_cut}
  \ElsIf{$\mathcal{Q}_p(x'_p) = \infty$}
  \State For each $r \in \mathcal{R}$ such that $\mathcal{Y}_{pr}(\eta'_{pr})=\emptyset$, generate the violated cut \eqref{feas_cut_sr}. If none is \phantom{.......................} violated, generate \eqref{feas_cut_mr} if $\bar{\chi}^{\text{LB}}_p(\widetilde{\chi}'_p, \eta'_p) > \sum_{i \in \mathcal{I}}\bar{\chi}'^{p}_{i}$, and \eqref{feas_cut_gen} otherwise. \label{alg:Benders:integer_feas_cut}
  \EndIf
  \EndFor
  \If{$f(x') + \sum_{p \in \mathcal{P}} \gamma^p\mathcal{Q}_p(x'_p) < \text{UB}$}
  \State  Update $\text{UB} \gets f(x') + \sum_{p \in \mathcal{P}} \gamma^p\mathcal{Q}_p(x'_p)$;\ Set $x'$ as the incumbent solution\label{alg:Benders:update_incumbent}
  \EndIf
  \EndIf \label{alg:Benders:enter_IP_phase:end}
\EndWhile\\
\Return the incumbent solution, with optimal value UB
\end{algorithmic}
\end{algorithm}

\subsection{Arc selection algorithm}
\label{sec:Algorithms:Heuristic}
Exact and heuristic algorithms serve complementary roles in this work. Our LBBD algorithm is designed to expand the frontier of instances for which provably optimal solutions can be obtained. It achieves optimality for instances of practical interest to transit authorities, who typically plan electrification in phases by considering a subset of routes at a time. However, as the planning scope broadens to long-term, citywide electrification, even finding a feasible solution becomes challenging, and heuristic methods become essential. 

We propose an easily implementable algorithm operating directly on the extensive formulation of the BFEP. This method, which we call the arc selection algorithm, exploits our modeling of BEB schedules as circulations on dense graphs. The goal is to sparsify these graphs in a systematic way to build compact restrictions of the problem. We identify a restricted set of relevant arcs by considering the active variables in the optimal solutions of a collection of single-route problems and LP relaxations of the BFEP. The algorithm then solves a sequence of restricted problems where these arcs are progressively introduced. Finally, the original extensive formulation is solved, which provides a global optimality gap and allows one to use the arc selection algorithm as an exact method. In each step, the incumbent solution is used as a warm start.

\begin{revblock}
The arc selection algorithm is closely related to the Kernel Search heuristic \citep{angelelli2010kernel, guastaroba2017adaptive}, which solves a sequence of restricted MIPs built around a \emph{kernel} of variables presumed to be nonzero in an optimal solution. The kernel is initialized with the active variables of an LP relaxation, and the remaining variables are introduced in groups ordered by their reduced costs. Our algorithm follows the same restriction principle, but exploits the structure of the BFEP: only depot BEB scheduling variables are restricted in the intermediate formulations, and they are reintroduced based on their activity in both LP and IP auxiliary problems.
\end{revblock}

\begin{revblock}The arc selection algorithm, detailed in\end{revblock} Algorithm \ref{alg:arc_selection}, comprises four phases. In the first phase (steps \ref{alg:arc_selection:LP1}--\ref{alg:arc_selection:collect_arcs}), three sets of arcs associated with the flow variables $\{(w^p_r, v^p_r, z^p_r)\}_{(p,r) \in \mathcal{P}\times \mathcal{R}}$ of the BEB scheduling model are constructed. The first two sets $\mathcal{E}^{\text{LP1}}$ and $\mathcal{E}^{\text{LP2}}$ correspond to the active arcs from the LP relaxation of the BFEP, solved without \begin{revblock}and then\end{revblock} with additional constraints on early charging. \begin{revblock}Early charging designates charging operations initiated before the full depletion of the battery, i.e., from a state $s>0$. Charging from the fully depleted state tends to be cost-efficient, as the deadhead costs of the charging trip are amortized over the largest number of service hours. We therefore want the restricted problems to include relevant arcs of this type. Solving the LP relaxation with $\mathcal{S}^z_b=\{0\}$, so that only charging trips from the fully depleted state are allowed, is a simple way to identify such arcs.\end{revblock} To construct the third set $\mathcal{E}^{\text{SR}}$, we solve, for each period $p \in \mathcal{P}$, each route $r \in \mathcal{R}$, each depot BEB type $b \in \mathcal{B}$, and each size $m \in \{0,1,\dots,\max_{t \in \mathcal{T}}d_r^{pt}-1 \}$ of the combined fleet of on-route BEBs and conventional buses, a single-route minimum-cost scheduling problem in which the charger capacity constraints \eqref{model:operations:linking_constraints} are relaxed, the depot BEBs are restricted to type $b$, and their number is fixed to the minimum feasible fleet size $\bar{\eta}^{\text{LB}}_{pr}(m)$. Each arc that is active in at least one optimal solution is included in $\mathcal{E}^{\text{SR}}$. In the second (step \ref{alg:arc_selection:IP1}) and third (step \ref{alg:arc_selection:IP2}) phases, we solve the extensive formulation of the BFEP associated with the edge set $\mathcal{E}^{\text{LP1}} \cup \mathcal{E}^{\text{LP2}}$, and then with the edge set $\mathcal{E}^{\text{LP1}} \cup \mathcal{E}^{\text{LP2}} \cup \mathcal{E}^{\text{SR}}$. Finally, in the fourth phase (step \ref{alg:arc_selection:EX}), we solve the warm-started extensive formulation without restrictions. In our implementation, 75\% of the time budget allocated to the heuristic steps \ref{alg:arc_selection:LP1}--\ref{alg:arc_selection:IP2} goes to the second restricted IP (step \ref{alg:arc_selection:IP2}).

\begin{algorithm}[H]
\caption{Arc selection algorithm}
\label{alg:arc_selection}
\begin{algorithmic}[1]
\State \label{alg:arc_selection:LP1} Solve the LP relaxation of the extensive formulation, let $\mathcal{E}^{\text{LP1}}$ be the set of active arcs $(w,v,z)$
\State \label{alg:arc_selection:LP2} Solve the LP relaxation of the extensive formulation without early charging, i.e., with $\mathcal{S}^z_b = \{0\}$ $\forall b \in \mathcal{B}$, let $\mathcal{E}^{\text{LP2}}$ be the set of active arcs $(w,v,z)$
\For{$(p,r,b) \in \mathcal{P} \times \mathcal{R} \times \mathcal{B}$}
    \For{$m \in \{0,1,\dots,\max_{t \in \mathcal{T}}d_r^{pt}-1 \}$}
    \State With only depot BEBs of type $b$, i.e., $\mathcal{B} \gets \{b\}$, solve \eqref{model:min_fleet}, with optimal value $\bar{\eta}^{\text{LB}}_{pr}(m)$
    \State Solve \begin{revblock}\eqref{model:min_fleet} with the objective function $c^{y\top}_{pr} y_{pr}$ and the additional constraint $\bar{\eta}^p_{rb} = \bar{\eta}^{\text{LB}}_{pr}(m)$\end{revblock}
    \State Let $\mathcal{E}^{\text{SR}}_{prbm}$ be the set of active arcs $(w^p_{rb},v^p_{rb},z^p_{rb})$ in the optimal solution
    \EndFor
    \State Collect the active arcs for bus type $b$ on route $r$ in period $p$ as $\mathcal{E}^{\text{SR}}_{prb}:=\bigcup_{m = 0}^{\max_{t \in \mathcal{T}}d_r^{pt}-1} \mathcal{E}^{\text{SR}}_{prbm}$
\EndFor
\State \label{alg:arc_selection:collect_arcs} Collect the active arcs from the single-route problems as $\mathcal{E}^{\text{SR}}:=\bigcup_{(p,r,b) \in \mathcal{P} \times \mathcal{R} \times \mathcal{B}} \mathcal{E}^{\text{SR}}_{prb}$
\State \label{alg:arc_selection:IP1} Solve the extensive formulation with positive flows allowed on the arcs $\mathcal{E}^{\text{LP1}} \cup \mathcal{E}^{\text{LP2}}$. Save the incumbent solution $(x^1, y^1)$ when the time limit is reached 
\State \label{alg:arc_selection:IP2} Solve the extensive formulation with positive flows allowed on the arcs $\mathcal{E}^{\text{LP1}} \cup \mathcal{E}^{\text{LP2}} \cup \mathcal{E}^{\text{SR}}$, warm-started at $(x^1, y^1)$. Save the incumbent solution $(x^2, y^2)$ when the time limit is reached
\State \label{alg:arc_selection:EX} Solve the extensive formulation, warm-started at $(x^2, y^2)$. Save the incumbent solution $(x^3, y^3)$ when the time limit is reached, along with the bounds LB and UB on the optimal value\\
\Return the incumbent solution $(x^3, y^3)$, and the bounds LB and UB
\end{algorithmic}
\end{algorithm}

\section{Computational experiments}
\label{sec:Experiments}
This computational study evaluates: (i) the impact of our acceleration techniques on the performance of the LBBD framework; (ii) the performance of our algorithms as exact methods; and (iii) the scalability of our heuristic algorithm on citywide electrification planning instances. Section \ref{sec:Experiments:aceclerationLBBD} studies the acceleration techniques. Sections \ref{sec:Experiments:small} and \ref{sec:Experiments:large} present results for the partial and complete bus networks. Finally, in Section \ref{sec:Experiment:chicago}, we illustrate our model on the Chicago transit system. The experiments were conducted with a Python implementation and Gurobi 10.0.3 on a computing cluster node equipped with 48 Intel Xeon Platinum 8260 cores @ 2.40 GHz and 192 GB of RAM. The relative optimality tolerance of the algorithms was set to $0.01\%$.

Our instances are constructed from historical data made available by transit agencies on the \href{https://mobilitydatabase.org/}{Mobility Database}. We selected eight US cities, for which we built instances based on the real routes ($\mathcal{R}$), terminals ($\mathcal{J}$), and depots ($\mathcal{I}$) in use. We defined the minimum service requirement parameters on each route as the average number of buses in service during each hour ($T{=}24$) of a typical weekday, and assumed that this level would remain constant throughout the planning horizon ($d^{p}{=}d^{p'} \ \forall p,p' \in \mathcal{P}$). The model parameters used for our tests are detailed in \ref{ec:Experiments:instances}.

\subsection{Logic-based Benders decomposition accelerations}
\label{sec:Experiments:aceclerationLBBD}
We study the impact of each acceleration technique on the overall performance of our LBBD framework. Specifically, we evaluate the individual effect of replacing the closest cut selection method (CC) of Section \ref{sec:Algorithms:Benders:linear} by standard Benders cuts, and of removing the preprocessing (PP) of Section \ref{sec:Algorithms:Preprocessing}, the disaggregation method (DA) of Section \ref{sec:Algorithms:Benders:Disaggregation}, the master problem strengthening (MS) of Section \ref{sec:Algorithms:Benders:MPS}, and the improved monotone cuts (MC) of Section \ref{sec:Algorithms:Benders:monotone}. \begin{revblock}Finally, we consider the baseline LBBD implementation that uses none of these accelerations.\end{revblock} 

For each problem size $(|\mathcal{R}|, |\mathcal{P}|) \in \{(3,4),(6,6),(9,6)\}$, we consider a set of 10 instances composed of disjoint subsets of routes from the Atlanta network. We impose a time limit of four hours per instance for each method, with the default optimality tolerance of $0.01\%$. Table \ref{tab:accelerations_comparisons_gm} presents average results for the computing time, in seconds, the number of iterations (Iter), the number of Benders cuts (BCuts) and monotone cuts (MCuts) added to the model, and the number of binary indicators (Ind) needed to encode the latter. These metrics are reported as geometric means to better reflect central tendencies. We present the percentage of computing time spent on the master problem (MP), the linear cut separation problems (LP), and the integer monotone cut separation problems (IP). Finally, the number of instances solved to optimality, and the average optimality gap, taken as 0 for instances that terminated before the time limit, are reported. Instances that reached the time limit contribute to the overall results using the figures recorded at termination.

The results of Table \ref{tab:accelerations_comparisons_gm} show that our acceleration techniques drastically improve the performance of the LBBD framework. Indeed, they allow us to solve 29 instances to optimality within the time limit, compared to 11 for the baseline method. Moreover, for challenging instances, they reduce the average computing time and optimality gap by three orders of magnitude. 

\vspace{-0.5cm}
\begin{table}[H]
\centering
\caption{Performance of LBBD with different acceleration techniques}
\resizebox{1\textwidth}{!}{  
\setlength{\tabcolsep}{7pt}
\renewcommand{\arraystretch}{0.75}
\begin{tabular}{clrrrrrrrrrr}
\toprule
\multirow{2}{*}{\specialcell{Instances \\ $(|\mathcal{R}|, |\mathcal{P}|)$}} & \multirow{2}{*}{\specialcell{Accelerations}}  & \multicolumn{4}{c}{Summary}  & \multicolumn{3}{c}{Cuts} & \multicolumn{3}{c}{Time (\%)} \\
\cmidrule(lr){3-6} \cmidrule(lr){7-9} \cmidrule(lr){10-12}
&  &  \multicolumn{1}{c}{Opt} & \multicolumn{1}{c}{Gap} & \multicolumn{1}{r}{Time} & \multicolumn{1}{r}{Iter} & \multicolumn{1}{c}{BCuts} & \multicolumn{1}{c}{MCuts} & \multicolumn{1}{r}{Ind} & MP& LP & IP \\
\midrule 
      & All           & \textbf{10} &   \textbf{0.0000} &       11.1 &     \textbf{5.9} &     33.1 &    \textbf{1.7} &     \textbf{3.1} &     10.5 &    44.7 &    10.2 \\     
      & All $-$ PP       & \textbf{10} &   \textbf{0.0000} &       17.2 &    10.3 &     57.3 &    2.1 &     3.4 &    22.8 &    52.3 &    10.6 \\     
      & All $-$ DA       & \textbf{10} &   \textbf{0.0000} &       11.8 &     8.4 &     \textbf{19.2} &    1.9 &     3.3 &    15.8 &    37.1 &    12.2 \\     
(3,4) & All $-$ MS       & \textbf{10} &   \textbf{0.0000} &       18.6 &    11.4 &     90.4 &    1.9 &     3.3 &    7.8 &    65.9 &    8.0 \\     
      & All $-$ CC       & \textbf{10} &   \textbf{0.0000} &       12.1 &     9.4 &     58.2 &    2.1 &     4.5 &    29.7 &    19.6 &    21.0 \\     
      & All $-$ MC       & \textbf{10} &   \textbf{0.0000} &       \textbf{10.9} &     6.1 &     33.5 &    1.9 &     6.2 &    11.5 &    45.0 &    7.9\\     
      & None          & \textbf{10} &   \textbf{0.0000} &      174.3 &    65.0 &    219.3 &    2.4 &    27.2 &    58.3 &    24.8 &    15.7 \\     
\midrule       
      & All           & \textbf{10} &   \textbf{0.0000} &       \textbf{56.9} &    \textbf{10.3} &    113.3 &    1.9 &     4.9 &    29.7 &    38.3 &    13.3 \\        
      & All $-$ PP       &  9 &   0.0010 &      133.0 &    22.5 &    164.6 &    2.6 &     5.6 &    39.0 &    35.8 &    16.3 \\     
      & All $-$ DA       & \textbf{10} &   \textbf{0.0000} &      146.8 &    26.8 &    \textbf{101.8} &    2.2 &     6.6 &    57.6 &    28.3 &    5.6 \\     
(6,6) & All $-$ MS       & \textbf{10} &   \textbf{0.0000} &       69.8 &    14.2 &    241.5 &    \textbf{1.6} &     \textbf{3.2} &    19.9 &    58.2 &    8.7 \\      
      & All $-$ CC       & \textbf{10} &   \textbf{0.0000} &      108.1 &    21.0 &    230.1 &    2.4 &     7.2 &    63.1 &    13.5 &    14.6 \\     
      & All $-$ MC       & \textbf{10} &   \textbf{0.0000} &       57.9 &    10.8 &    114.6 &    2.1 &    17.6 &    31.7 &    38.3 &    11.3 \\     
      & None          &  1 &   0.7932 &    13747.4 &   209.3 &   1141.6 &    2.8 &    90.2 &    95.7 &    2.6 &    1.5 \\    
\midrule      
      & All           & \textbf{9} &    \textbf{0.0037} &      \textbf{372.8} &    \textbf{24.0} &    \textbf{247.9} &    3.5 &    13.1 &    54.5 &    29.2 &    7.7 \\       
      & All $-$ PP       & 8 &    0.0095 &      659.4 &    33.2 &    317.7 &    5.9 &    \textbf{10.4} &    62.4 &    26.2 &    7.5 \\     
      & All $-$ DA       & 2 &    0.0498 &    10306.1 &    80.3 &    308.2 &    4.1 &    27.5 &    95.7 &    2.9 &    0.9 \\  
(9,6) & All $-$ MS       & \textbf{9} &    0.0038 &      444.9 &    26.8 &    422.0 &    3.9 &    13.0 &    45.1 &    39.4 &    8.3 \\   
      & All $-$ CC       & 6 &    0.0334 &     2137.6 &    61.3 &    581.8 &    5.5 &    35.5 &    82.1 &    6.3 &    9.8 \\
      & All $-$ MC       & 8 &    0.0084 &      398.8 &    24.9 &    251.8 &    3.6 &    36.4 &    53.2 &    29.2 &    8.8 \\     
      & None          & 0 &    2.6173 &    14400.0 &   198.9 &   1190.7 &    \textbf{1.1} &    13.4 &    91.9 &    4.8 &    3.2 \\     
\bottomrule
\end{tabular}
}  
\label{tab:accelerations_comparisons_gm}
\end{table}
\vspace{-0.5cm}

Our ablation study reveals that our acceleration techniques act in complementarity and drastically improve the overall performance of the LBBD algorithm. 
\begin{itemize}[wide=0pt]
    \item (PP): The preprocessing steps divide the average solving time and the number of iterations by a factor of two across all the groups of instances. 
    \item (DA): Disaggregating the operational model based on single-route relaxations is our most impactful acceleration technique. It slightly increases the number of generated Benders cuts for smaller instances. However, in the last group, excluding this technique reduces from nine to two the number of solved instances and increases the average computing time by two orders of magnitude.
    \item (MS): Including a relaxation of the operational model in the master problem formulation is our most efficient acceleration technique for the smallest instances, where it divides by two the number of iterations and by three the number of generated Benders cuts. The master problem strengthening consistently remains beneficial, but is less impactful for large instances. 
    \item (CC): Using the closest cuts systematically reduces the number of generated cuts. For the last group, it divides the solving time and the number of iterations by almost six and three, respectively. 
    \item (MC): These accelerations have a limited impact for these instances, where the LP relaxation of the operational problems is generally tight and the algorithm can rely primarily on classical Benders cuts. Yet, in the last group, one instance is solved by the variants All, All$-$DA, and All$-$MS, which all generate two single-route fleet assignment cuts \eqref{feas_cut_sr}. This instance was not solved by the variant All$-$MC, which had instead generated five weaker generic feasibility cuts \eqref{feas_cut_gen} at termination.
\end{itemize}

\subsection{Exact methods on small instances}
\label{sec:Experiments:small}
In this section, we evaluate the ability of our LBBD algorithm to solve small instances to optimality. We compare LBBD to the extensive formulation (EX) and arc selection (AS) approaches. Our benchmark set is generated from the networks of Atlanta, Boston, Chicago, and Dallas by considering the electrification of a subset of routes over $|\mathcal{P}| \in \{2,4,6,8,10\}$ years, with and without on-route chargers and BEBs. We take $|\mathcal{R}| \in \{3,6,9,12,15\}$ routes when on-route charging is allowed, and $|\mathcal{R}| \in \{2,4,6,8,10\}$ otherwise, for a total of 200 instances. Due to the larger number of instances in this benchmark set compared to that of the previous section, the time limit is reduced to two hours per instance. \begin{revblock}For AS, the computing budget is divided equally between the heuristic phases (steps \ref{alg:arc_selection:LP1}--\ref{alg:arc_selection:IP2}) and the exact phase (step \ref{alg:arc_selection:EX})\end{revblock}. The preprocessing of Section \ref{sec:Algorithms:Preprocessing} is applied to all methods. 

\begin{figure}[H]
\caption{Solved instances and optimality gaps for partial networks}
    \centering
    \begin{subfigure}[b]{0.49\textwidth}
        \centering
        \includegraphics[width=\textwidth]{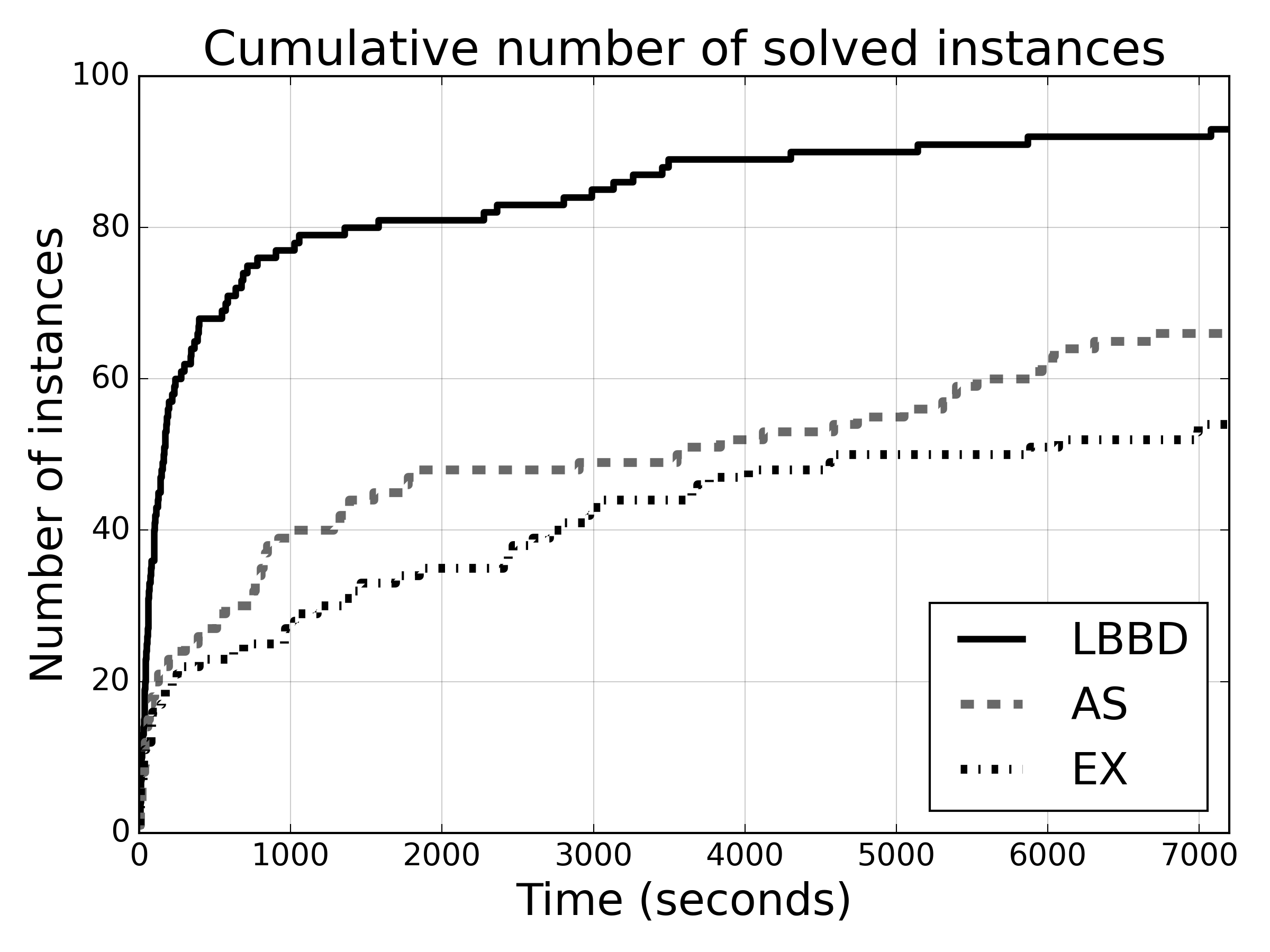}
        \label{fig:opt_gaps}
    \end{subfigure}
    \hfill
    \begin{subfigure}[b]{0.49\textwidth}
        \centering
        \includegraphics[width=\textwidth]{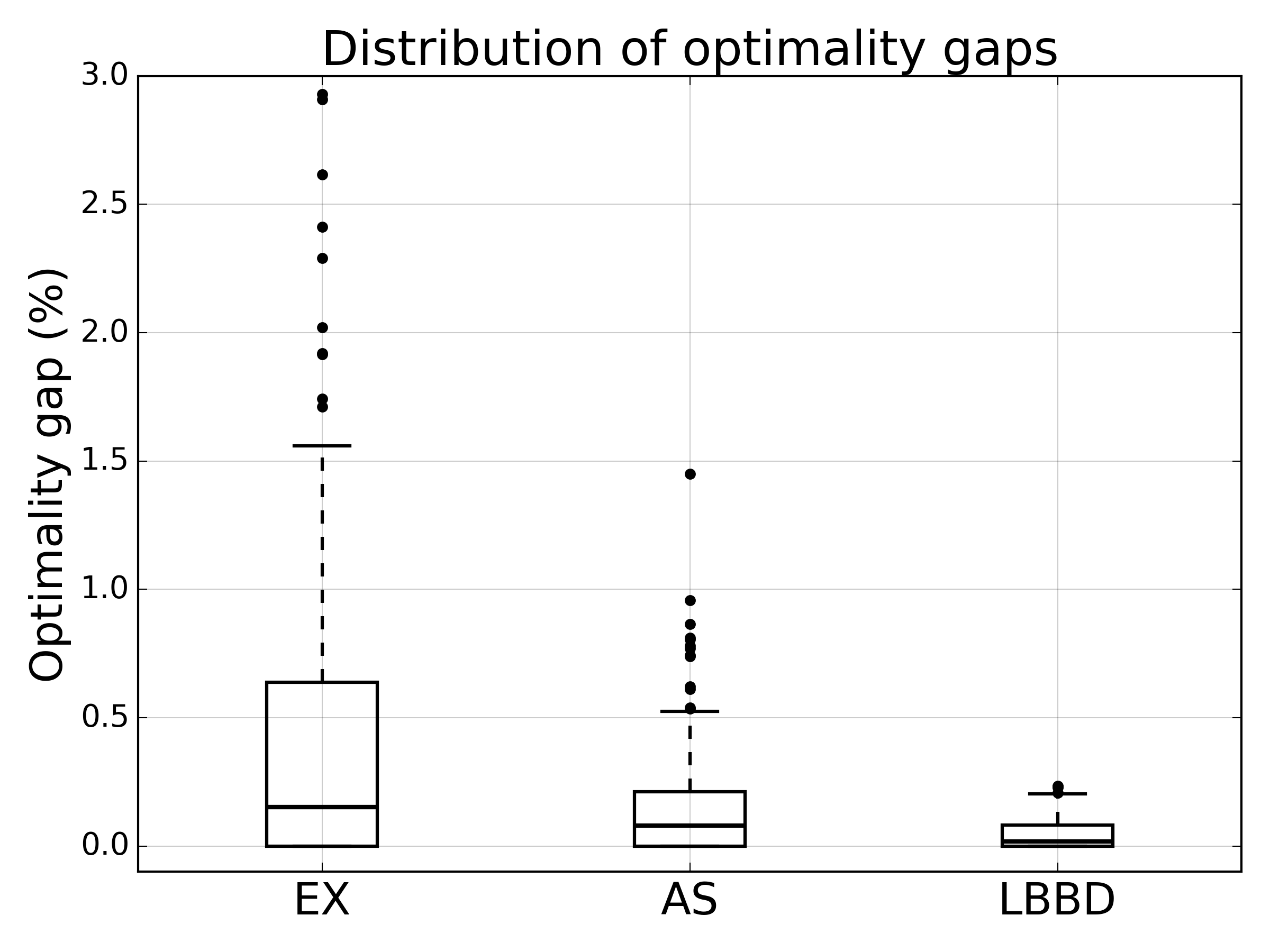}
        \label{fig:cum_solved}
    \end{subfigure}
    \captionsetup{skip=-8pt}
    \label{fig:results_partial_networks}
\end{figure}
\vspace{-1.5cm}

Figure~\ref{fig:results_partial_networks} shows, for each method, the number of instances solved to optimality within $t$ seconds (left), and a box plot of the final optimality gaps (right). Table \ref{tab:results_summary} separates the results by on-route charging availability and number of routes. We report the number of solved instances, the average optimality gap, and the count of instances on which each algorithm provides the best performance\begin{revblock}, defined as the shortest solving time or, for instances left unsolved by all methods, the smallest optimality gap at termination\end{revblock}.

\begin{table}[H]
\caption{Results per group of instances for partial networks}
  \centering
  \begin{minipage}{0.48\textwidth}
   \captionsetup{type=table,labelformat=empty}
   \vspace{-0.5cm}
  \caption*{With on-route charging}
  \vspace{-0.1cm}
   \centering
\resizebox{1\textwidth}{!}{%
  \setlength{\tabcolsep}{4pt}
  \begin{tabular}{cccc ccc ccc}
    \toprule
    \multirow{2}{*}{\specialcell{$|\mathcal{R}|$}}  & \multicolumn{3}{c}{EX} & \multicolumn{3}{c}{AS} & \multicolumn{3}{c}{LBBD} \\
    \cmidrule(lr){2-4} \cmidrule(lr){5-7} \cmidrule(lr){8-10}
                   & Opt & Gap & Best & Opt & Gap & Best & Opt & Gap & Best \\
    \midrule
    3      & 18 & 0.01  & 3    & 18 & 0.01  & 0    & \textbf{19} & \textbf{0.00} & \textbf{17} \\
    6      & 13 & 0.18  & 1    & 13 & 0.09  & 0    & \textbf{17} & \textbf{0.00} & \textbf{19} \\
    9      & 2  & 0.29  & 0    & 8  & 0.09  & 1    & \textbf{13} & \textbf{0.01} & \textbf{19} \\
    12     & 1  & 0.42  & 0    & 3  & 0.19  & 1    & \textbf{7}  & \textbf{0.04} & \textbf{19} \\
    15     & 1  & 0.46$^*$ & 0 & \textbf{4} & 0.22  & 4    & \textbf{4}  & \textbf{0.06} & \textbf{16} \\
    \midrule
    Tot/avg & 35 & 0.27  & 4    & 46 & 0.12  & 6    & \textbf{60} & \textbf{0.02} & \textbf{90} \\
    \bottomrule
  \end{tabular}%
}

   \vspace{-0.2cm}
   \caption*{\hspace{-0.4cm}\footnotesize{\textmd{$^*$ excludes an instance for which no solution was found \hfill}}}
  \end{minipage}%
  \hfill
  \begin{minipage}{0.474\textwidth}
   \captionsetup{type=table,labelformat=empty}
   \vspace{0.021cm}
  \caption*{Without on-route charging}
  \vspace{-0.1cm}
    \centering
\resizebox{1\textwidth}{!}{%
  \setlength{\tabcolsep}{4pt}
  \begin{tabular}{cccc ccc ccc}
    \toprule
    \multirow{2}{*}{\specialcell{$|\mathcal{R}|$}} & \multicolumn{3}{c}{EX} & \multicolumn{3}{c}{AS} & \multicolumn{3}{c}{LBBD} \\
    \cmidrule(lr){2-4} \cmidrule(lr){5-7} \cmidrule(lr){8-10}
     & Opt & Gap & Best & Opt & Gap & Best & Opt & Gap & Best \\
    \midrule
    2      & 12  & 0.04 & 9   & 11  & 0.04 & 0   & \textbf{17} & \textbf{0.01} & \textbf{11} \\
    4      & 5   & 0.33 & 2   & 6   & 0.24 & 0   & \textbf{11} & \textbf{0.03} & \textbf{18} \\
    6      & 1   & 0.56 & 0   & 2   & 0.20 & 1   & \textbf{5}  & \textbf{0.06} & \textbf{19} \\
    8      & \textbf{1} & 0.89 & 0   & \textbf{1} & 0.25 & 4   & \textbf{1}  & \textbf{0.10} & \textbf{16} \\
    10     & 0   & 1.12 & 0   & 0   & 0.24 & 5   & 0   & \textbf{0.15} & \textbf{15} \\
    \midrule
    Tot/avg & 19 & 0.59 & 11  & 20  & 0.19 & 10  & \textbf{34} & \textbf{0.07} & \textbf{79} \\
    \bottomrule
  \end{tabular}%
} 
    \vspace{-0.2cm}
    \caption*{\footnotesize{\phantom{$^*$ excludes an instance for which no solution was found}}}
  \end{minipage}
\label{tab:results_summary}
\end{table}
\vspace{-0.5cm}
\vspace{-0.5cm}

The results show that LBBD significantly outperforms the other methods. It solves 94 instances, compared to 54 for EX and 66 for AS, requires only 2.5\% and 5.3\% of the time limit to match the number of instances solved by the other two methods, and achieves the best performance in 169 out of 200 instances. The largest optimality gap obtained with LBBD is 0.23\%, whereas 50 instances exhibit an optimality gap of 0.64\% or higher (up to 2.93\%) with EX, and of 0.21\% or higher (up to 1.45\%) with AS. AS maintains the most stable performance as the number of routes increases, particularly when on-route BEBs are prohibited. In this case, meeting electrification targets requires larger fleets of depot BEBs, and the scheduling simplifications performed in the first heuristic phase greatly reduce the computational effort needed to identify good-quality solutions. Nevertheless, AS provides limited improvements over EX in the number of instances solved to optimality.

\subsection{Heuristics on large instances}
\label{sec:Experiments:large}
In this section, we evaluate our heuristic method AS on large-scale electrification planning problems. We construct our instances from the complete bus networks of eight US cities, with a transition horizon of $|\mathcal{P}| = 10$ years. We also consider a simpler policy restriction (PR) heuristic, where the arc selection phase of algorithm AS is replaced by the most efficient a priori arc elimination rules identified in preliminary experiments: during service times, the depot BEBs are only allowed to start charging from the fully depleted battery state (i.e., $z^{pt}_{rbis}=0$ if $s > 0$ and $d^{pt}_r > 0$), and can only idle in the fully charged state (i.e., $v^{pt}_{rbs}=0$ if $s< s_b$ and $d^{pt}_r > 0$), while they cannot be in service when there is no demand (i.e., $w^{pt}_{rbs}=0$ if $d^{pt}_r = 0$). The problem is first solved with these restrictions, and the incumbent solution obtained from the restricted model is then used to warm-start the original formulation. This approach is intended to serve as a fair baseline method against AS. For both methods, 75\% of the time budget is allocated to the heuristic phase.

\begin{table}[H]
\centering
\caption{Results for complete networks}
\resizebox{1\textwidth}{!}{%
  \setlength{\tabcolsep}{7pt}
  \renewcommand{\arraystretch}{0.75}
  \begin{tabular}{l rrr rrr rrr rrr}
    \toprule
    \multirow{3}{*}{\specialcell{City}} & \multicolumn{6}{c}{With on-route charging} & \multicolumn{6}{c}{Without on-route charging} \\
	\cmidrule(lr){2-7} \cmidrule(lr){8-13}
    & \multicolumn{2}{c}{EX} & \multicolumn{2}{c}{PR} & \multicolumn{2}{c}{AS} & \multicolumn{2}{c}{EX} & \multicolumn{2}{c}{PR} & \multicolumn{2}{c}{AS} \\
    \cmidrule(lr){2-3} \cmidrule(lr){4-5} \cmidrule(lr){6-7} \cmidrule(lr){8-9} \cmidrule(lr){10-11} \cmidrule(lr){12-13} 
    & \multicolumn{1}{c}{UB} & \multicolumn{1}{c}{Gap} & \multicolumn{1}{c}{UB} & \multicolumn{1}{c}{Gap} &\multicolumn{1}{c}{UB} & \multicolumn{1}{c}{Gap} &\multicolumn{1}{c}{UB} & \multicolumn{1}{c}{Gap} &\multicolumn{1}{c}{UB} & \multicolumn{1}{c}{Gap} &\multicolumn{1}{c}{UB} & \multicolumn{1}{c}{Gap}   \\
    \midrule
    Atlanta     & \multicolumn{1}{c}{$-$}    & \multicolumn{1}{c}{$-$}    & 964.36   & 0.63                     & \textbf{962.32}         & \textbf{0.42} & \multicolumn{1}{c}{$-$} & \multicolumn{1}{c}{$-$} & 1103.09 & 6.45 & \textbf{1045.61} & \textbf{1.30} \\
    Boston      & \multicolumn{1}{c}{$-$}    & \multicolumn{1}{c}{$-$}    & 2384.71  & \multicolumn{1}{c}{$*$}  & \textbf{2379.87}        & \textbf{0.81} & \multicolumn{1}{c}{*}   & \multicolumn{1}{c}{*}   & 2561.54 & 7.15 & \textbf{2397.25} & \textbf{0.78} \\
    Chicago     & \multicolumn{1}{c}{$-$}    & \multicolumn{1}{c}{$-$}    & 2734.13  & 0.83                     & \textbf{2727.72}        & \textbf{0.59} & \multicolumn{1}{c}{$-$} & \multicolumn{1}{c}{$-$} & 2890.68 & 5.41 & \textbf{2744.84} & \textbf{0.39} \\
    Dallas      & \multicolumn{1}{c}{$-$}    & \multicolumn{1}{c}{$-$}    & 1205.69  & 0.85                     & \textbf{1200.61}        & \textbf{0.44} & \multicolumn{1}{c}{$-$} & \multicolumn{1}{c}{$-$} & 1319.92 & 4.86 & \textbf{1262.80} & \textbf{0.56} \\
    Detroit     & 439.06                     & 0.83                       & 439.23   & 0.60                     & \textbf{438.46}         & \textbf{0.43} & \multicolumn{1}{c}{$-$} & \multicolumn{1}{c}{$-$} & 459.80  & 2.24 & \textbf{451.84}  & \textbf{0.51} \\
    Houston     & 1972.56                    & 17.51                      & 1644.14  & 1.03                     & \textbf{1630.17}        & \textbf{0.18} & \multicolumn{1}{c}{$-$} & \multicolumn{1}{c}{$-$} & 1763.37 & 6.23 & \textbf{1659.81} & \textbf{0.38} \\
    Las Vegas   & \multicolumn{1}{c}{$-$}    & \multicolumn{1}{c}{$-$}    & 634.02   & 0.09                     & \textbf{634.01}         & \textbf{0.08} & \multicolumn{1}{c}{$-$} & \multicolumn{1}{c}{$-$} & 690.75  & 1.01 & \textbf{686.69}  & \textbf{0.43} \\
    Los Angeles & 4078.27                    & 18.75                      & 3344.34  & 0.95                     & \textbf{3325.93}        & \textbf{0.40} & \multicolumn{1}{c}{$-$} & \multicolumn{1}{c}{$-$} & 3556.79 & 6.04 & \textbf{3361.52} & \textbf{0.59} \\
    \bottomrule
  \end{tabular}%
}  
\vspace{-0.2cm}
{\captionsetup{type=table,labelformat=empty}
\caption*{\footnotesize{$-$ no feasible solution was found within the eight-hour time limit, $*$ out-of-memory error }}}
\label{tab:results_full_unique}
\end{table}
\vspace{-0.5cm}

Table \ref{tab:results_full_unique} reports upper bounds (millions of dollars) and optimality gaps (percentage), for EX, PR, and AS, on instances including and excluding on-route charging. The results indicate that sparsifying the depot BEBs scheduling graphs can provide good solutions for otherwise intractable instances. AS dominates, most noticeably when on-route charging is unavailable. Since both PR and AS limit the efficiency of the depot BEBs by restricting the flexibility of their schedules, but do not affect on-route BEBs, it was expected that they would provide solutions of lower quality in the second group of instances. The relatively stable performance of AS across instances with and without on-route charging confirms that the active arcs from the collection of subproblems we designed suffice to construct higher-quality schedules than the a priori restriction rules of PR.

Interestingly, AS achieves optimality gaps that are not much higher for complete networks than for the partial networks of the previous section. This is in part because the lower bound Gurobi achieves at the root node tends to become better as instance size increases, although it usually stops making progress once the solver starts branching on large instances. In Benders decomposition, the initial master problem shows the opposite trend. Moreover, solving each cut separation problem \eqref{model:closest_cut} can take several minutes, making each iteration very costly. \begin{revblock}Consequently, when we attempted to run LBBD on the complete network instances, it usually failed to find feasible solutions within the eight-hour time limit. At this scale, sparsifying the operational graphs becomes essential in practice, which restricts us to heuristic approaches. Our exact and heuristic methods thus play complementary roles: LBBD dominates on the small-scale instances representative of partial electrification projects, while AS extends the reach of our framework to citywide planning.\end{revblock} Further discussions and detailed results for each phase of algorithms AS and PR are provided in \ref{ec:heuristics}.

\subsection{Case study of Chicago} \label{sec:Experiment:chicago}
This section presents a case study of the Chicago bus network. We examine the solution obtained in Section~\ref{sec:Experiments:large} from the AS heuristic with on-route charging. We discuss (i) investment timing; (ii) fleet composition and vehicle utilization; and (iii) charger deployment. \begin{revblock}\ref{ec:rolling_horizon} extends the case study with a rolling-horizon implementation of our model in the context of price uncertainty\end{revblock}.

\begin{revblock}
\paragraph{Purpose and limitations.}
This case study illustrates the managerial insights enabled by our model under a parsimonious and interpretable parameterization, rather than prescribing an implementation plan for CTA. Actual planning applications should calibrate the model to agency-specific data, including local financial assumptions, electrification targets, the existing fleet and infrastructure, projected demand growth, and route-dependent energy consumption and terminal layover patterns.\end{revblock}

\paragraph{Solution summary.}
Overall, the initial fleet of 1495 conventional buses is replaced by a fleet of the same size consisting of 788 on-route BEBs, 690 short-range depot BEBs, and 17 long-range depot BEBs, supported by 106 on-route chargers and 157 depot chargers, for \$1.52B in investments and \$1.66B in operational costs over 10 years (undiscounted totals). Due to the cheaper operating costs of BEBs compared to conventional buses, the operational costs decrease throughout the transition. The initial annual operational costs of \$240.4M reach \$154.5M in year 5, and \$141.2M in year 10.
\begin{figure}[H]
\caption{Planning solution and costs summary}
    \label{fig:chicago_planning}
    \centering
    \begin{subfigure}{0.49\textwidth}
        \includegraphics[width=\linewidth]{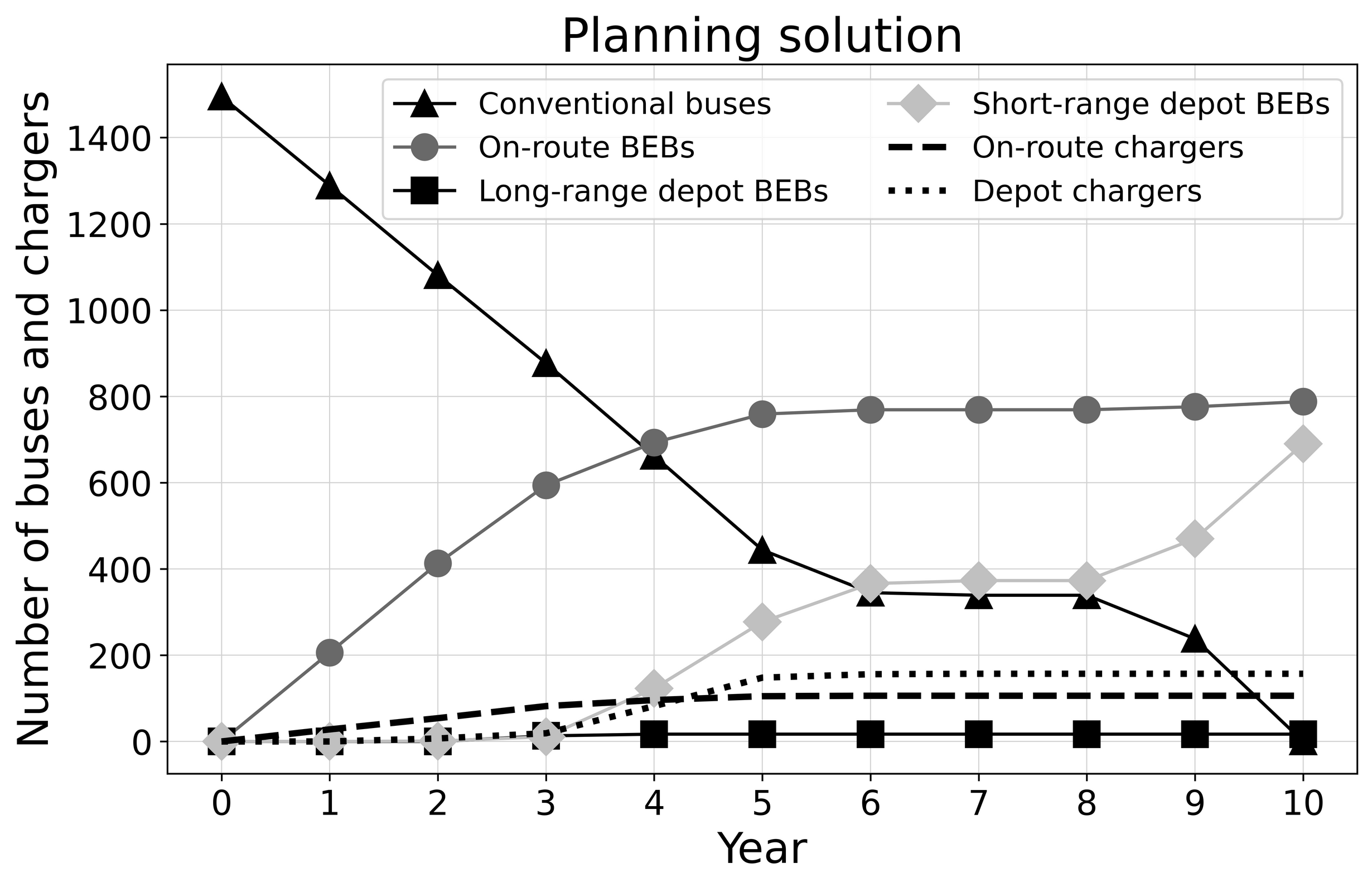}
        \label{fig:chicago_planning:1}
    \end{subfigure}
    \hfill
    \begin{subfigure}{0.482\textwidth}
        \includegraphics[width=\linewidth]{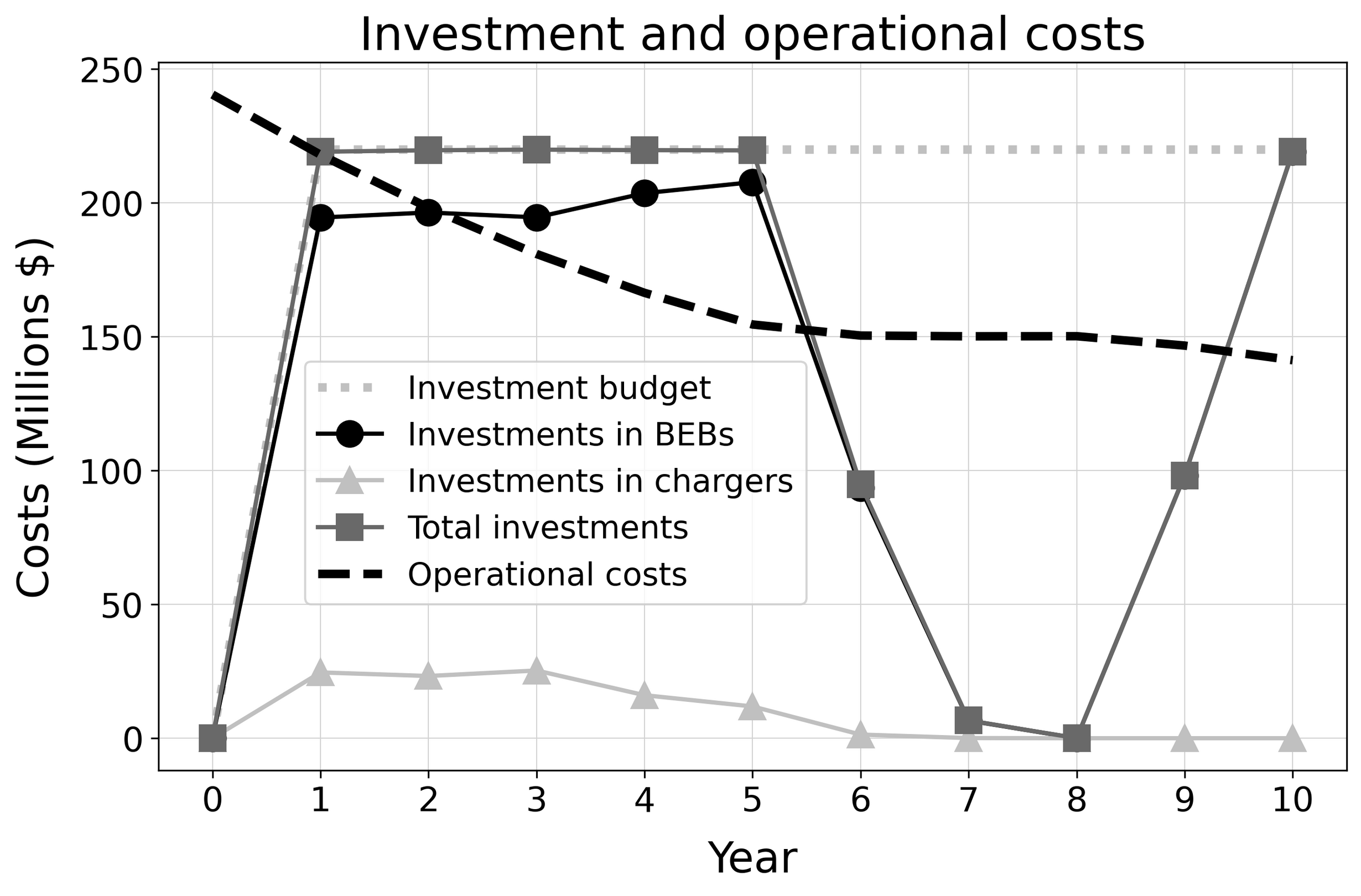}
        \label{fig:chicago_planning:2}
    \end{subfigure}
    \captionsetup{skip=-8pt}
\end{figure}
\vspace{-1.65cm}

\paragraph{Investment timing and profitability.}
In this instance, operating conditions are stationary across investment periods (objective coefficients and service requirements are constant), and electrification targets are enforced only in the last period. In this setting, the timing of investments is driven by their profitability. For each period $p \in \mathcal{P}$ with new investments, let
\begin{equation}
    \mathrm{ROI}_p := \frac{O_p(x_{p-1}) - O_p(x_p)}{I_p(x_p-x_{p-1})}
\end{equation}
denote the yearly operational savings generated per dollar of capital expenditure. Advancing an investment by one period yields one additional year of operational savings, but also requires paying its capital cost one year earlier in the discounted objective. This trade-off induces a profitability threshold at $1-\gamma$. We can verify that if the investments implemented in period $p$ can be shifted to period $p-1$ (respectively $p+1$) without violating the budget constraints, then this shift improves the discounted objective if and only if $\mathrm{ROI}_p > 1-\gamma$ (respectively $< 1-\gamma$). As a result, the optimal solution follows a three-phase structure: an initial phase where high-return investments are front-loaded, a middle phase where no investments occur, and a final phase where the remaining investments are driven by the electrification target and delayed toward the end of the horizon.

The solution is consistent with this structure. The first phase spans years 1--6, with returns ranging from $\mathrm{ROI}_{1}=10.26\%$ to $\mathrm{ROI}_{6}=4.36\%$ (here, $1-\gamma = 4\%$). The yearly budget of \$200M is almost fully used in the first five periods, and \$94.8M is spent in year 6. We also observe some investments in the following year, with $\mathrm{ROI}_{7}=4.30\%$, highlighting that the solution is slightly suboptimal. We can verify that advancing these investments to year 6 and adjusting the operations accordingly improves the value of the solution by \$16.5k (less than 0.001\% of the total costs). After a period without investment, the last phase covers years 9 and 10. As these last investments fall below the profitability threshold ($\mathrm{ROI}_9 = 3.57\%$ and $\mathrm{ROI}_{10} =2.50\%$), they are delayed as much as possible.

\paragraph{Fleet composition and vehicle utilization.}

\begin{figure}[H]
\caption{Number of buses in service per hour by year}
\label{fig:chicago_schedule}
    \centering
    \begin{subfigure}{0.32\textwidth}
        \includegraphics[width=\linewidth]{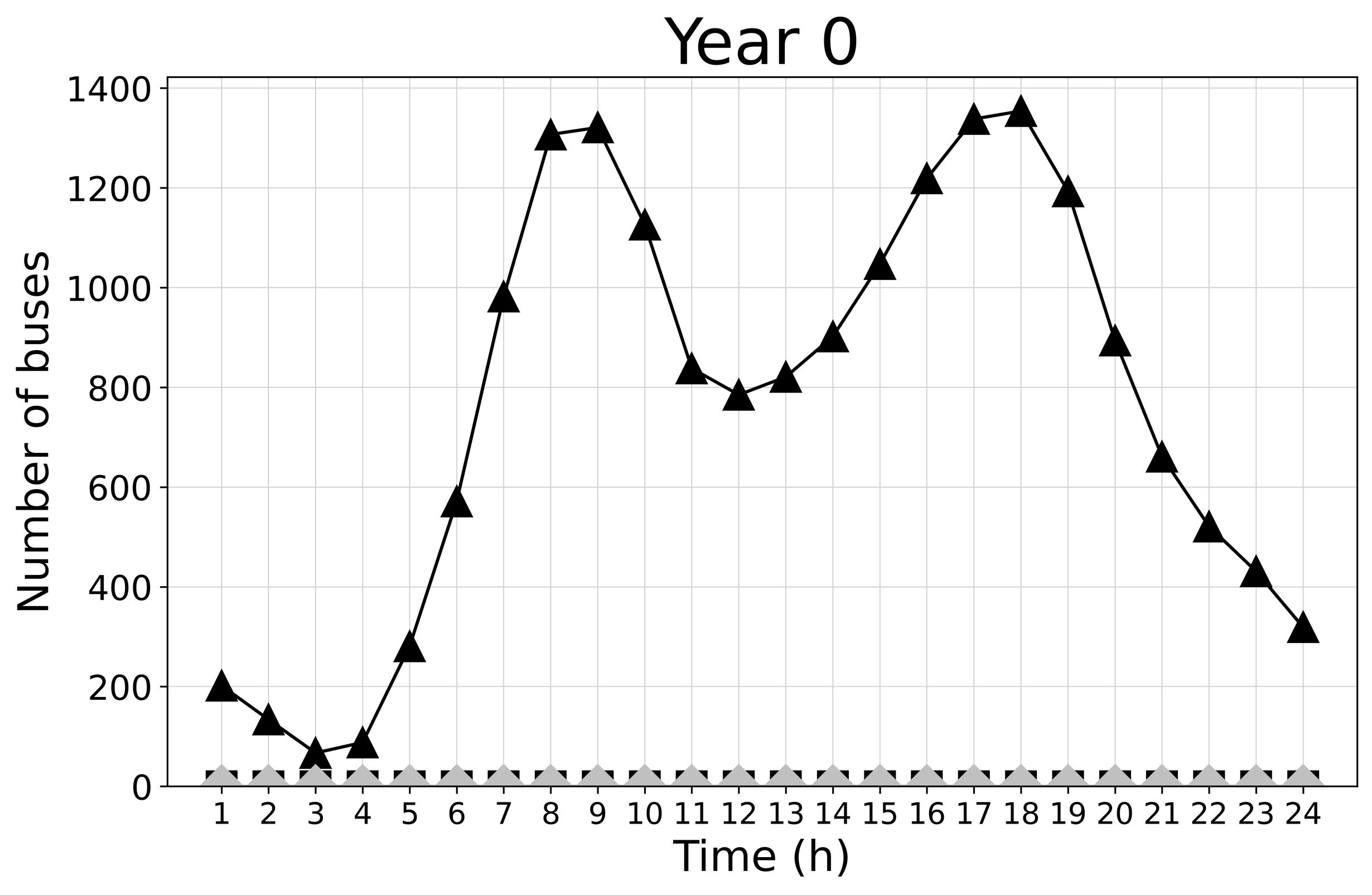}
        \label{fig:chicago_schedule:1}
    \end{subfigure}
    \hfill
    \begin{subfigure}{0.32\textwidth}
        \includegraphics[width=\linewidth]{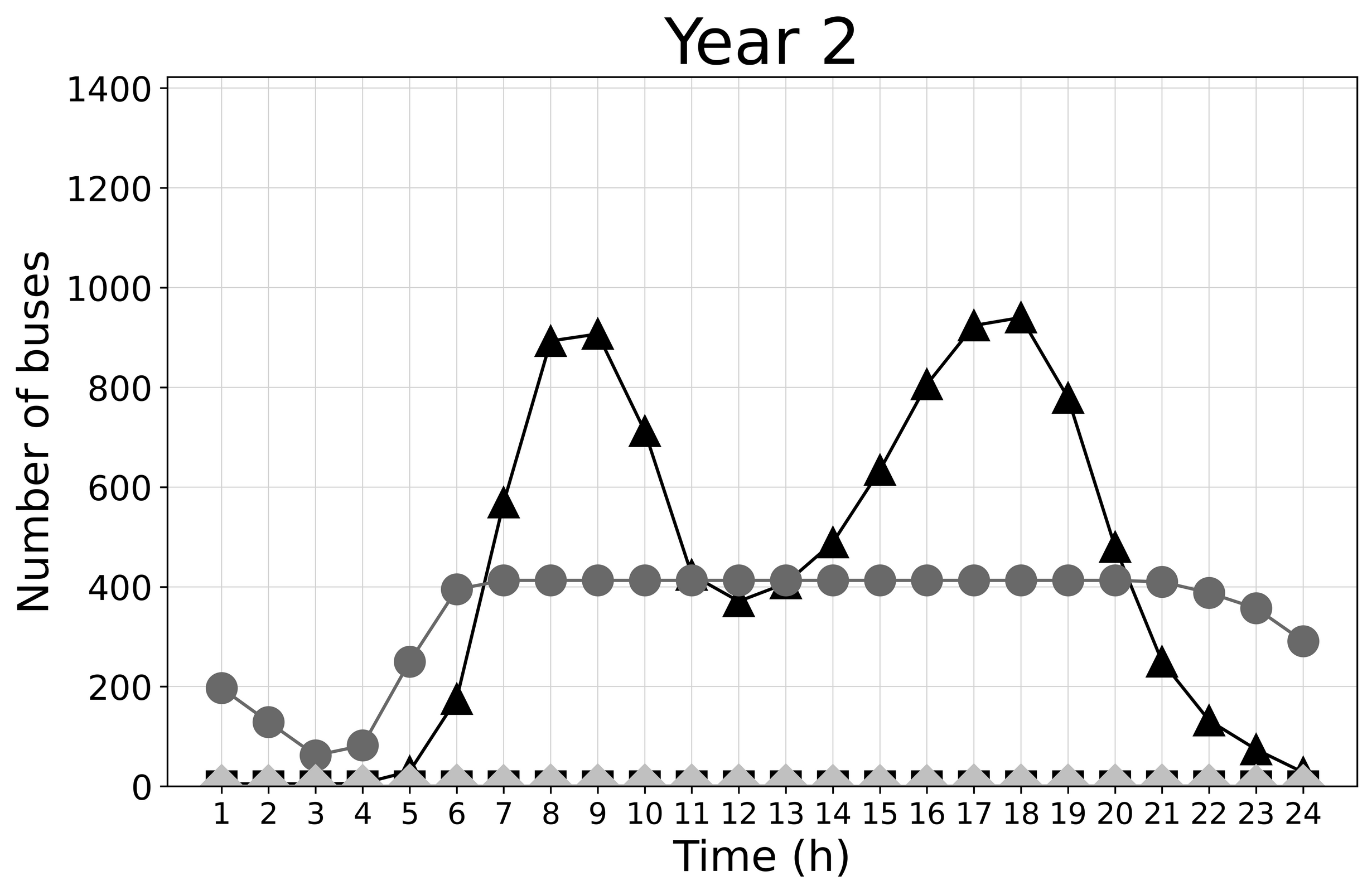}
        \label{fig:chicago_schedule:2}
    \end{subfigure}
    \hfill
    \begin{subfigure}{0.32\textwidth}
        \includegraphics[width=\linewidth]{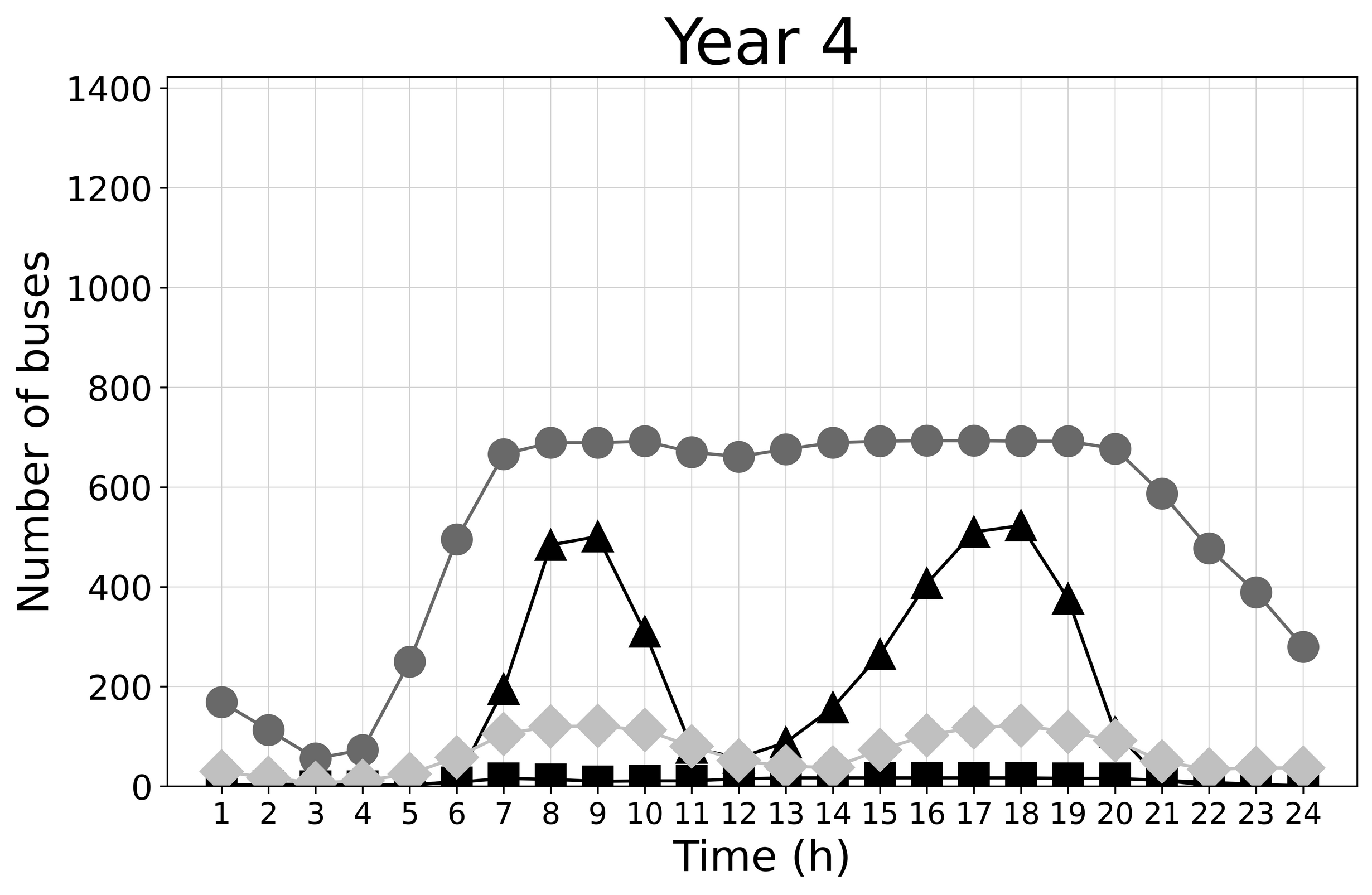}
        \label{fig:chicago_schedule:3}
    \end{subfigure}

    \vspace{-0.2cm}
    \begin{subfigure}{0.32\textwidth}
        \includegraphics[width=\linewidth]{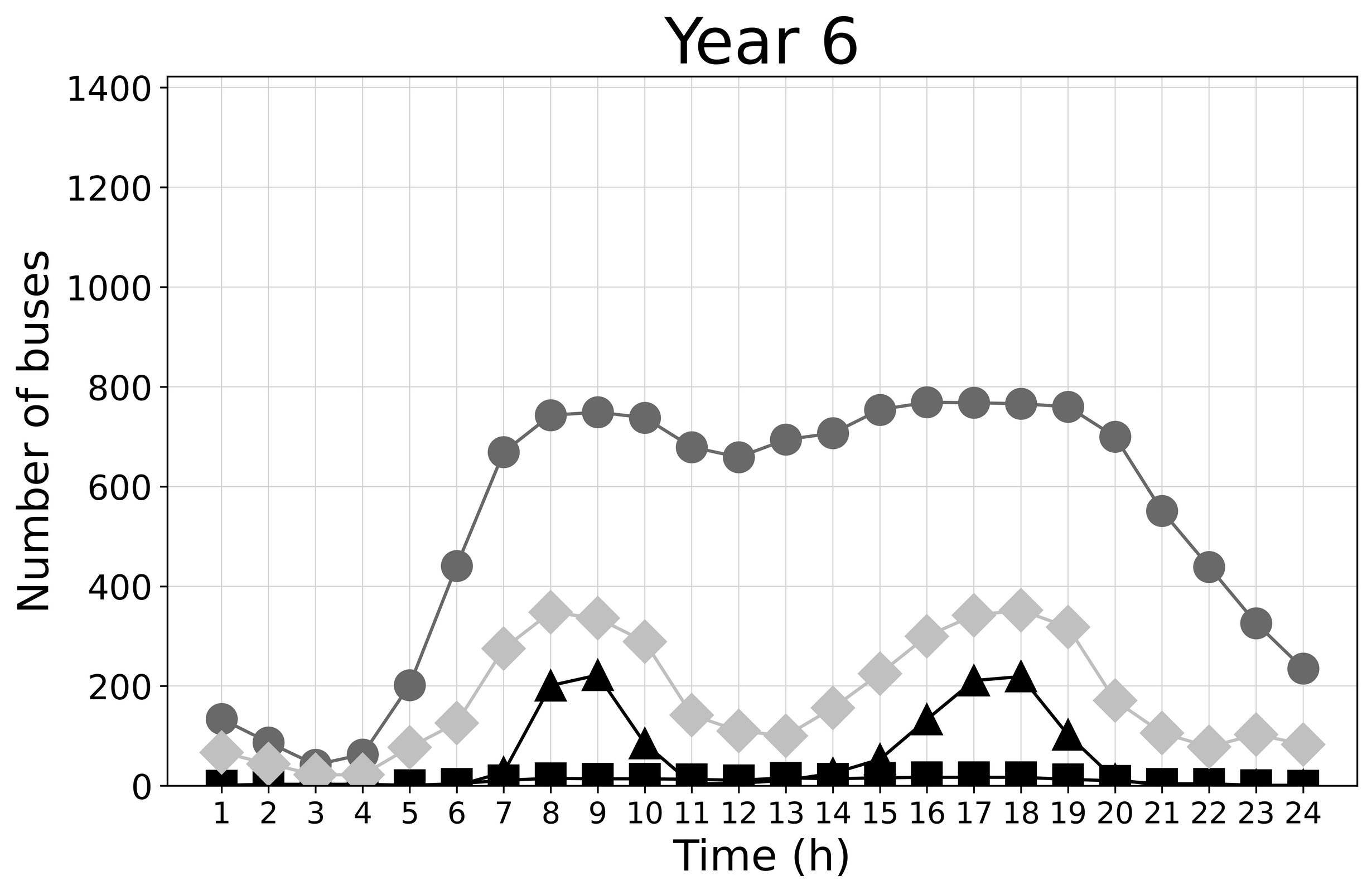}
        \label{fig:chicago_schedule:4}
    \end{subfigure}
    \hfill
    \begin{subfigure}{0.32\textwidth}
        \includegraphics[width=\linewidth]{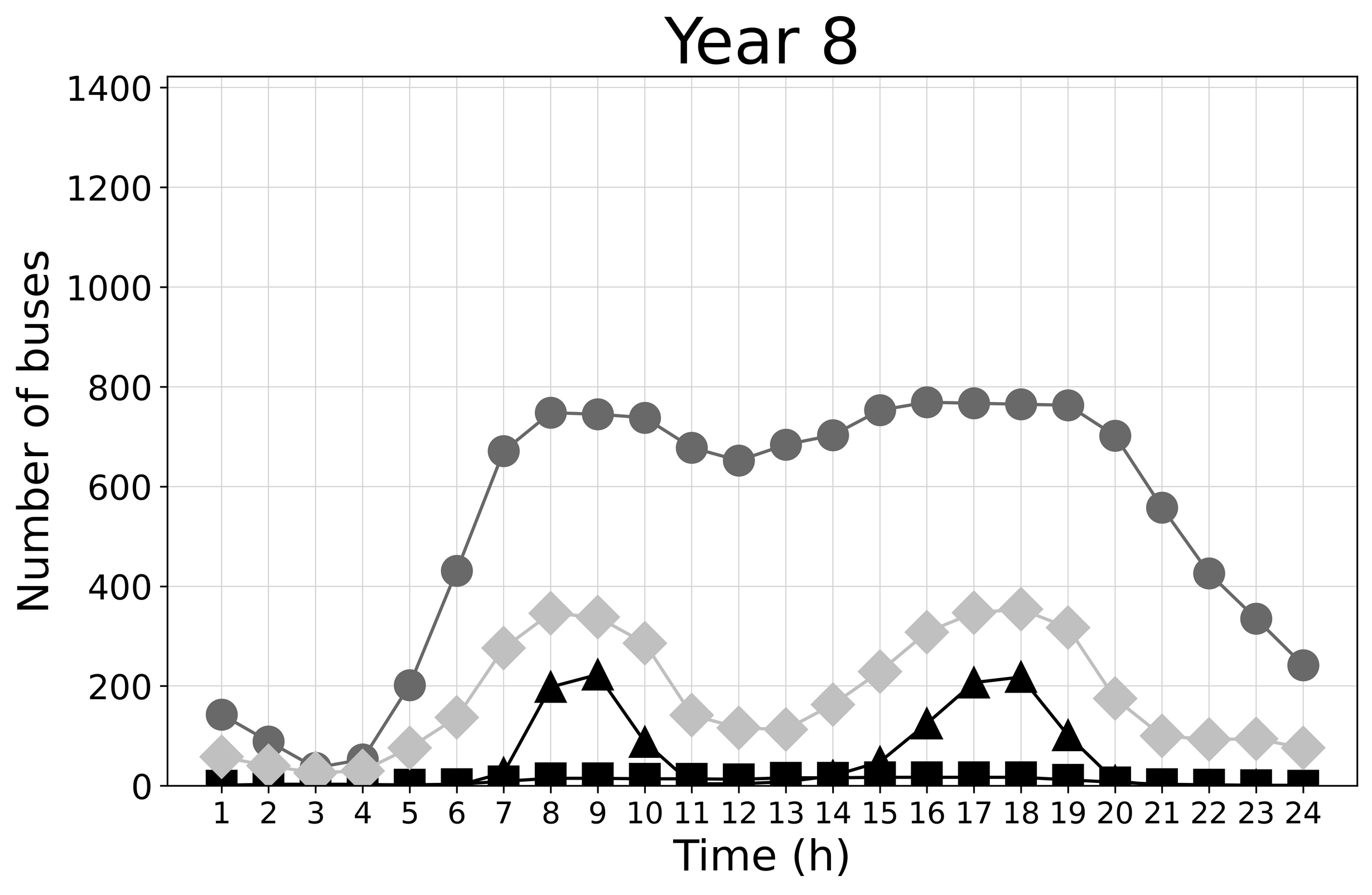}
        \label{fig:chicago_schedule:5}
    \end{subfigure}
    \hfill
    \begin{subfigure}{0.32\textwidth}
        \includegraphics[width=\linewidth]{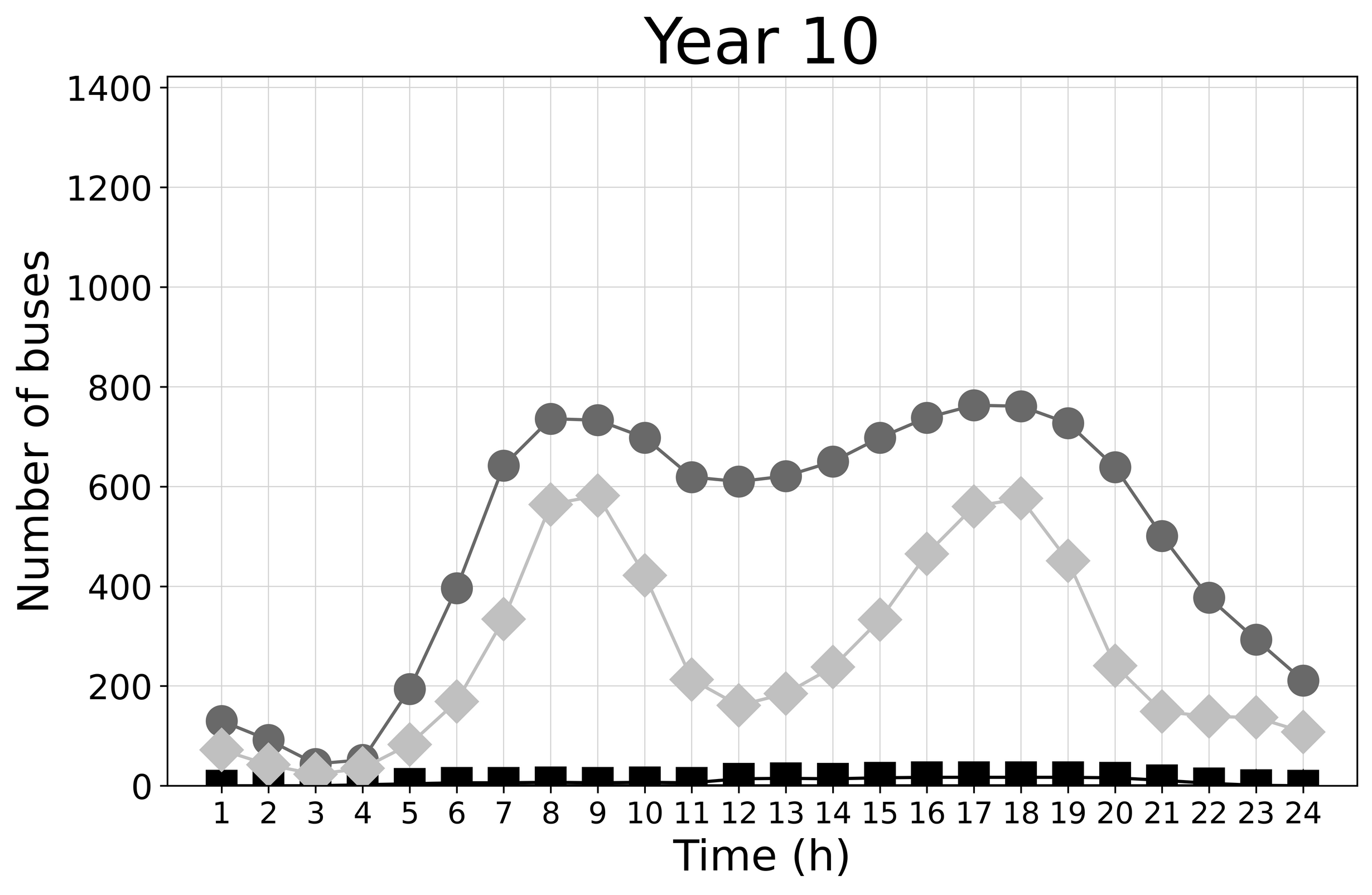}
        \label{fig:chicago_schedule:6}
    \end{subfigure}

    \vspace{-0.4cm}
    \begin{subfigure}{\textwidth}
        \includegraphics[width=\linewidth]{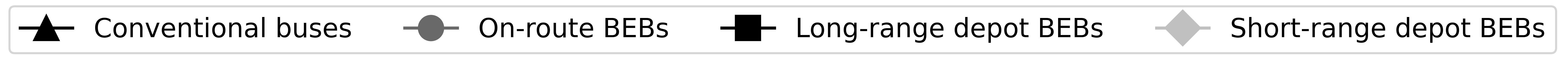}
        \label{fig:chicago_schedule:legend}
    \end{subfigure}
\end{figure}
\vspace{-1.3cm}

In the first three years, all the investments are made in on-route BEBs and fast chargers. Depot BEBs are progressively introduced starting from year 4, and represent nearly all of the BEBs acquired in the last phase of investments. This can be explained by the usage profile presented in Figure \ref{fig:chicago_schedule}. In year 2, the 413 on-route BEBs are in service 84.2\% of the time, and are all in service for 14 consecutive hours starting from 5 am. In comparison, depot BEBs can be in service at most 66.7\% of the time, since we assumed they require one hour of charging for every two hours of operation. Since the hourly operating costs of conventional buses are by far the highest, on-route BEBs are preferable to depot BEBs as long as conventional buses are still needed to satisfy the base load. After that, short-range depot BEBs, which are the cheapest to acquire and operate, become the most cost-efficient option. In the final schedule, they are primarily in service during peak hours, and recharge mainly at midday and overnight. 

\begin{figure}[H]
\caption{Charger deployment}
\vspace{-0.3cm}
\label{fig:chicago_chargers}
    \centering
    \begin{subfigure}{0.32\textwidth}
        \captionsetup{font=small}
        \caption*{Year 2}
        \includegraphics[width=\linewidth]{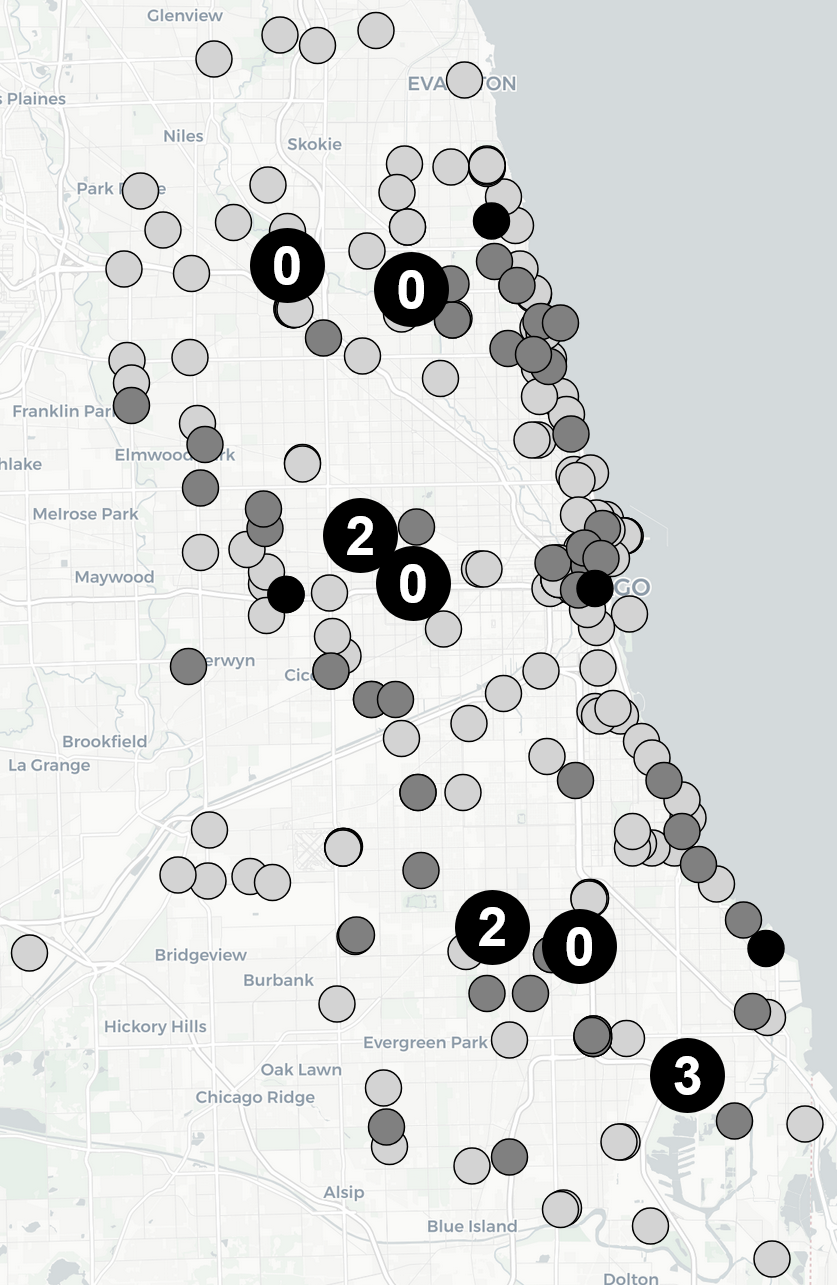}
        \captionsetup{skip=0pt} 
    \end{subfigure}
    \hfill
    \begin{subfigure}{0.32\textwidth}
        \captionsetup{font=small}
        \caption*{Year 4}
        \includegraphics[width=\linewidth]{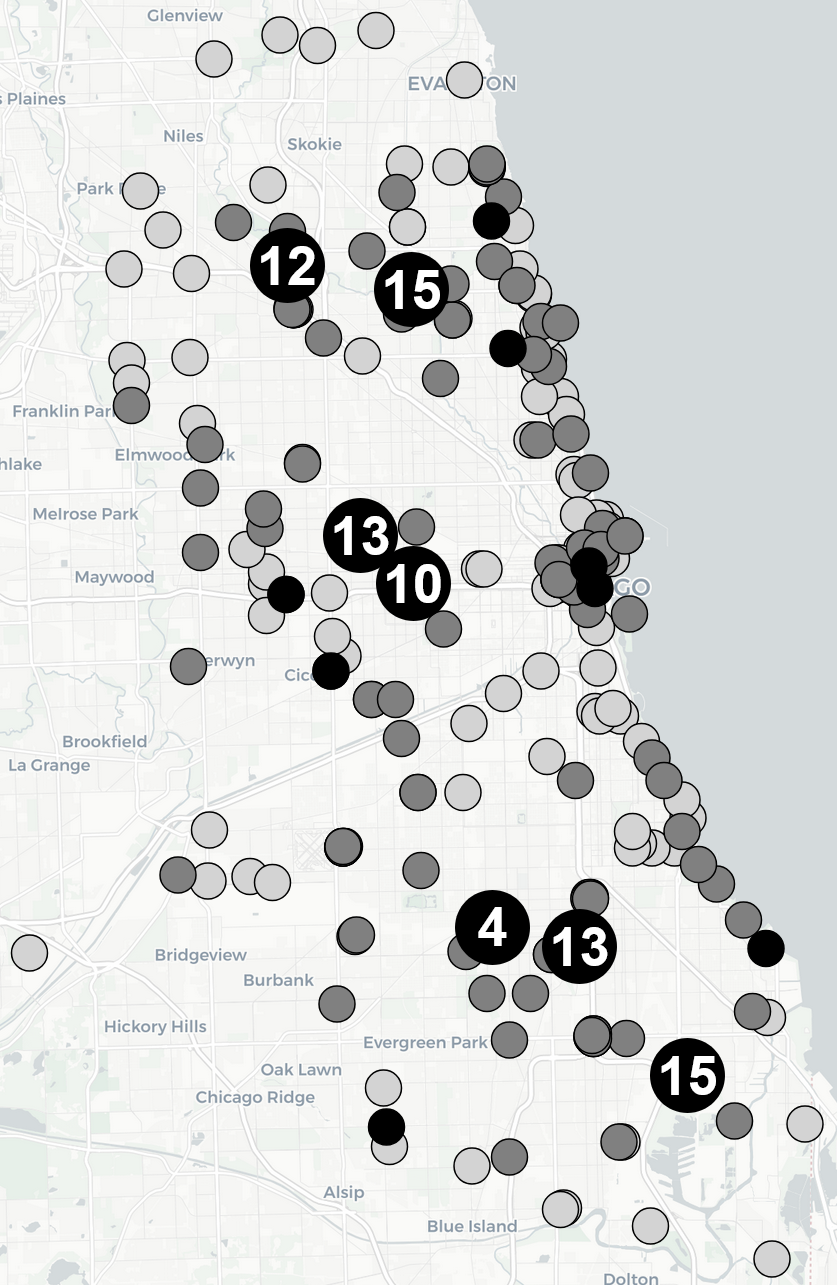}
        \captionsetup{skip=0pt} 
    \end{subfigure}
    \hfill
    \begin{subfigure}{0.32\textwidth}
        \captionsetup{font=small}
        \caption*{Year 6}
        \includegraphics[width=\linewidth]{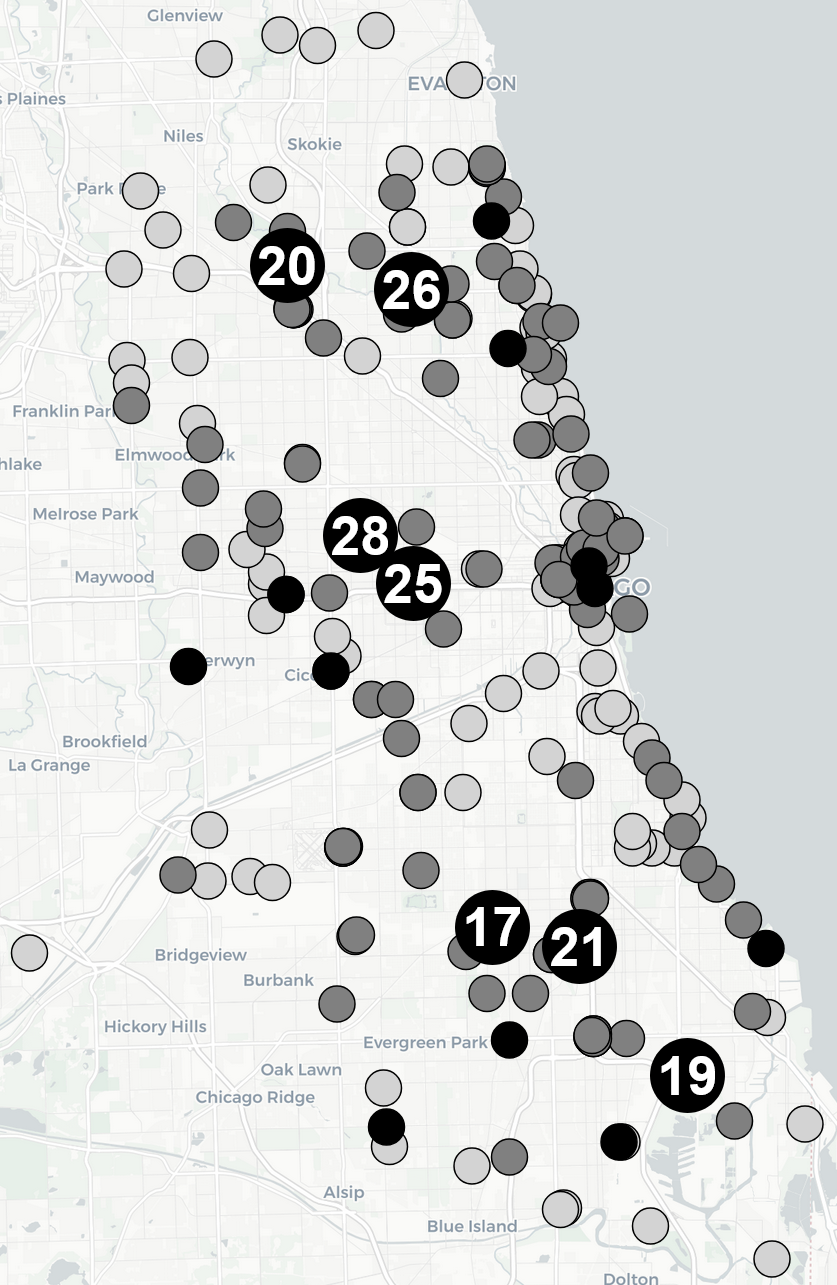}
        \captionsetup{skip=0pt} 
    \end{subfigure}
    \vspace{0.8cm}
    \begin{subfigure}{\textwidth}
        \includegraphics[width=\linewidth]{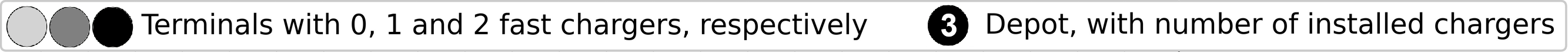}
        \label{fig:chicago_chargers:legend}
    \end{subfigure}
\end{figure}
\vspace{-2.0cm}

\paragraph{Charger deployment.}
Figure \ref{fig:chicago_chargers} shows that fast chargers are primarily deployed in densely connected areas. Terminals located on the city outskirts are never equipped with fast chargers, meaning that the routes outside of the city center are mostly serviced by depot BEBs. In this solution, no chargers are installed after the first investment phase. Indeed, although more chargers could reduce the operational costs by increasing the scheduling flexibility of the depot BEBs introduced in years 9 and 10, these savings in the final years do not offset additional capital investment in the objective.

\section{Conclusion}
\label{sec:Conclusion}
\paragraph{Summary.} This paper presents an optimization framework for bus fleet electrification planning. The model is a two-stage integer program with integer subproblems, where the first stage controls the fleet composition and charger deployment over multiple investment periods, while the second stage optimizes daily operations through a flow-based scheduling model. We introduce: (i) an exact logic-based Benders decomposition equipped with preprocessing, master problem strengthening, and cut separation from multiple relaxations; and (ii) a restriction heuristic that progressively relaxes additional scheduling constraints selected from auxiliary problems. Our acceleration techniques yield speedups of three orders of magnitude compared to a baseline logic-based Benders decomposition, and our heuristic delivers high-quality solutions for citywide instances.

\paragraph{Policy relevance.} Beyond these methodological contributions, the proposed framework can provide practical value as a decision-support tool for both transit agencies and policymakers. For agencies, it produces an actionable transition plan grounded in a granular representation of operational costs and feasibility. For policymakers, it provides a quantitative basis to assess whether electrification is likely to be financially self-sustaining under prevailing technology costs and operating conditions, or whether targeted support is needed to accelerate adoption. Such assessments could ultimately help guide the design of electrification targets and fleet renewal policies.

\emph{\begin{revblock}Limitations\end{revblock} and extensions.} \begin{revblock}This work relies on simplifying assumptions. Some are structural, as our algorithms exploit them directly, while others admit extensions that preserve the key structural properties of our model. The no-interlining assumption, which restricts operational flexibility, is structural, as the single-route subproblems used in our algorithms assume that each vehicle serves a unique route. Similarly, the hourly discretization of a single representative day per investment period reflects a tractability limitation of time-expanded formulations, whose size grows linearly with the number of intervals. In contrast, as we detail in \ref{ec:partial-charging}\end{revblock}, richer charging policies than those assumed in the paper can be represented within the same flow-based operational model. \begin{revblock}Our model is also deterministic. As illustrated in \ref{ec:rolling_horizon}, deploying it in a rolling horizon mitigates this limitation by adapting the plan as information is revealed. However, this scheme reacts to new information rather than explicitly modeling uncertainty. Uncertainty in operational inputs (e.g., energy prices and labor costs) could be modeled via scenarios, yielding\end{revblock} separable subproblems that would benefit from decomposition while preserving the monotonicity of the operational value functions used by our exact method. \begin{revblock}Integrating\end{revblock} power-system considerations into the model, including distribution capacity limits and expansion decisions, vehicle-to-grid participation, and market-responsive charging\begin{revblock}, could also enrich the strategic and operational problems while preserving their core structure\end{revblock}. Algorithmically, embedding advanced EVSP solution methods into our logic-based Benders framework could enable the use of high-fidelity scheduling models for strategic planning or, conversely, allow key strategic decisions to be incorporated into EVSP formulations.



%
%
%



\bibliographystyle{informs2014} 
\bibliography{references} 


\EquationErrorMsgdimen=500pt
\EquationErrorBoxdimen=500pt
\setlength{\maxdseq}{500pt} 

\ECSwitch
\ECHead{Appendix to: \\ Bus Fleet Electrification Planning Through Logic-Based Benders Decomposition and Restriction Heuristics}

\begin{revblock}
\noindent\emph{Outline.} The appendix comprises ten sections. \ref{ec:Proofs} presents omitted results on the properties of the model. \ref{ec:Example_operational_model} illustrates the operational model of Section~\ref{sec:model:operations}. \ref{ec:benders_theory} reviews Benders cut selection methods and proves Proposition~\ref{prop3:closest_deepest}. \ref{ec:sec:Monotone_cuts_implementation} details the implementation of the monotone cuts of Section~\ref{sec:Algorithms:Benders:monotone}. \ref{ec:Experiments:instances} describes the benchmark instances and model parameters used in Section~\ref{sec:Experiments}. \ref{ec:Benders_cuts} empirically compares the reviewed Benders cut selection methods. \ref{ec:heuristics} presents detailed results on the restriction heuristics of Section~\ref{sec:Experiments:large}. \ref{ec:partial-charging} shows how the charging assumptions of Section~\ref{sec:Model:assumptions} can be relaxed within our framework and quantifies their impact on fleet sizing. \ref{ec:rolling_horizon} extends the Chicago case study of Section \ref{sec:Experiment:chicago} to illustrate a rolling-horizon implementation of the model under price uncertainty. \ref{ec:notation} provides a notation table.\end{revblock}

\section{Properties of the BFEP} \label{ec:Proofs}
In this section, \begin{revblock} we show that the LP relaxation of the operational problem \eqref{model:operations_detailed} is loose (Proposition~\ref{prop:lp_vs_ip})\end{revblock}, and prove the validity of the dual reduction based on terminal dominance (Proposition~\ref{prop1:bound_dominated_terminals}).

\begin{revblock}
\begin{proposition}\label{prop:lp_vs_ip} 
The operational problem \eqref{model:operations_detailed} can be infeasible while its LP relaxation is feasible for the same strategic solution $x_p\in\mathcal{X}_p$. In particular, the minimum fleet size required for feasibility can be 50\% larger for \eqref{model:operations_detailed} than for its LP relaxation.
\end{proposition}
\begin{proof}{Proof.}
Consider an instance with one period, $T=3$ intervals, one route, one depot, and only depot BEBs of the same type $b$ with capacity $s_b=2$. Dropping the fixed indices $p,r,b,i$ from the notation, consider constant service level requirements $d^0=d^1=d^2=1$ and charging times $\kappa_{0}=\kappa_{1}=2$. Consider $\bar\eta=2$ BEBs and $\bar{\chi}$ large enough for \eqref{model:operations_detailed:depot_capacity} to be nonbinding.

\emph{LP feasibility.} The fractional solution $w^t_{1}=w^t_{2}=z^t_{0}=\tfrac{1}{2}$ for all $t\in\mathcal{T}$, with all other variables zero, satisfies the service, flow-balance and fleet-size constraints \eqref{model:operations_detailed:service_level}--\eqref{model:operations_detailed:fleet_depot} by direct substitution.

\emph{IP infeasibility.} Suppose for contradiction that an integer solution exists. Set $Z_0:=\sum_{t \in \mathcal{T}} z^t_0$ and $Z_1:=\sum_{t \in \mathcal{T}} z^t_1$, both nonnegative integers. Summing the flow-balance equations~\eqref{model:operations_detailed:depot_flow_full} over $t\in\mathcal{T}$, the idling terms telescope and cancel by cyclicity, leaving $\sum_{t \in \mathcal{T}} w^t_2=Z_0+Z_1$. Similarly, \eqref{model:operations_detailed:depot_flow_empty} gives $\sum_{t \in \mathcal{T}} w^t_1=Z_0$. The total number of service trips is therefore $2Z_0+Z_1$, each occupying a bus for one interval, and the service requirement~\eqref{model:operations_detailed:service_level} imposes $2Z_0+Z_1\geq 3$. Each charging trip $z^t_s$ occupies a bus for $\kappa_s = 2$ intervals. Moreover, summing \eqref{model:operations_detailed:depot_flow_full}--\eqref{model:operations_detailed:depot_flow_empty} over $s$ shows that the left-hand side of \eqref{model:operations_detailed:fleet_depot} takes the same value in every interval, so this constraint also holds with any $t \in \mathcal{T}$ in place of $0$. Summing the resulting $T$ inequalities gives $\sum_{t \in \mathcal{T}}\sum_{s \in \mathcal{S}_b} v^t_s + (2Z_0+Z_1) + 2(Z_0+Z_1) \leq T\bar\eta = 6$, hence $4Z_0 + 3Z_1 \leq 6$. Subtracting this from three times the service inequality $2Z_0+Z_1\geq 3$ yields $2Z_0\geq 3$, hence $Z_0\geq 2$ by integrality, contradicting $4Z_0+3Z_1\leq 6$. Conversely, the circulation $w^0_2=w^1_1=z^2_0=v^1_2=w^2_2=v^0_1=z^1_1=1$, with all other variables zero, is feasible for $\bar{\eta}=3$. The minimum feasible fleet size is therefore exactly three. \halmos\end{proof}
\end{revblock}

\begin{proof}{Proof of Proposition \ref{prop1:bound_dominated_terminals}.}
Consider a solution such that $\widetilde{\chi}^{p}_{j} < \widetilde{\chi}^{\text{UB}}_{j}$ and $\widetilde{\chi}^{p}_{j'} > 0$ for some $j \in \mathcal{J},\ j' \in \mathcal{J}(j), \ p \in \mathcal{P}$. By the monotonicity constraints \eqref{model:BFEP:invest_charg_mono}, these conditions hold for a consecutive sequence of periods $\{p_1,\dots,p_2\}$, for $p_1 \leq p_2$. Since $\mathcal{R}(j) \supseteq \mathcal{R}(j')$, we can modify the solution without affecting the investment and operational costs by setting $\widetilde{\chi}^{p'}_{j} \leftarrow \widetilde{\chi}^{p'}_{j} + 1$ and $\widetilde{\chi}^{p'}_{j'} \leftarrow \widetilde{\chi}^{p'}_{j'} - 1$ for all $p' \in \{p_1,\dots,p_2\}$, and reassigning to the additional charger of terminal $j$ the charging activities of $\min\left\{\rho, \sum_{r \in \mathcal{R}(j')} \widetilde{w}^{p't}_{rj'}\right\}$ BEBs that relied on terminal $j'$ during each interval $t \in \mathcal{T}$ and period $p' \in \{p_1,\dots,p_2\}$ in the original solution. These steps can be repeated until no charger remains in locations dominated by terminals with nonzero residual capacity, i.e., until $\widetilde{\chi}^{p}_{j'} = 0$ $\forall j' \in \mathcal{J}(j)$ if $\widetilde{\chi}^{p}_{j} < \widetilde{\chi}^{\text{UB}}_{j}$ for each period $p \in \mathcal{P}$ and each terminal $j \in \mathcal{J}$. This process terminates after a finite number of steps since the dominance relation of definition \ref{def:bound_dominated_terminals} is a strict partial order.

From there, we can assume that terminal $j$ has reached its maximum hosting capacity $\widetilde{\chi}^{\text{UB}}_{j}$ before considering installing chargers at a dominated location $j' \in \mathcal{J}(j)$. The number of vehicles required to satisfy the service requirements \eqref{model:operations_detailed:service_level} on bus routes connected to terminal $j$ is at most $\max_{p \in \mathcal{P}, t \in \mathcal{T}} \sum_{r \in \mathcal{R}(j) } d_r^{pt}$. The demand for on-route chargers that cannot be satisfied by terminal $j$ is thus upper bounded by $\max\left\{0, \max_{p \in \mathcal{P}, t \in \mathcal{T}} \sum_{r \in \mathcal{R}(j) } d_r^{pt} - \rho \widetilde{\chi}^{\text{UB}}_{j} \right\}$. Since dominated terminals are only connected to a subset of the routes $\mathcal{R}(j)$, this upper bound is in particular valid for the residual demand that can be supplied by location $j' \in \mathcal{J}(j)$. The right-hand side of inequality \eqref{eq:bound_dominated_terminals} is the smallest integer number of chargers that suffices to satisfy the charger capacity constraints \eqref{model:operations_detailed:terminal_capacity} at terminal $j'$ given this upper bound on charging demand. \hfill
\halmos\end{proof}

\section{Illustration of the operational model} \label{ec:Example_operational_model}
Consider an operational horizon of $T=6$ intervals, one route, one BEB type, and one depot, i.e., $\mathcal{R} = \{r\}$, $\mathcal{B} = \{b\}$, and $\mathcal{I}=\{i\}$. Fix a period $p \in \mathcal{P}$ and assume that the fleet is composed of $\bar{\eta}^p_{rb}=3$ depot BEBs with capacity $s_b = 3$, and that the depot is equipped with $\bar{\chi}^p_{i}=2$ chargers. Starting from state $s=2$, one interval suffices to fully recharge the bus, i.e., $\kappa_{rbi2}=1$, compared to two intervals for states $s\in \{1,0\}$, i.e., $\kappa_{rbi1} = \kappa_{rbi0} =~2$. The service level requirements are $d^{p}_r=(d^{p0}_r, d^{p1}_r, d^{p2}_r, d^{p3}_r, d^{p4}_r, d^{p5}_r) = (2,3,2,1,1,1)$. 

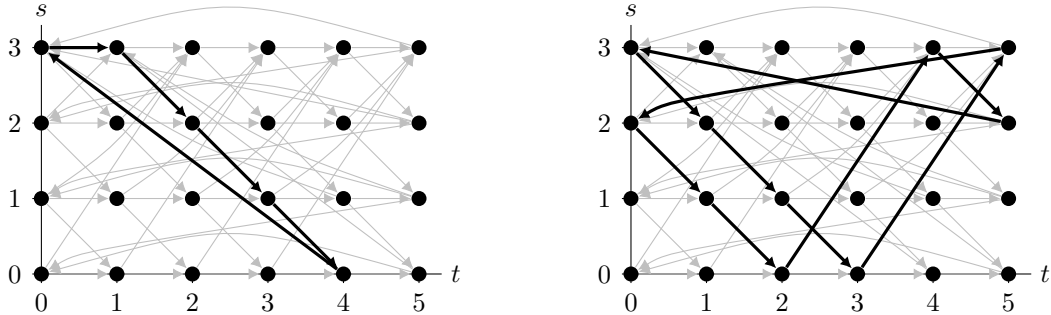
\begin{figure}[htbp]
  \centering
  \caption{Solution consisting of two cyclic schedules, with periods of one day (left) and two days (right)}
  \begin{subfigure}[b]{0.40\textwidth}
    \centering
    \resizebox{\linewidth}{!}{
    \begin{tikzpicture}[>={Latex[width=1.5mm,length=1.5mm]}, node distance=1.5cm]
  
  \draw[-] (0,0) -- (5.3,0);
  \draw[-] (0,0) -- (0,3.3);
  
  \node[right] at (5.3,0) {$t$};
  \node[above] at (0,3.3) {$s$};

  \foreach \x in {0,1,...,5} {
    \draw (\x,0) -- (\x,-0.12); 
    \node[below=4pt] at (\x,0) {\x};
  }
  \foreach \y in {0,1,...,3} {
    \draw (0,\y) -- (-0.12,\y); 
    \node[left=4pt] at (0,\y) {\y};
  }

  \foreach \x in {0,1,...,5} {
    \foreach \y in {0,1,...,3} {
      \node[circle, fill, inner sep=2pt] (\x\y) at (\x,\y) {};
    }
  }

  \draw[->, lightgray] (00) -- (23); 
  \draw[->, lightgray] (01) .. controls (1, 1.7) .. (23);
  \draw[->, lightgray] (02) -- (13);
  \draw[->, lightgray] (10) -- (33);
  \draw[->, lightgray] (11) .. controls (2, 1.7) .. (33);
  \draw[->, lightgray] (12) -- (23);
  \draw[->, lightgray] (20) -- (43);
  \draw[->, lightgray] (21) .. controls (3, 1.7) .. (43);
  \draw[->, lightgray] (22) -- (33);
  \draw[->, lightgray] (30) -- (53);
  \draw[->, lightgray] (31) .. controls (4, 1.7) .. (53);
  \draw[->, lightgray] (32) -- (43);
  \draw[->, lightgray] (40) -- (03);
  \draw[->, lightgray] (41) .. controls (2, 1.7) .. (03);
  \draw[->, lightgray] (42) -- (53);
  \draw[->, lightgray] (50) -- (13);
  \draw[->, lightgray] (51) .. controls (3, 1.7) .. (13);
  \draw[->, lightgray] (52) -- (03);

    \draw[->, lightgray] (01) -- (10);
    \draw[->, lightgray] (11) -- (20);
    \draw[->, lightgray] (21) -- (30);
    \draw[->, lightgray] (41) -- (50);
    \draw[->, lightgray] (02) -- (11);
    \draw[->, lightgray] (12) -- (21);
    \draw[->, lightgray] (32) -- (41);
    \draw[->, lightgray] (42) -- (51);
    \draw[->, lightgray] (03) -- (12);
    \draw[->, lightgray] (23) -- (32);
    \draw[->, lightgray] (33) -- (42);
    \draw[->, lightgray] (43) -- (52);
    \draw[->, lightgray] (51) .. controls (0.5,0.3) .. (00);
    \draw[->, lightgray] (52) .. controls (0.5,1.3) .. (01);
    \draw[->, lightgray] (53) .. controls (0.5,2.3) .. (02);

    \draw[->, lightgray] (00) -- (10);
    \draw[->, lightgray] (10) -- (20);
    \draw[->, lightgray] (20) -- (30);
    \draw[->, lightgray] (30) -- (40);
    \draw[->, lightgray] (40) -- (50);
    \draw[->, lightgray] (01) -- (11);
    \draw[->, lightgray] (11) -- (21);
    \draw[->, lightgray] (21) -- (31);
    \draw[->, lightgray] (31) -- (41);
    \draw[->, lightgray] (41) -- (51);
    \draw[->, lightgray] (02) -- (12);
    \draw[->, lightgray] (12) -- (22);
    \draw[->, lightgray] (22) -- (32);
    \draw[->, lightgray] (32) -- (42);
    \draw[->, lightgray] (42) -- (52);
    \draw[->, lightgray] (13) -- (23);
    \draw[->, lightgray] (23) -- (33);
    \draw[->, lightgray] (33) -- (43);
    \draw[->, lightgray] (43) -- (53);
    \draw[->, lightgray] (50) .. controls (2.5,0.7) .. (00);
    \draw[->, lightgray] (51) .. controls (2.5,1.7) .. (01);
    \draw[->, lightgray] (52) .. controls (2.5,2.7) .. (02);
    \draw[->, lightgray] (53) .. controls (2.5,3.7) .. (03);

    \draw[->, black, very thick] (03) -- (13);
    \draw[->, black, very thick] (13) -- (22);
    \draw[->, black, very thick] (22) -- (31);
    \draw[->, black, very thick] (31) -- (40);
    \draw[->, black, very thick] (40) -- (03);
\end{tikzpicture}
    }
  \end{subfigure}
  \hspace{1cm}
  \begin{subfigure}[b]{0.40\textwidth}
    \centering
    \resizebox{\linewidth}{!}{
    \begin{tikzpicture}[>={Latex[width=1.5mm,length=1.5mm]}, node distance=1.5cm]
  
  \draw[-] (0,0) -- (5.3,0);
  \draw[-] (0,0) -- (0,3.3);
  
  \node[right] at (5.3,0) {$t$};
  \node[above] at (0,3.3) {$s$};

  \foreach \x in {0,1,...,5} {
    \draw (\x,0) -- (\x,-0.12); 
    \node[below=4pt] at (\x,0) {\x};
  }
  \foreach \y in {0,1,...,3} {
    \draw (0,\y) -- (-0.12,\y); 
    \node[left=4pt] at (0,\y) {\y};
  }

  \foreach \x in {0,1,...,5} {
    \foreach \y in {0,1,...,3} {
      \node[circle, fill, inner sep=2pt] (\x\y) at (\x,\y) {};
    }
  }

  \draw[->, lightgray] (00) -- (23); 
  \draw[->, lightgray] (01) .. controls (1, 1.7) .. (23);
  \draw[->, lightgray] (02) -- (13);
  \draw[->, lightgray] (10) -- (33);
  \draw[->, lightgray] (11) .. controls (2, 1.7) .. (33);
  \draw[->, lightgray] (12) -- (23);
  \draw[->, lightgray] (20) -- (43);
  \draw[->, lightgray] (21) .. controls (3, 1.7) .. (43);
  \draw[->, lightgray] (22) -- (33);
  \draw[->, lightgray] (30) -- (53);
  \draw[->, lightgray] (31) .. controls (4, 1.7) .. (53);
  \draw[->, lightgray] (32) -- (43);
  \draw[->, lightgray] (40) -- (03);
  \draw[->, lightgray] (41) .. controls (2, 1.7) .. (03);
  \draw[->, lightgray] (42) -- (53);
  \draw[->, lightgray] (50) -- (13);
  \draw[->, lightgray] (51) .. controls (3, 1.7) .. (13);
  \draw[->, lightgray] (52) -- (03);

    \draw[->,lightgray] (01) -- (10);
    \draw[->,lightgray] (11) -- (20);
    \draw[->,lightgray] (21) -- (30);
    \draw[->,lightgray] (31) -- (40);
    \draw[->,lightgray] (41) -- (50);
    \draw[->,lightgray] (02) -- (11);
    \draw[->,lightgray] (12) -- (21);
    \draw[->,lightgray] (22) -- (31);
    \draw[->,lightgray] (32) -- (41);
    \draw[->,lightgray] (42) -- (51);
    \draw[->,lightgray] (03) -- (12);
    \draw[->,lightgray] (13) -- (22);
    \draw[->,lightgray] (23) -- (32);
    \draw[->,lightgray] (33) -- (42);
    \draw[->,lightgray] (43) -- (52);
    \draw[->,lightgray] (51) .. controls (0.5,0.3) .. (00);
    \draw[->,lightgray] (52) .. controls (0.5,1.3) .. (01);
    \draw[->,lightgray] (53) .. controls (0.5,2.3) .. (02);

    \draw[->,lightgray] (00) -- (10);
    \draw[->,lightgray] (10) -- (20);
    \draw[->,lightgray] (20) -- (30);
    \draw[->,lightgray] (30) -- (40);
    \draw[->,lightgray] (40) -- (50);
    \draw[->,lightgray] (01) -- (11);
    \draw[->,lightgray] (11) -- (21);
    \draw[->,lightgray] (21) -- (31);
    \draw[->,lightgray] (31) -- (41);
    \draw[->,lightgray] (41) -- (51);
    \draw[->,lightgray] (02) -- (12);
    \draw[->,lightgray] (12) -- (22);
    \draw[->,lightgray] (22) -- (32);
    \draw[->,lightgray] (32) -- (42);
    \draw[->,lightgray] (42) -- (52);
    \draw[->,lightgray] (03) -- (13);
    \draw[->,lightgray] (13) -- (23);
    \draw[->,lightgray] (23) -- (33);
    \draw[->,lightgray] (33) -- (43);
    \draw[->,lightgray] (43) -- (53);
    \draw[->,lightgray] (50) .. controls (2.5,0.7) .. (00);
    \draw[->,lightgray] (51) .. controls (2.5,1.7) .. (01);
    \draw[->,lightgray] (52) .. controls (2.5,2.7) .. (02);
    \draw[->,lightgray] (53) .. controls (2.5,3.7) .. (03);

    \draw[->, black, very thick] (03) -- (12);
    \draw[->, black, very thick] (12) -- (21);
    \draw[->, black, very thick] (21) -- (30);
    \draw[->, black, very thick] (30) -- (53);
    \draw[->, black, very thick] (53) .. controls (0.5,2.3) .. (02);
    \draw[->, black, very thick] (02) -- (11);
    \draw[->, black, very thick] (11) -- (20);
    \draw[->, black, very thick] (20) -- (43);
    \draw[->, black, very thick] (43) -- (52);
    \draw[->, black, very thick] (52) -- (03);
  
\end{tikzpicture}
    }
  \end{subfigure}
  \label{fig:graph_cycles}
\end{figure}

\paragraph{Bus schedule.}
Figure \ref{fig:graph_cycles} depicts two cycles in which the arcs in gray have flow 0, and the arcs in black have flow~1. In the first schedule, a bus is in the fully charged state $s = 3$ at time $t = 0$ and idles for one interval, reaching state $(t, s) = (1, 3)$. It is then in service for three intervals, reaching state $(t, s) = (4, 0)$. Finally, it charges for two intervals, completing the cycle by returning to state $(t, s) = (0, 3)$. Since this cycle spans $T$ intervals, a single bus can repeat the corresponding schedule daily. The second schedule can be interpreted similarly. However, since it spans $2T$ intervals, it must be executed by two buses staggered by $T$ intervals. At time $t=0$, one bus is thus in state $s=2$, whereas the second bus is in state $s=3$. After one day, their roles are inverted. It follows that each arc along the cycle is traversed once per day, and the operations repeat daily at the fleet level.

\paragraph{Solution.} The circulation obtained by summing these two cycles yields a feasible solution for the operational problem \eqref{model:operations_detailed}. It corresponds to setting the service variables $w^{pt}_{rbs}$ such that $(t,s) \in \{(1,3), (2,2), (3,1), (0,3), (1,2), (2,1), (5,3), (0,2), (1,1), (4,3)\}$, the idling variables $v^{pt}_{rbs}$ such that $(t,s) = (0,3)$, and the charging variables $z^{pt}_{rbis}$ such that $(t,s) \in \{(4,0), (2,0), (3,0), (5,2)\}$ to 1, and all the other operational variables to 0. The service level matches the minimum requirements imposed by constraints \eqref{model:operations_detailed:service_level}. The charger usage is 0 in intervals $t \in \{0,1\}$, 1 in interval $t=2$, and 2 in intervals $t \in \{3,4,5\}$, which respects the charger capacity constraints \eqref{model:operations_detailed:depot_capacity} for $\bar{\chi}^p_{i}=2$. As the solution is composed of a set of cycles with unit flows, the flow-balance equations \eqref{model:operations_detailed:depot_flow_full}--\eqref{model:operations_detailed:depot_flow_empty} are satisfied. The active variables contributing to the fleet size constraint \eqref{model:operations_detailed:fleet_depot} are $w^{p0}_{rb3}$, $w^{p0}_{rb2}$, and $v^{p0}_{rb3}$, and the limit of $\bar{\eta}^p_{rb}=3$ vehicles is thus respected. Constraints \eqref{model:operations_detailed:terminal_capacity} and \eqref{model:operations_detailed:fleet_conv_route} are trivially satisfied due to the absence of conventional buses and on-route BEBs in the solution.

\begin{revblock}
\paragraph{Fleet-level vs. bus-level cyclicity.}
This example shows that allowing staggered multi-day rotations can reduce fleet size requirements compared to bus-level cyclicity. Indeed, we can observe from Figure~\ref{fig:graph_cycles} that any cycle that repeats daily contains at most three service arcs in this instance. Because the total service requirement is $\sum_{t \in \mathcal T} d^{pt}_r = 10$, at least four BEBs would thus be needed under bus-level cyclicity, compared to the three BEBs used in the schedule above.
\end{revblock}

\section{Review of Benders cut selection methods}\label{ec:benders_theory}
\begin{revblock}
This section reviews the main Benders cut selection methods from the literature and proves Proposition~\ref{prop3:closest_deepest}. We adopt a generic notation to keep the development independent of the BFEP. Consider the two-stage problem:
\begingroup
\setspacing
\allowdisplaybreaks
\begin{align}
\label{model_proof:P} &\min_{x \in \mathcal{X} \subseteq \mathbb{R}^{n_1}_+} \quad  f^{\top}x + \mathcal{Q}(x),
\end{align}
\endgroup
with subproblem value function:
\begingroup
\setspacing
\allowdisplaybreaks
\begin{subequations}
\label{model_proof:SP}
\begin{align}
\label{model_proof:SP:obj} \mathcal{Q}(x) := &\min_{y \in \mathbb{R}^{n}_+} \quad  c^{\top}y \\ 
\label{model_proof:SP:con}  &\text{s.t. }  Ay \geq b-Bx,
\end{align}
\end{subequations}
\endgroup
where $c \geq 0$. In epigraphic form, \eqref{model_proof:P} reads:
\begingroup
\setspacing
\allowdisplaybreaks
\begin{subequations}
\label{model_proof:EP}
\begin{align}
\label{model_proof:EP:obj} &\min_{x \in \mathcal{X}, \ \theta \in \mathbb{R}_+} \quad  f^{\top}x + \theta \\ 
\label{model_proof:EP:con}  &\text{s.t. }  \theta \geq \mathcal{Q}(x).
\end{align}
\end{subequations}
\endgroup

\paragraph{Separation and Benders reformulation.} A solution $(x', \theta')$ violates constraint \eqref{model_proof:EP:con} if and only if the feasibility problem $\{y \in \mathbb{R}^n_+ : Ay \geq b-Bx', \ c^{\top}y \leq \theta'\}$ is empty, which by LP duality is equivalent to the separation problem:
\begingroup
\setspacing
\allowdisplaybreaks
\begin{align}
\label{model_proof:EDSP} &\max_{(\pi, \pi_0) \in \Pi} \quad \pi^{\top}(b-Bx') - \pi_0\theta',
\end{align}
\endgroup
being unbounded, where the cone $\Pi := \{(\pi, \pi_0) \in \mathbb{R}^{m}_+ \times \mathbb{R}_+ : A^{\top}\pi - \pi_0 c \leq 0 \}$ contains the origin, so that \eqref{model_proof:EDSP} is always feasible. Each $(\pi, \pi_0) \in \Pi$ yields a valid Benders cut $\pi^{\top}(b-Bx) - \pi_0\theta \leq 0$, and \eqref{model_proof:P} admits the equivalent reformulation:
\begingroup
\setspacing
\allowdisplaybreaks
\begin{subequations}
\label{model_proof:EMP}
\begin{align}
\label{model_proof:EMP:obj} &\min_{x \in \mathcal{X}, \ \theta \in \mathbb{R}_+} \quad f^{\top}x + \theta \\ 
\label{model_proof:EMP:con} &\text{s.t. } \pi^{\top}(b-Bx) - \theta\pi_0 \leq 0 \quad \forall (\pi, \pi_0) \in \Pi.
\end{align}
\end{subequations}
\endgroup
A cut \eqref{model_proof:EMP:con} with $\pi_0 = 0$ is a feasibility cut $\pi^{\top}(b - Bx) \leq 0$; a cut with $\pi_0 > 0$ is an optimality cut $\theta \geq \pi^{\top}(b - Bx)/\pi_0$. We are now ready to review the main Benders cut selection methods.

\paragraph{Cut selection methods.} The objective and feasible cone of \eqref{model_proof:EDSP} are positive homogeneous in $(\pi, \pi_0)$, so each violated cut corresponds to a ray of $\Pi$, from which a representative can be selected by imposing a normalization constraint: following \citet{fischetti2010note}, one chooses a positive homogeneous function $h: \Pi \rightarrow \mathbb{R}$ and solves:
\begingroup
\setspacing
\allowdisplaybreaks
\begin{subequations}
\label{model_proof:select} 
\begin{align}
\label{model_proof:select:obj} 
&\max_{(\pi, \pi_0) \in \Pi} \quad \pi^{\top}(b-Bx') - \pi_0\theta' \\ 
\label{model_proof:select:con}  &\text{s.t. }  h(\pi, \pi_0) = 1.
\end{align}
\end{subequations}
\endgroup

The \emph{standard} cut \citep{Benders1962} instead solves the dual of subproblem \eqref{model_proof:SP} at $x'$, i.e., problem \eqref{model_proof:EDSP} with $\pi_0$ fixed to one and no normalization, and takes an arbitrary optimal vertex or, if unbounded, an arbitrary extreme ray; it is prone to weak cuts under dual degeneracy. The \emph{Magnanti--Wong} cut \citep{magnanti1981accelerating} selects, among the optimal dual solutions at $x'$, one attaining the best bound at a \emph{core point} (a point in the relative interior of the convex hull of $\mathcal{X}$); it is \emph{Pareto optimal} (nondominated) but does not separate feasibility cuts. The \emph{MIS} cut \citep{fischetti2010note} uses $h_{\mathrm{MIS}}(\pi, \pi_0) := \pi_0 + \sum_{i \in S(B)} \pi_i$, where $S(B)$ is the row support of $B$, to target a minimal infeasible subsystem of the separation problem \citep{gleeson1990identifying}. The \emph{$\ell_p$ deepest} cut \citep{hosseini2025deepest} uses $h_{\ell_p}(\pi, \pi_0) := \lVert (\pi^{\top}B, \pi_0) \rVert_p$, the $\ell_p$ norm of the coefficient vector of the cut, and returns the violated cut whose hyperplane is farthest from $(x', \theta')$ in the distance induced by the dual norm $\ell_q$, where $1/p + 1/q = 1$. 

\paragraph{Closest and Conforti--Wolsey deepest cuts.} Let $\mathcal{E} := \{(x, \theta) \in \mathbb{R}^{n_1}_+ \times \mathbb{R}_+ : \theta \geq \mathcal{Q}(x)\}$ denote the epigraph of the subproblem value function which, by the separation argument above, contains exactly the solutions violating no Benders cut \eqref{model_proof:EMP:con}. Consider a guiding point $(x^o, \theta^o) \in \mathcal{E}$ and assume that $(x', \theta')$ violates at least one cut. To separate a point from a polyhedron, \citet{conforti2019facet} propose a technique that identifies a face of the polyhedron crossed by the segment joining the point to a fixed point in the relative interior of the polyhedron. Applying this technique to $\mathcal{E}$, with $(x^o, \theta^o)$ in its relative interior, \citet{hosseini2025deepest} show that the resulting separation problem falls within the deepest cut framework and corresponds to the selection problem \eqref{model_proof:select} where $h$ is taken as the Conforti--Wolsey pseudonorm $h_{CW}(\pi, \pi_0) := \pi^{\top}B(x^o - x') + \pi_0(\theta^o - \theta')$.

\citet{seo2022closest} instead selects the violated cut whose hyperplane intersects the half-line $\{(x^o, \theta^o) + \beta(x' - x^o, \theta' - \theta^o) \mid \beta \geq 0\}$ nearest to $(x^o, \theta^o)$; setting the cut to equality at the parameterized point and solving for the intersection parameter $\beta$ yields the closest cut separation problem:
\begingroup
\setspacing
\allowdisplaybreaks
\begin{subequations}
\label{model_proof:CC}
\begin{align}\label{model_proof:CC:obj} &\min_{(\pi, \pi_0) \in \Pi} \quad  \beta = \frac{-\pi^{\top}(b-Bx^o) + \theta^o\pi_0}{\pi^{\top}B(x^o - x') + \pi_0(\theta^o - \theta')} \\ 
\label{model_proof:CC:con}  &\text{s.t. } \pi^{\top}B(x^o - x') + \pi_0(\theta^o - \theta') > 0.
\end{align}
\end{subequations}
\endgroup
Since $(x^o, \theta^o)$ violates no cut, the numerator of \eqref{model_proof:CC:obj} is nonnegative and equals the denominator minus the violation $\pi^{\top}(b-Bx') - \pi_0 \theta'$; hence every violated cut is feasible for \eqref{model_proof:CC} with $\beta < 1$, every non-violated feasible solution has $\beta \geq 1$, and problem \eqref{model_proof:CC} returns a violated cut. With this setup in place, the equivalence of the closest cut method of \cite{seo2022closest} and the Conforti–Wolsey deepest cut method of \cite{hosseini2025deepest} is easily established.

\proof{Proof of Proposition~\ref{prop3:closest_deepest}.}
By maximizing $1-\beta$ instead of minimizing $\beta$, the closest cut selection problem \eqref{model_proof:CC} can equivalently be expressed as:
\begingroup
\setspacing
\allowdisplaybreaks
\begin{subequations}
\label{model_proof:CC2}
\begin{align}
\label{model_proof:CC2:obj} &\max_{(\pi, \pi_0) \in \Pi} \quad 1-\beta = \frac{\pi^{\top}(b-Bx') - \theta'\pi_0 }{\pi^{\top}B(x^o - x') + \pi_0(\theta^o - \theta')} \\ 
\label{model_proof:CC2:con} &\text{s.t. } \pi^{\top}B(x^o - x') + \pi_0(\theta^o - \theta') > 0.
\end{align}
\end{subequations}
\endgroup

Since the objective \eqref{model_proof:CC2:obj} and the left-hand side of the normalization constraint \eqref{model_proof:CC2:con} are both positive homogeneous in the dual variables $(\pi, \pi_0)$, the normalization constraint can be replaced by $\pi^{\top}B(x^o - x') + \pi_0(\theta^o - \theta') = 1$, upon which the objective \eqref{model_proof:CC2:obj} reduces to the violation $\pi^{\top}(b-Bx') - \theta'\pi_0$. The resulting formulation is \eqref{model_proof:select} with $h = h_{CW}$, i.e., the Conforti--Wolsey deepest cut selection problem. \Halmos
\endproof


\paragraph{Optimal line-shifting cuts.}
A closely related cut selection method is that of \citet{glomb2026novel}. Given an upper bound $u$ on the optimal value of \eqref{model_proof:P}, they associate with each cut \eqref{model_proof:EMP:con} its \emph{solution candidate set}, the set of points $x \in \mathcal{X}$ for which some $(x, \theta)$ with $f^{\top}x + \theta \leq u$ satisfies the cut. Their \emph{optimal line-shifting} cut excludes from this set as much as possible of the segment joining the iterate $x' \in \mathcal{X}$ to a point $x^o \in \operatorname{conv}(\mathcal{X})$. This criterion is precisely the geometric idea of the closest cuts, and the two separation problems in fact coincide. A point $x$ with $f^{\top}x \leq u$ belongs to the candidate set of a given cut if and only if $(x, u - f^{\top}x)$ satisfies it, and the separation problem in \citet[Problem~12]{glomb2026novel} is, in our notation, the selection problem \eqref{model_proof:select} with $h = h_{CW}$ for the separated point $(x', u - f^{\top}x')$ and the guiding point $(x^o, u - f^{\top}x^o)$. The methods thus differ only in their anchors. Closest cuts separate the master solution $(x', \theta')$ and admit any guiding point in $\mathcal{E}$, whereas both optimal line-shifting anchors are tied to a single scalar bound on the value of the full problem. In a Benders reformulation approximating several value functions, as the BFEP master problem does through the disaggregated variables $\Theta_p$ and $\theta_{pr}$, this anchoring requires allocating the bound across periods and routes. This design choice is not covered by \cite{glomb2026novel}, and we therefore exclude the optimal line-shifting from the computational comparison of \ref{ec:Benders_cuts}.
\end{revblock}

\section{Implementation of monotone cuts}\label{ec:sec:Monotone_cuts_implementation} The monotone cuts \eqref{LBBD_MP:opt_cut}, \eqref{feas_cut_sr}, \eqref{feas_cut_mr}, and \eqref{feas_cut_gen} are active when a componentwise inequality on integer vectors is respected. For a decision vector $\nu \in \mathbb{Z}^m_+$ and a fixed solution $\nu' \in \mathbb{Z}^m_+$, the condition $\mathbf{1}\{\nu \preceq \nu'\}=0$ holds if and only if the disjunction $\lor_{l=1}^m(\nu_l \geq \nu'_l+1)$ is satisfied or, equivalently, if $\sum_{l=1}^m\mathbf{1}\{\nu_l \geq \nu'_l+1\}\geq 1$. Adding such disjunctive constraints to the relaxed master problem requires keeping track of each clause using a binary variable \citep{balas1979disjunctive}. In our implementation, for each period $p \in \mathcal{P}$ and each component $x_{pl}$ of the vector $x_p$, we maintain a pool $\mathcal{A}_{pl}$ of indicators $a_{plk} \in \{0,1\}$, each bounded by a constraint $(k{+}1)a_{plk} \leq x_{pl}$ enforcing that $a_{plk} \leq \mathbf{1}\{x_{pl} \geq k{+}1\}$. Initially, each pool $\mathcal{A}_{pl}$ is empty, and $a_{plk}$ is generated when $x'_{pl}=k$ appears as the threshold value in a monotone cut. A feasibility cut \eqref{LBBD_MP:feas_cut} is then reformulated as $\sum_{l = 1}^n a_{pl(x'_{pl})} \geq 1$, and an optimality cut \eqref{LBBD_MP:opt_cut} becomes $\Theta_p \geq \mathcal{Q}_p(x'_p)(1 - \sum_{l = 1}^n a_{pl(x'_{pl})})$. Similarly, a pool of binary variables is maintained for each sum of decision variables appearing in the feasibility cuts \eqref{feas_cut_sr} and \eqref{feas_cut_mr}.

In contrast with standard implementations \citep[e.g.,][]{hooker1999mixed, liu2024generalized}, where a new collection of $n$ binary variables is generated each time a monotone cut is added to the master problem, we allow the indicators $a_{plk}$ to appear in all the cuts where the same threshold on $x_{pl}$ is encountered. In practice, this makes a significant difference in the number of variables added to the master problem. Indeed, monotone cuts are mostly needed in late iterations of the algorithm, where near-optimal solutions that only differ in a few components are sequentially visited and bounded by monotone optimality cuts. By reusing previously generated binary indicators, such cuts can be generated by introducing a small fraction of the indicators required by the standard implementation. 

\section{Instances and model parameters}\label{ec:Experiments:instances}
\begin{table}[htbp]
\centering
\caption{Characteristics of complete network instances}
\resizebox{0.7\textwidth}{!}{  
\renewcommand{\arraystretch}{0.75}
\begin{tabular}{lccccccc}
        \toprule
        City & \# Depots & \# Terminals & \# Routes & \# Buses & Service rate (\%) \\
        \midrule
        Atlanta      & \p5  &  115 &  110 & \p441  & 68.5 \\
        Boston       & \p9  &  214 &  166 &  1286  & 53.1 \\
        Chicago      & \p7  &  235 &  126 &  1495  & 51.3 \\
        Dallas       & \p3  &  152 &  149 & \p623  & 55.8 \\
        Detroit      & \p4  & \p78 & \p42 & \p214  & 60.6 \\
        Houston      & \p8  &  135 &  117 & \p941  & 47.3 \\
        Las Vegas    & \p2  & \p75 & \p39 & \p280  & 71.4 \\
        Los Angeles  &  12  &  188 &  140 &  1817  & 51.9 \\
        \bottomrule
    \end{tabular}
}  
\label{tab:instances}
\end{table}
We based all our model parameters on the baseline scenario presented in \citet{johnson2020financial}. We consider 70 kW AC plug-in depot chargers, with a cost of $\$60.05$k, including lifetime maintenance. The on-route chargers represent 325 kW DC chargers with an acquisition cost of $\$877.59$k. Also, we consider $|\mathcal{B}|=2$ models of depot BEBs, respectively with $s_1 = 6$ and $s_2 = 12$ hours of operational capacity, all requiring $\kappa_{rbis} = \rev{\lceil\frac{s_b-s}{2}\rceil}$ time intervals \begin{revblock}to travel back and forth between their route and any depot and fully charge from state $s \in \mathcal{S}^z_b$\end{revblock}. We assume that each terminal location has an installation capacity of $\widetilde{\chi}_j^{\mathrm{UB}}=2$ on-route chargers, and that each installed on-route charger can serve up to $\rho=8$ on-route BEBs per hourly time interval. \begin{revblock} The value $\rho=8$ follows the baseline fast-charging scenario of \citet{johnson2020financial}, which assumes one route fast charger for up to eight vehicles. Equivalently, with our one-hour discretization, this corresponds to an average allocation of $60/8=7.5$ minutes of charger capacity per BEB. For a 325 kW fast charger, this represents approximately 40.6 kWh of charger output per BEB, or enough energy for roughly 20 miles of service at a reference consumption rate of 2 kWh/mile. \end{revblock} The unit cost of short-range depot BEBs is set to $\$943$k per unit, compared to $\$1093$k for long-range depot BEBs and on-route BEBS. Salvage revenues from the retirement of conventional buses are ignored. The operating costs of conventional buses were estimated to be \$50 per hour of service, compared to an average of \$29 and \$31 for depot and on-route BEBs, respectively, after factoring charging and deadhead costs. Deadhead costs depend on the distance from each route to the selected charging location in the case of depot BEBs, and are assumed to be negligible for on-route BEBs since fast chargers are located at bus terminals. In addition, a fixed yearly maintenance cost of \$10k per unit is applied to the conventional buses. The operational costs are computed based on 250 yearly repetitions of the simulated weekday of service, and a yearly discount factor of 4\% ($\gamma = 0.96$) is used. In all the experiments, we assume the initial fleet to be composed exclusively of conventional buses, which must all be retired by the end of the planning horizon ($\widehat{\eta}^{\text{UB}}_P = 0$). No explicit BEB acquisition target is imposed ($\eta_p^{\text{LB}}=0$ for each $p \in \mathcal{P}$). For each instance, we estimate the minimum investment budget needed to satisfy the electrification targets by solving the LP relaxation of the problem with an alternative objective. This budget is then multiplied by 1.5 (2.5 for instances without on-route charging), and equally divided over the $P$ investment periods to obtain the yearly budget parameters $I^{\text{UB}}_p$. Table \ref{tab:instances} summarizes the characteristics of each city. The number of buses refers to the initial state of the system (year $p=0$), defined as the sum of peak demands across routes $\sum_{r \in \mathcal{R}} \max_{t \in \mathcal{T}} d^{0t}_r$. The service rate is the ratio of average service requirements $\frac{1}{T}\sum_{r \in \mathcal{R}}\sum_{t \in \mathcal{T}}d^{0t}_r$ to fleet size in year $0$.

\section{Empirical comparison of Benders cut selection methods} \label{ec:Benders_cuts}
In Table \ref{tab:cuts_comparison_gm}, we evaluate our logic-based Benders decomposition algorithm with the standard Benders cuts \citep{Benders1962}, the MIS cuts \citep{fischetti2010note}, the MW cuts \citep{magnanti1981accelerating}, the closest cuts \citep{seo2022closest} (equivalently, the Conforti--Wolsey deepest cuts, by Proposition \ref{prop3:closest_deepest}), and the $\ell_1$ deepest cuts \citep{hosseini2025deepest}. We also performed preliminary experiments with the $\ell_2$ deepest cuts, but solving the resulting nonlinear cut selection problems proved to be prohibitively expensive. The experiments are performed on the instances presented in the ablation study of Section \ref{sec:Experiments:aceclerationLBBD}, and the same metrics are reported. For the MW and closest cuts, we use the procedure proposed by \cite{seo2022closest} to identify a guiding point $(x^o_p, \Theta^o_p) \in \text{conv}(\mathcal{X}_p) \times \mathbb{R}_{+}$ that does not violate any Benders cut. First, we solve the LP relaxation of problem \eqref{model:BFEP} in extensive form, which gives a solution $(x'_p, y'_p)$ for each period $p \in \mathcal{P}$. We then set $\Theta^o_p$ to the objective value $\sum_{r \in \mathcal{R}}c^{y\top}_{pr}y'_{pr}$ of solution $y'_p$ for the operational problem \eqref{model:operations}, similarly define $\theta^o_{pr} = c^{y\top}_{pr}y'_{pr}$ for single-route problems, and take $x^o_p = x'_p + \epsilon \boldsymbol{1}$ for $\epsilon=0.1$. The added perturbation slightly relaxes the operational constraints \eqref{model:operations:linking_constraints}--\eqref{model:operations:route_constraints}, which makes it more likely for $x^o$ to be a core point.

\begin{table}[H]
\centering
\caption{Performance of LBBD with different Benders cut selection methods}
\resizebox{1\textwidth}{!}{  
\setlength{\tabcolsep}{8pt}
\renewcommand{\arraystretch}{0.75}
\begin{tabular}{clrrrrrrrrrr}
\toprule
\multirow{2}{*}{\specialcell{Instances \\ $(|\mathcal{R}|, |\mathcal{P}|)$}} &  \multirow{2}{*}{\specialcell{Cut \ \\ type}}  & \multicolumn{4}{c}{Summary}  & \multicolumn{3}{c}{Cuts} & \multicolumn{3}{c}{Time (\%)} \\
\cmidrule(lr){3-6} \cmidrule(lr){7-9} \cmidrule(lr){10-12}
&  &  \multicolumn{1}{c}{Opt} & \multicolumn{1}{c}{Gap} & \multicolumn{1}{r}{Time} & \multicolumn{1}{r}{Iter} & \multicolumn{1}{c}{BCuts} & \multicolumn{1}{c}{MCuts} & \multicolumn{1}{r}{Ind} & MP& LP & IP \\
\midrule       
      & Standard  &     \textbf{10} &   \textbf{0.0000} &    12.1 &    9.4 &     58.2 &    2.1 &    4.5 &    29.7 &    19.6 &    21.0 \\     
      & MIS       &     \textbf{10} &   \textbf{0.0000} &    17.5 &    7.4 &     40.3 &    1.8 &    3.3 &     9.7 &    63.3 &     8.6 \\     
(3,4) & Closest   &     \textbf{10} &   \textbf{0.0000} &    \textbf{11.1} &    \textbf{5.9} &     33.1 &    \textbf{1.7} &    \textbf{3.1} &    10.5 &    44.7 &    10.2 \\     
      & MW        &     \textbf{10} &   \textbf{0.0000} &    11.2 &    6.1 &     \textbf{30.9} &    1.9 &    4.0 &    11.2 &    37.6 &    15.2 \\     
      & $\ell_1$ deepest  &     \textbf{10} &   \textbf{0.0000} &    16.9 &    7.3 &     39.4 &    1.9 &    3.3 &     8.9 &    62.3 &     9.7 \\     
\midrule      
      & Standard  &     \textbf{10} &   \textbf{0.0000} &   108.1 &   21.0 &    230.1 &    2.4 &    7.2 &    63.1 &    13.5 &    14.6 \\     
      & MIS       &     \textbf{10} &   \textbf{0.0000} &    79.7 &   12.6 &    137.0 &    2.0 &    4.9 &    26.5 &    55.3 &     8.7 \\     
(6,6) & Closest   &     \textbf{10} &   \textbf{0.0000} &    \textbf{56.9} &   \textbf{10.3} &    113.3 &    1.9 &    4.9 &    29.7 &    38.3 &    13.3 \\     
      & MW        &     \textbf{10} &   \textbf{0.0000} &    61.7 &   11.2 &    \textbf{105.6} &    2.1 &    6.4 &    31.4 &    37.7 &    13.2 \\     
      & $\ell_1$ deepest   &     \textbf{10} &   \textbf{0.0000} &    76.5 &   11.8 &    131.4 &    \textbf{1.7} &    \textbf{3.3} &    25.7 &    55.3 &     9.2 \\     
\midrule     
      & Standard  &     6 &    0.0334 &  2137.6 &   61.3 &    581.8 &    5.5 &   35.5 &    82.1 &     6.3 &     9.8 \\     
      & MIS       &     \textbf{9} &    0.0051 &   601.3 &   27.9 &    296.9 &    \textbf{3.3} &    \textbf{9.4} &    53.5 &    35.3 &     6.8 \\     
(9,6) & Closest   &     \textbf{9} &    \textbf{0.0037} &   \textbf{372.8} &   \textbf{24.0} &    \textbf{247.9} &    3.5 &   13.1 &    54.5 &    29.2 &     7.7 \\     
      & MW        &     8 &    0.0071 &   474.9 &   27.3 &    257.0 &    4.2 &   17.8 &    53.6 &    27.1 &    10.9 \\     
      & $\ell_1$ deepest   &     \textbf{9} &    0.0048 &   567.3 &   27.2 &    286.4 &    4.1 &   14.6 &    52.3 &    38.0 &     5.4 \\     
\bottomrule
\end{tabular}
}  
\label{tab:cuts_comparison_gm}
\end{table}
\vspace{-0.5cm}

The results show that the Benders cut selection method has a significant impact on the overall performance of the algorithm. For $(|\mathcal{R}|, |\mathcal{P}|) = (9,6)$, nine instances are solved to optimality with the closest, MIS, and $\ell_1$ deepest cuts, compared to eight with the MW cuts, and six with the standard cuts. The closest cuts provide the best performance across all groups of instances. They divide the average number of iterations and generated cuts by more than two compared to standard cuts, and provide a speedup factor exceeding two orders of magnitude for some instances. Although separating stronger Benders cuts requires more computational effort than standard cuts, solving the master problem is the most computationally expensive step of our algorithm for challenging instances, hence the importance of limiting the number of iterations and generated cuts.

\section{Detailed results on restriction heuristics} \label{ec:heuristics}
\begin{figure}[H]
\caption{Average number of arcs in restriction heuristics models - Complete networks}
\label{fig:arcs_heuristics}
    \centering
    \begin{subfigure}[b]{0.49\textwidth}
        \centering
        \includegraphics[width=\textwidth]{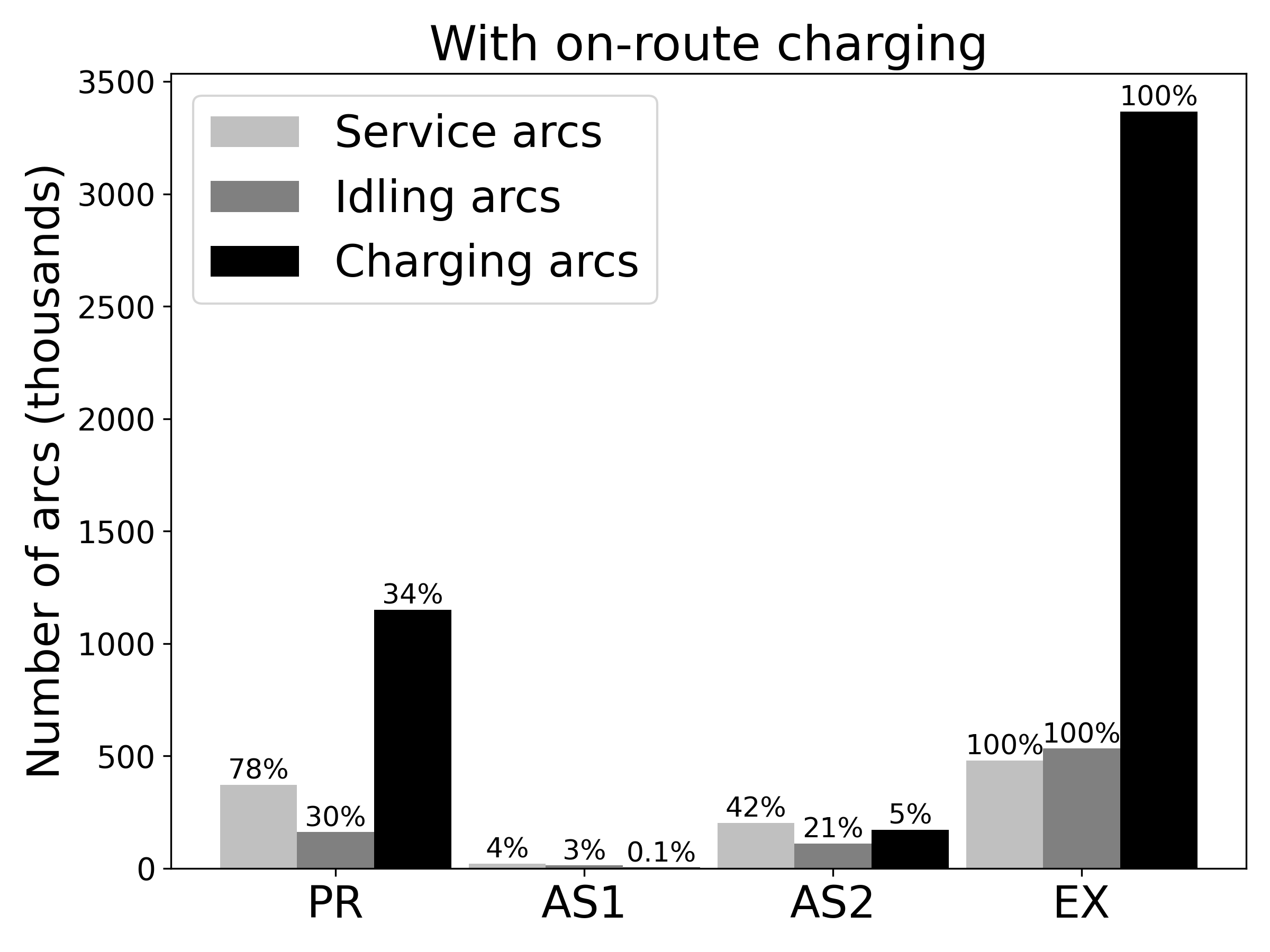}
    \end{subfigure}
    \hfill
    \begin{subfigure}[b]{0.49\textwidth}
        \centering
        \includegraphics[width=\textwidth]{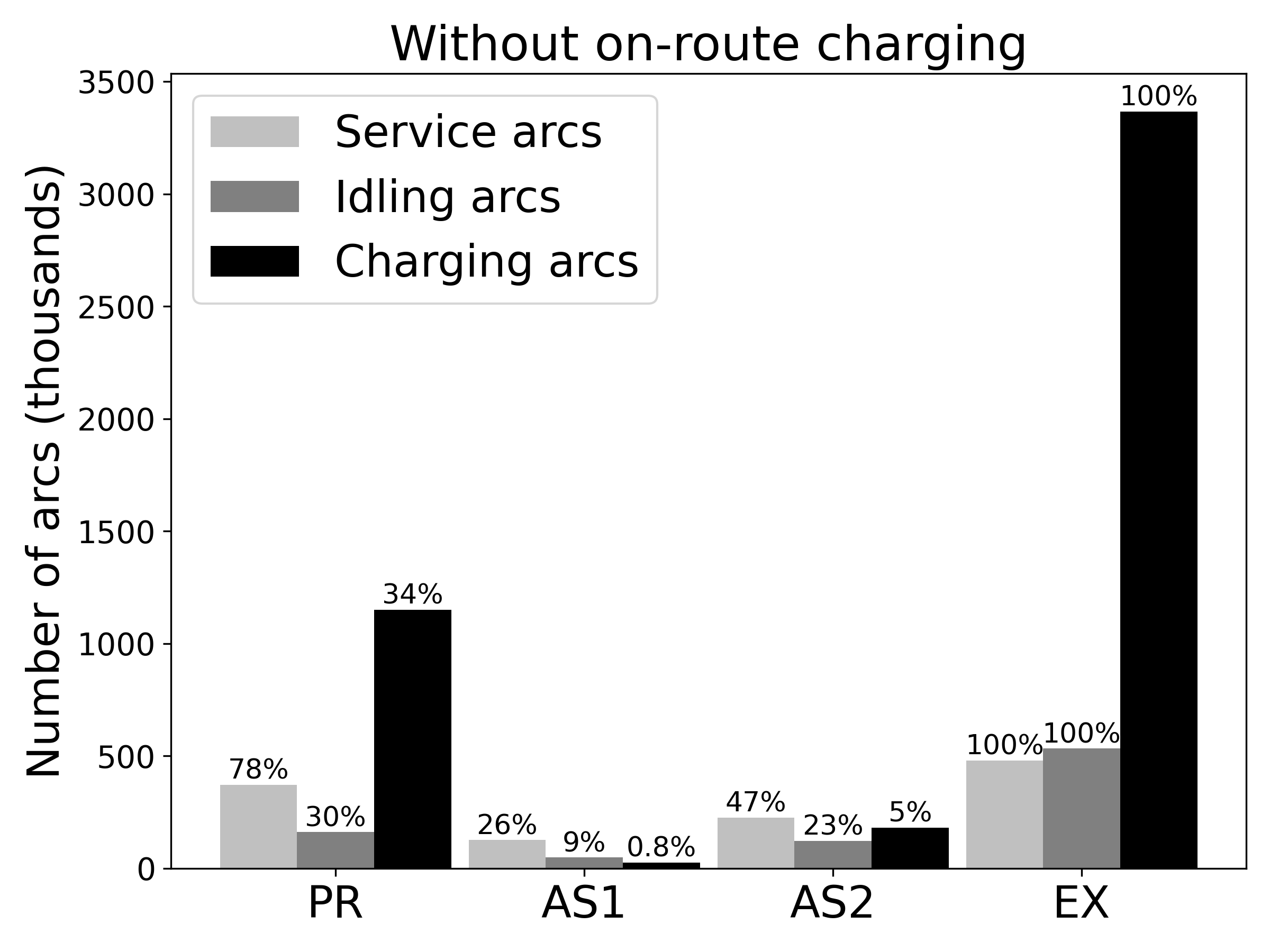}
    \end{subfigure}
    \captionsetup{skip=2pt}
\end{figure}

Figure \ref{fig:arcs_heuristics} presents the total number of arcs in the depot BEB scheduling graphs for the extensive formulation (EX), the policy restriction model (PR), and the first (AS1) and second (AS2) restricted models solved in the arc selection algorithm (steps \ref{alg:arc_selection:IP1} and \ref{alg:arc_selection:IP2} of Algorithm \ref{alg:arc_selection}). The values are averaged over the complete network instances of Section \ref{sec:Experiments:large}. The fraction of each type of arcs (service arcs $w$, idling arcs $v$, charging arcs $z$) retained in each model is also reported.

\begin{figure}[H]
\caption{Average objective bounds of restriction heuristics models - Complete networks}
\label{fig:objective_heuristics}
    \centering
    \begin{subfigure}[b]{0.49\textwidth}
        \centering
        \includegraphics[width=\textwidth]{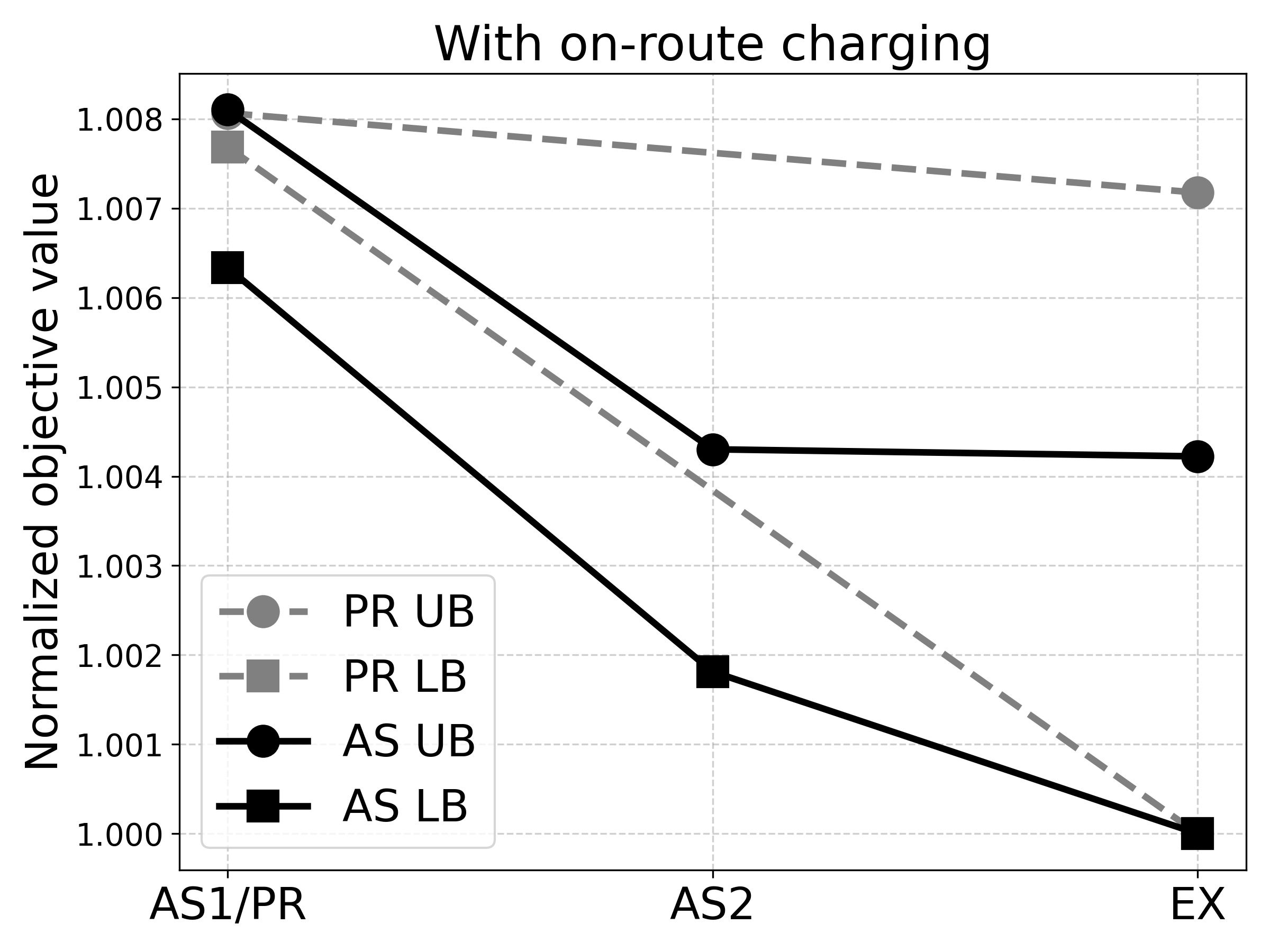}
    \end{subfigure}
    \hfill
    \begin{subfigure}[b]{0.49\textwidth}
        \centering
        \includegraphics[width=\textwidth]{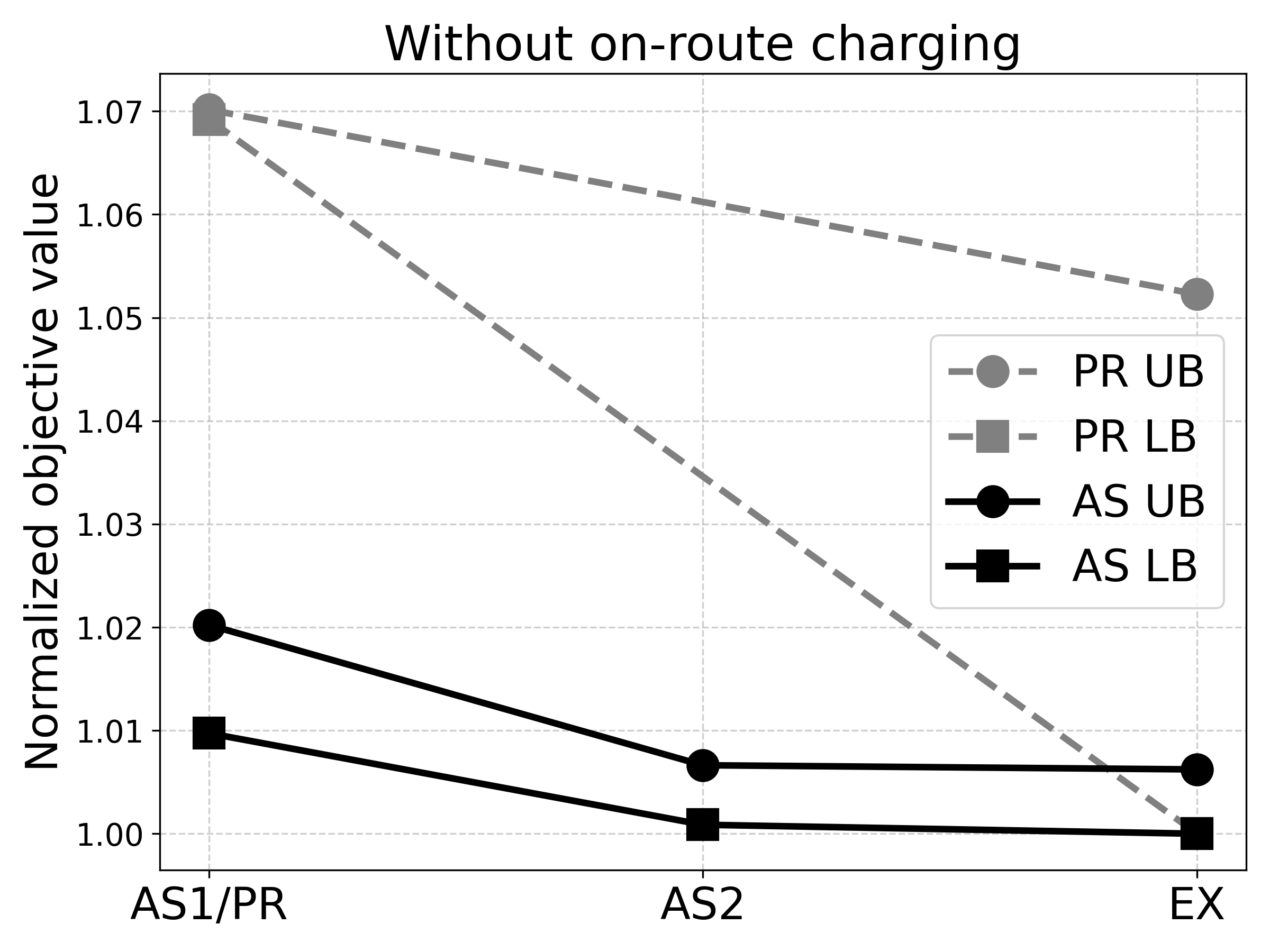}
    \end{subfigure}
    \captionsetup{skip=2pt}
\end{figure}

\vspace{-0.2cm}
For the same instances, Figure \ref{fig:objective_heuristics} presents the average lower and upper bounds on the objective value of each model solved in the PR and AS algorithms. The values are normalized based on the best-known lower bound for each instance. Note that time limits of 1.5 hours and 4.5 hours are respectively given to models AS1 and AS2, whereas 6 hours are given to model PR. In both the AS and PR algorithms, the warm-started EX model is then solved for two hours. 

\vspace{-0.2cm}
Although models AS1 and AS2 are much sparser, they capture better solutions than PR, especially for instances that rely exclusively on depot charging. AS1 retains 0.9\% and 4.6\% of the arcs in the first and second groups of instances, respectively. For AS2, these values are 11.1\% and 12.1\%, respectively. PR, whose a priori selection rules retain the same arcs whether on-route BEBs are available or not, retains 38.4\% of the arcs in the model. In the second group of instances, AS1 provides feasible solutions with optimality gaps of 2\% for the unrestricted model, compared to 7\% for PR. The PR model can consistently be solved to near-optimality, confirming that the a priori restriction rules eliminate the best feasible solutions from the feasible set. By contrast, models AS1 and AS2 are more difficult to solve, but the good upper bounds they provide confirm the existence of high-quality solutions that use only their sparse depot BEB operational graphs. In algorithm PR, the two-hour computing budget allocated to the warm-started model EX suffices to improve the heuristic solution provided by the restricted model. In contrast, the best solution identified by model AS2 is almost never improved in the last phase of algorithm AS. In this case, the main purpose of solving the warm-started model EX is thus to obtain an optimality gap for the original problem. 

\begin{revblock}
\section{Charging assumptions and their impact on fleet sizing}
\label{ec:partial-charging}

Throughout the paper, we assumed that (i) BEBs can be partitioned into depot BEBs and on-route BEBs; and (ii) depot BEBs cannot be partially charged. In this section, we show how these assumptions can be relaxed within our framework, and we discuss their impact on fleet sizing.

\paragraph{Model extension.}
Partial charging can be integrated in the scheduling graphs of model \eqref{model:operations_detailed} by replacing each charging arc $z^{pt}_{rbis}$ by a collection of arcs $\{z^{pt}_{rbiss'}\}_{s'=s+1}^{s_b}$, one per target state of charge $s' \geq s+1$. Similarly, allowing depot BEBs to opportunity-charge at terminals amounts to adding, in parallel to each idling arc $v^{pt}_{rbs}$, a collection of arcs $\{\widehat{w}^{pt}_{rbsj}\}_{j \in \mathcal{J}(r)}$ representing service intervals during which the BEB charges at terminal $j$ to maintain its state of charge. In both cases, the new arcs enter the existing flow-balance constraints \eqref{model:operations_detailed:depot_flow_full}--\eqref{model:operations_detailed:depot_flow_empty} and fleet-size constraints \eqref{model:operations_detailed:fleet_depot}; partial charging arcs additionally enter the depot charger usage constraints \eqref{model:operations_detailed:depot_capacity}, and on-route charging arcs enter the terminal capacity constraints \eqref{model:operations_detailed:terminal_capacity}. Partial charging thus multiplies the number of charging arcs per depot BEB of type $b$ by $(s_b+1)/2$, and allowing depot BEBs to opportunity-charge is equivalent to multiplying the number of idling arcs by $|\mathcal{J}(r)|+1$ for buses assigned to route $r$.

\paragraph{Impact on fleet sizing.}
To quantify the impact of partial charging on fleet sizing, we solve the single-route fleet-minimization subproblem \eqref{model:min_fleet}, restricted to a single depot BEB type, under both full charging only and partial charging enabled. Summing across routes yields a lower bound on the number of depot BEBs required for citywide electrification using a homogeneous fleet. We also report the on-route BEB lower bound, which equals $\sum_{r \in \mathcal{R}} \max_{t\in \mathcal{T}} d^{t}_{r}$ as on-route BEBs can be in service without interruption. Table~\ref{tab:fleet-bounds} compares these bounds with the year-10 fully electrified fleet obtained for each city in Section~\ref{sec:Experiments:large} by solving the BFEP with the AS heuristic.

\begin{table}[H]
\rev{
\centering
\caption{\rev{Lower bounds on fleet size by bus type and charging policy, compared with BFEP solutions}}
\resizebox{\textwidth}{!}{
\renewcommand{\arraystretch}{0.75}
\begin{tabular}{l c cc cc cccc ccc}
\toprule
& \multicolumn{5}{c}{Lower bounds on fleet size for feasibility} & \multicolumn{7}{c}{Fleet size from BFEP solutions} \\
\cmidrule(lr){2-6} \cmidrule(lr){7-13}
& \multirow{2}{*}{\shortstack{\\ On-route \\ \\ BEB}} & \multicolumn{2}{c}{Short-range BEB} & \multicolumn{2}{c}{Long-range BEB} & \multicolumn{4}{c}{With on-route charging} & \multicolumn{3}{c}{Without on-route} \\
\cmidrule(lr){3-4} \cmidrule(lr){5-6} \cmidrule(lr){7-10} \cmidrule(lr){11-13}
City & & Full & Partial & Full & Partial & Total & on-route & short & long & Total & short & long \\
\midrule
Atlanta     & \p441 & \p581 & \p558 & \p534 & \p520 &  \p441 & 343 &  \pp97 &  \p1 & \p542 & \p436 &  106 \\
Boston      &  1286 &  1439 &  1399 &  1328 &  1306 &   1286 & 662 &  \p622 &  \p2 &  1340 &  1011 &  329 \\
Chicago     &  1495 &  1639 &  1587 &  1520 &  1503 &   1495 & 788 &  \p690 &   17 &  1522 &  1135 &  387 \\
Dallas      & \p623 & \p744 & \p719 & \p689 & \p658 &  \p623 & 336 &  \p281 &  \p6 & \p689 & \p549 &  140 \\
Detroit     & \p214 & \p263 & \p255 & \p235 & \p226 &  \p214 & 145 &  \pp61 &  \p8 & \p237 & \p152 & \p85 \\
Houston     & \p941 &  1036 &  1015 & \p973 & \p967 &  \p941 & 317 &  \p610 &   14 & \p974 & \p755 &  219 \\
Las Vegas   & \p280 & \p378 & \p357 & \p329 & \p322 &  \p280 & 227 &  \pp30 &   23 & \p331 & \p168 &  163 \\
Los Angeles &  1817 &  1966 &  1925 &  1852 &  1834 &   1817 & 780 &   1024 &   13 &  1857 &  1469 &  388 \\
\bottomrule
\end{tabular}
}
\label{tab:fleet-bounds}
}
\end{table}

We observe that when on-route charging is available, the fleet size from BFEP solutions always attains the on-route BEB lower bound. In our instances, hybrid depot/on-route use thus cannot reduce the fleet further, and may only shift the fleet composition across BEB types. Allowing partial charging, in contrast, reduces the short-range lower bound by 2.9\% across the eight cities, and the long-range lower bound by 1.7\%. When on-route charging is unavailable, the BFEP fleet closely follows the long-range, full-charge lower bound, exceeding it by 0.43\% across the eight cities. We may thus approximate the impact of partial charging on optimal fleet sizing by the 1.7\% reduction in this lower bound. We conclude that the operational restrictions on charging adopted throughout the paper have a negligible impact on fleet sizing when on-route charging is available, since the BFEP fleet already attains the on-route bound. When using only depot charging, allowing partial charging may yield small but noticeable reductions. In this setting, for planning purposes, the full-charge-only restriction thus yields conservative fleet size requirements.

\section{Rolling-horizon implementation of the BFEP}\label{ec:rolling_horizon}

Since the transition is inherently sequential, the model should be deployed in a rolling horizon when parameters are uncertain. At the start of each year $p$, the agency re-solves the problem with its current best estimates of the input parameters, fixes the decisions already implemented in years $1$ to $p-1$, commits only to the year-$p$ decisions, and repeats the following year. Future decisions thus remain revisable as costs, supply conditions, and other parameters become better known.

We illustrate this deployment on the Chicago case study using a stylized model of BEB acquisition-cost uncertainty. Holding the feasible set $\mathcal{X}$ fixed, we scale every BEB objective coefficient in year $p$ by a multiplier $\zeta^p$ following a random walk, with $\zeta^0 = 1$ and $\zeta^p = \zeta^{p-1}\mu^{\,p}\exp(\varepsilon^p)$. The term $\varepsilon^p \stackrel{\text{iid}}{\sim} \mathcal{N}(0,0.05^2)$ is unanticipated year-on-year noise (about $5\%$ annual volatility), while $\mu^{\,p}$ is a structural multiplier representing a market event whose forward impact becomes known once it occurs. We consider two such events, both materializing in year~5: (i) a technology breakthrough (\emph{Scenario~A}), with $\mu = (1,1,1,1,1,0.93,0.93,0.93,0.93,0.93)$, announcing a compounding $7\%$ annual decline from year~6, and (ii) a supply-chain disruption (\emph{Scenario~B}), with $\mu = (1,1,1,1,1.30,0.885,0.87,1,1,1)$, raising prices $30\%$ and $15\%$ above baseline in years~5 and~6 before reverting. Each year the agency observes the realized $\zeta^p$. It extrapolates the current level flat for all future years until the event occurs, then from year~5 on, the structural multipliers are known and it forecasts year-$q$ prices by compounding them from the latest observation, as $\zeta^p\prod_{p'=p+1}^{q}\mu^{\,p'}$. 

For each scenario we sample four price trajectories and compare three plans on each: the \emph{baseline} solution of Figure \ref{fig:chicago_planning}, solved once in year~1 and never revised, the \emph{rolling-horizon} plan solved annually under the information structure above, and the \emph{perfect-foresight} benchmark solved once with the realized trajectory known in advance, which bounds the savings attainable on each path. Figure~\ref{fig:price_scenarios} shows the price paths and Figure~\ref{fig:rh_sol} the resulting fleet-sizing trajectories. Relative to the baseline, rolling horizon defers the year-5 investments in both scenarios: as far as possible in Scenario~A, to exploit the announced decline, and to years~6 and~7 in Scenario~B, once the temporary spike has passed. Table~\ref{tab:simul_costs} reports the average realized cost: rolling horizon saves roughly $3.5\%$ over the baseline, recovering $77\%$ of the perfect-foresight savings in Scenario~A and $99\%$ in Scenario~B.
\end{revblock}

\begin{figure}[H]
\caption{\begin{revblock}Simulated price trajectories\end{revblock}}
    \label{fig:price_scenarios}
    \centering
    \begin{subfigure}{0.49\textwidth}
        \includegraphics[width=\linewidth]{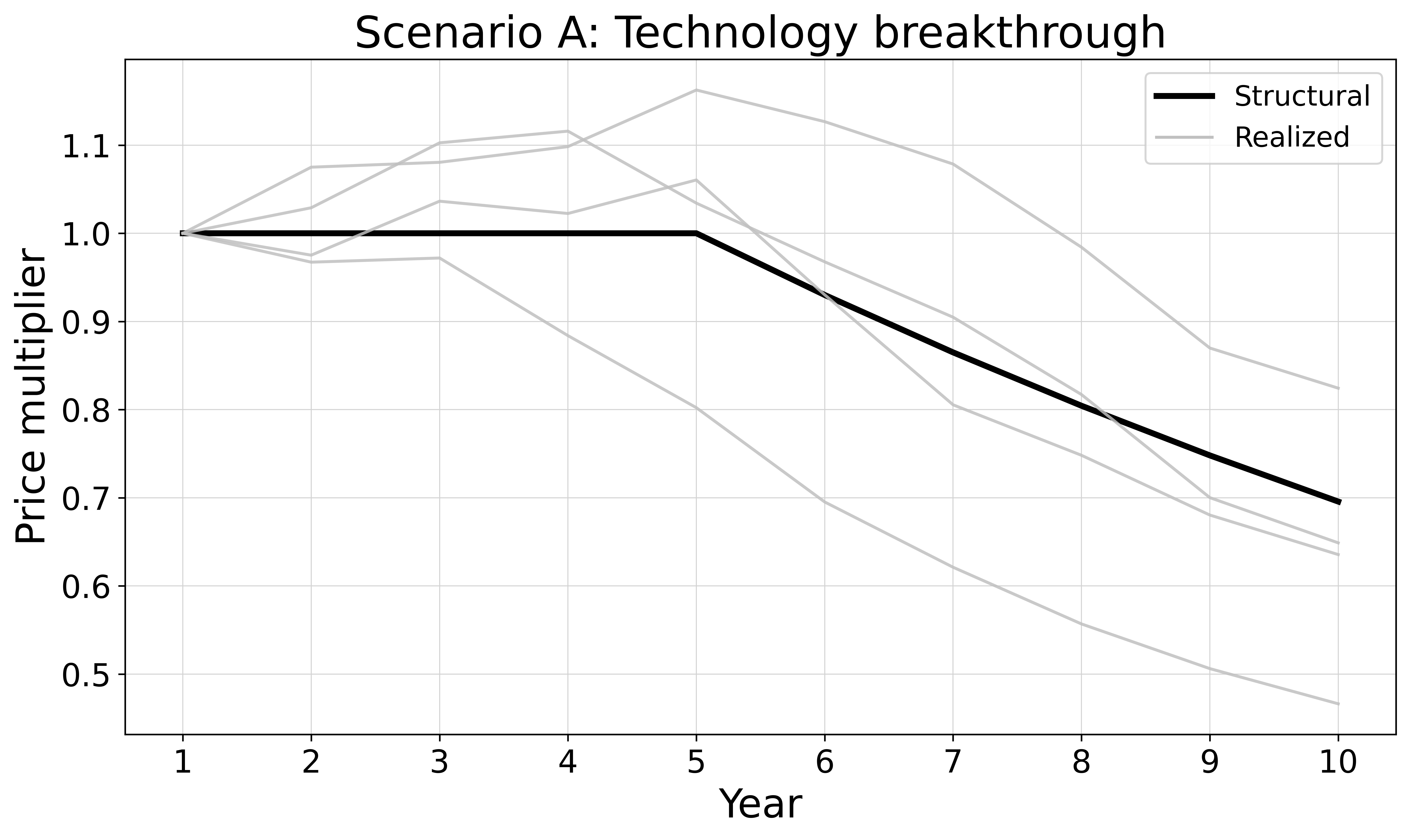}
    \end{subfigure}
    \hfill
    \begin{subfigure}{0.482\textwidth}
        \includegraphics[width=\linewidth]{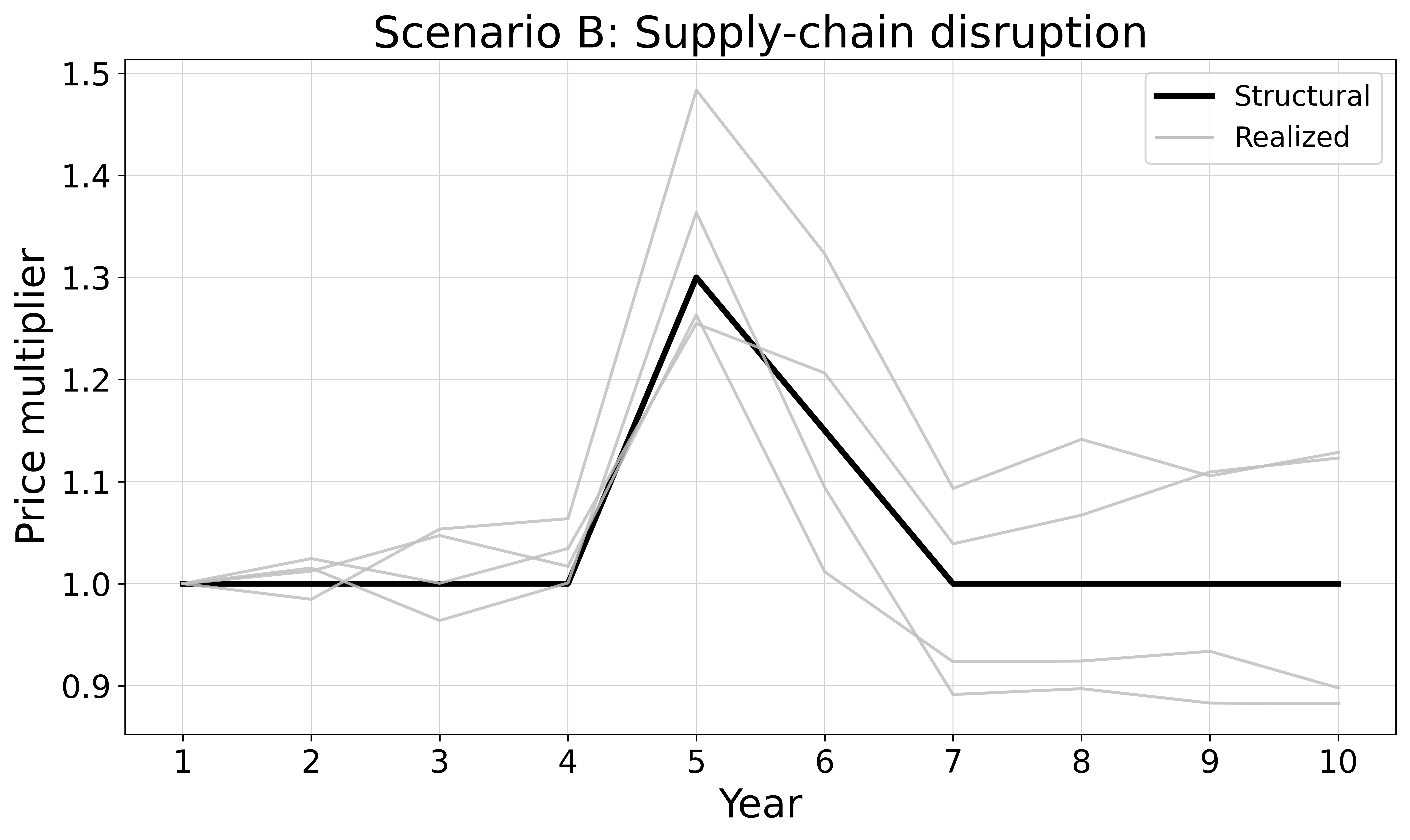}
    \end{subfigure}
    \captionsetup{skip=-8pt}
\end{figure}

\vspace{-1.1cm}

\begin{figure}[H]
\caption{\begin{revblock}Fleet sizing decisions per scenario and price trajectory\end{revblock}}
\label{fig:rh_sol}
    \centering
    \begin{subfigure}{0.49\textwidth}
        \includegraphics[width=\linewidth]{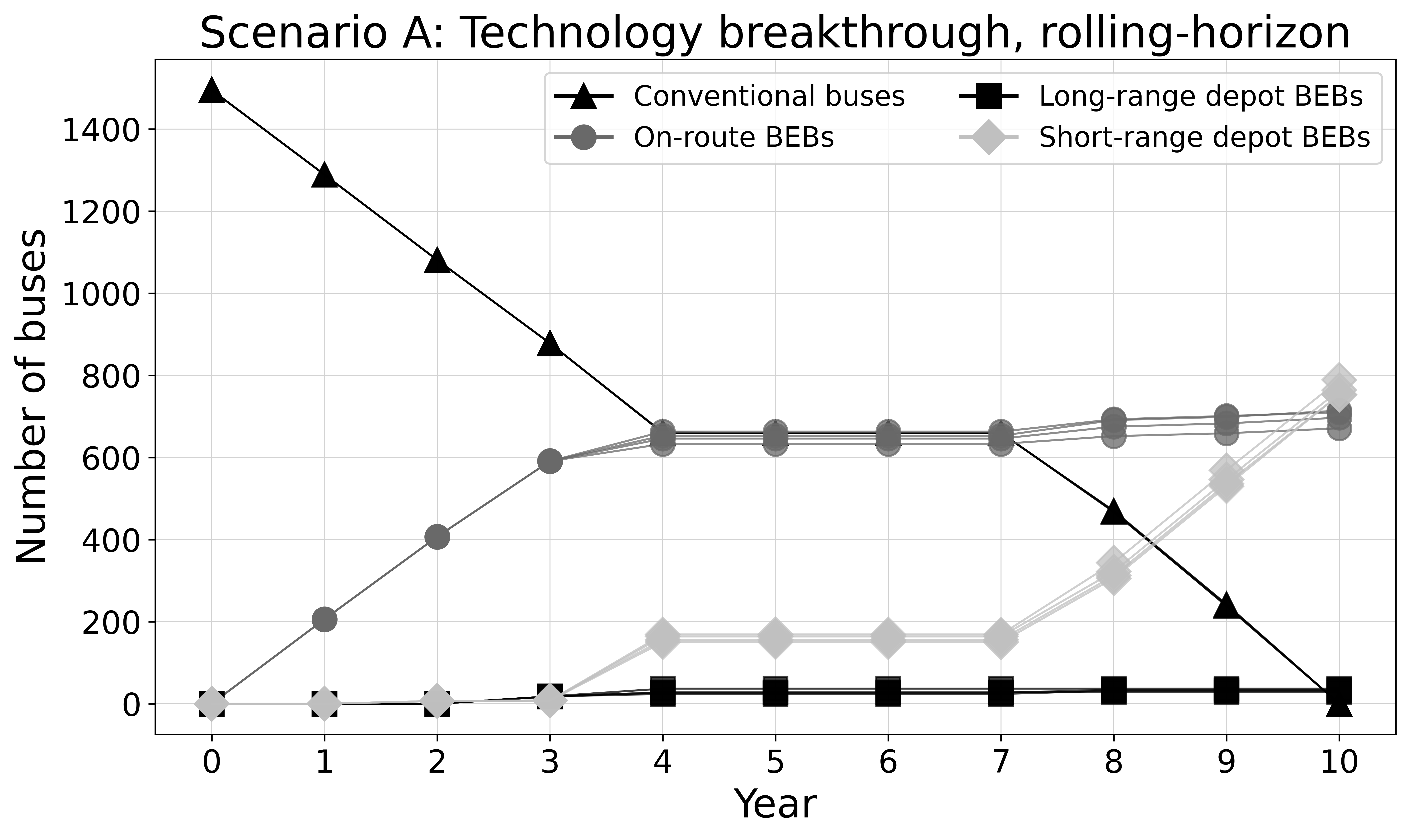}
    \end{subfigure}
    \hfill
    \begin{subfigure}{0.49\textwidth}
        \includegraphics[width=\linewidth]{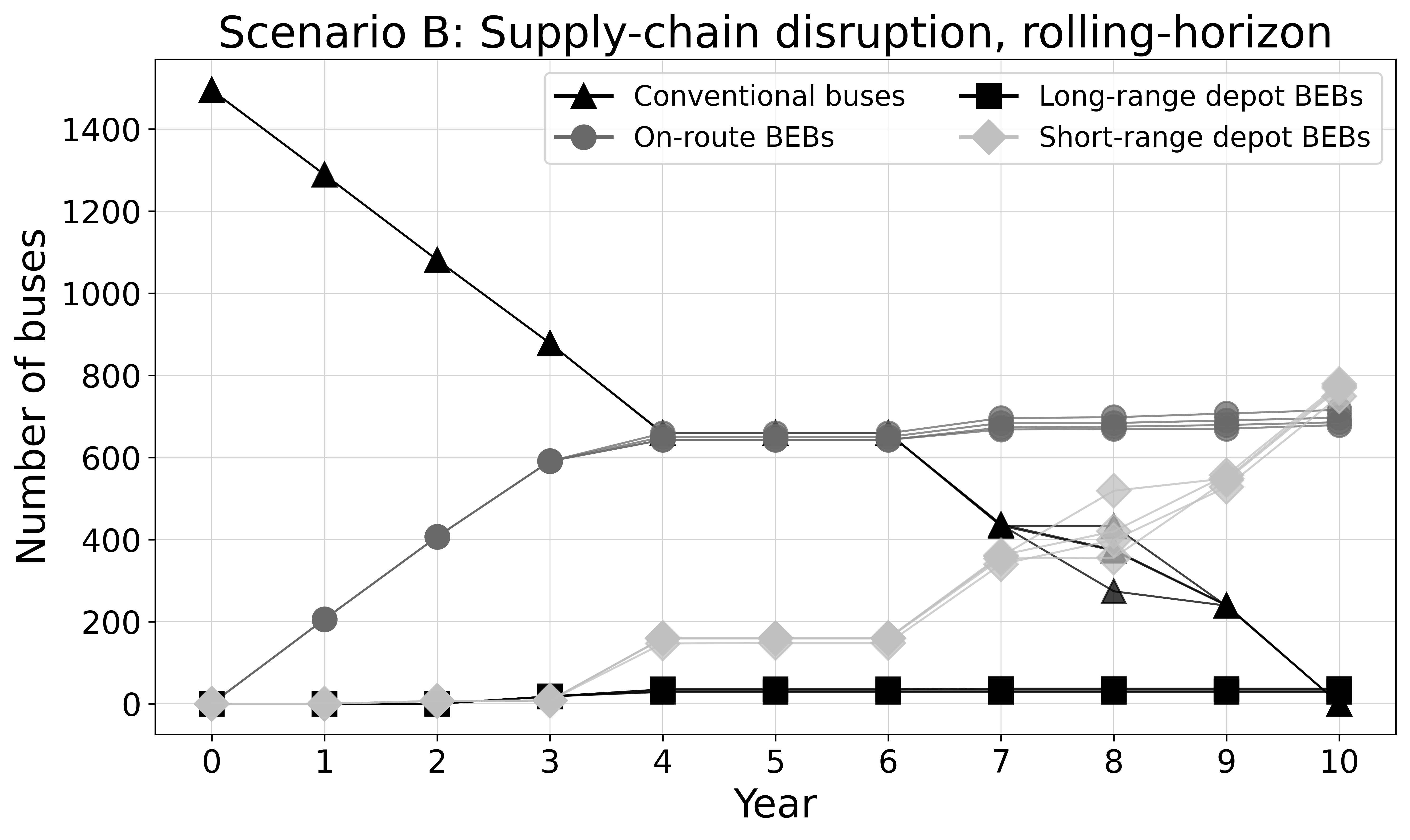}
    \end{subfigure}

    \begin{subfigure}{0.49\textwidth}
        \includegraphics[width=\linewidth]{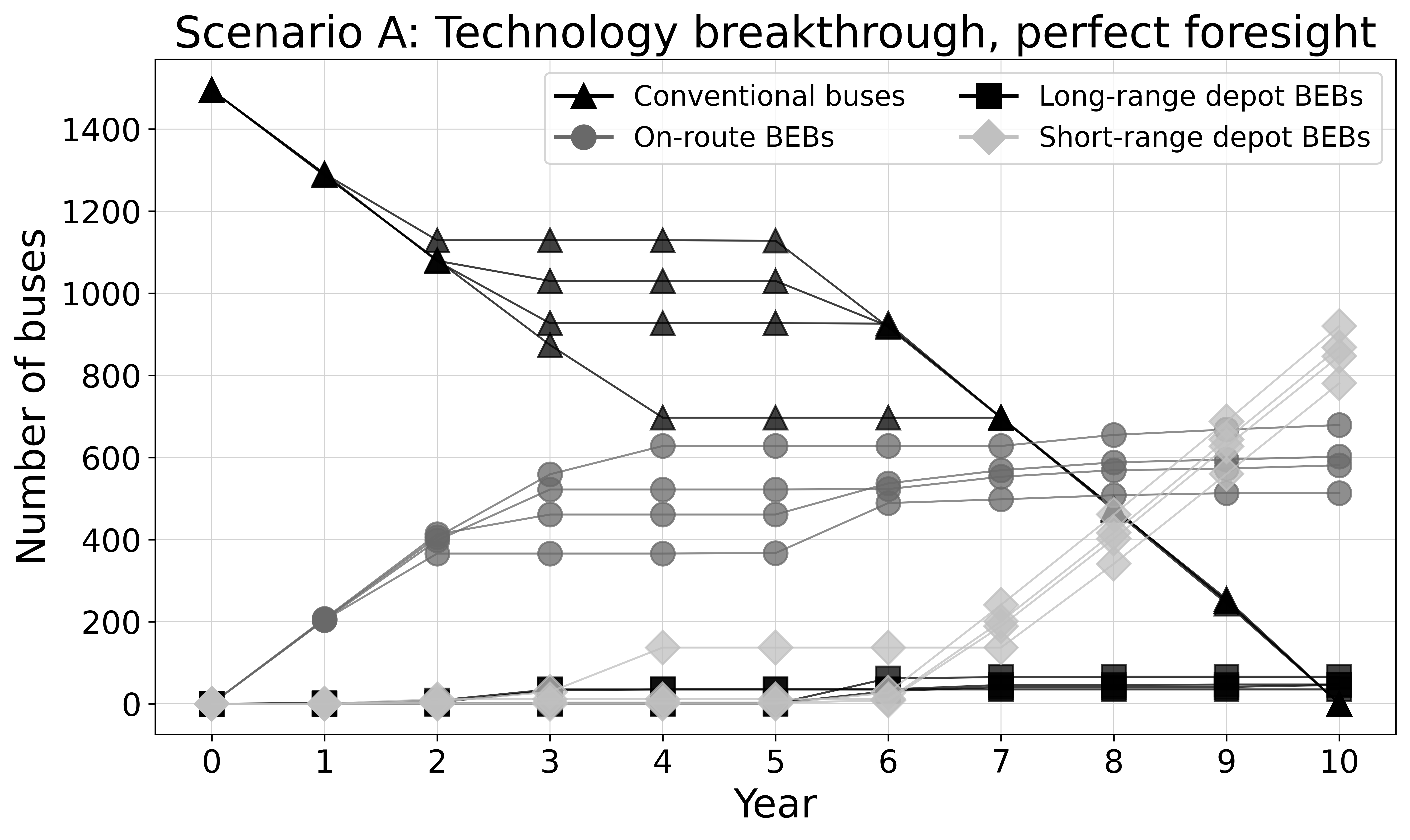}
    \end{subfigure}
    \hfill
    \begin{subfigure}{0.49\textwidth}
        \includegraphics[width=\linewidth]{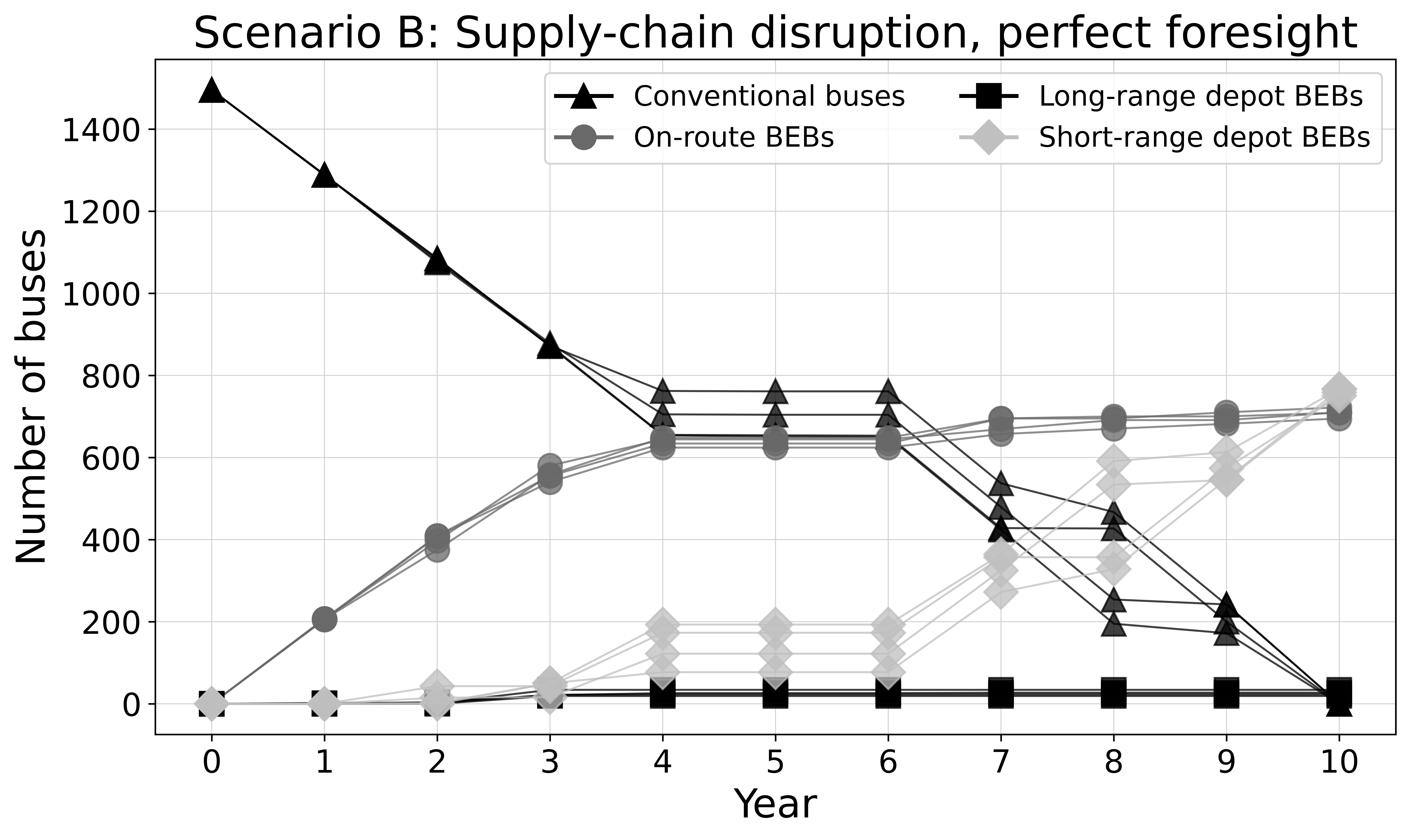}
    \end{subfigure}
\end{figure}

\vspace{-1.1cm}

\begin{table}[H]
\centering
\begin{revblock}
\caption{\begin{revblock}Average realized cost (in millions of dollars) and savings over baseline.\end{revblock}}
\setlength{\tabcolsep}{7pt}
\resizebox{\textwidth}{!}{
\renewcommand{\arraystretch}{0.75}
\begin{tabular}{l c cc cc}
  \toprule
  \multirow{2}{*}{Scenario} & \multicolumn{1}{c}{Baseline (no replan)} & \multicolumn{2}{c}{Rolling horizon} & \multicolumn{2}{c}{Perfect foresight} \\
  \cmidrule(lr){2-2} \cmidrule(lr){3-4} \cmidrule(lr){5-6}
  & \multicolumn{1}{c}{cost} & \multicolumn{1}{c}{cost} & \multicolumn{1}{c}{savings (\%)} & \multicolumn{1}{c}{cost} & \multicolumn{1}{c}{savings (\%)} \\
  \midrule
  A: Technology breakthrough & 2703.22 & 2609.55 & 3.47 & 2581.65 & 4.50 \\
  B: Supply-chain disruption & 2837.35 & 2734.44 & 3.63 & 2733.80 & 3.65 \\
  \bottomrule
\end{tabular}
\label{tab:simul_costs}
}
\end{revblock}
\end{table}

\vspace{-0.3cm}

\begin{revblock}
In conclusion, the rolling-horizon implementation provides a conceptually simple and computationally tractable alternative to explicitly optimizing the BFEP under uncertainty. In practical deployments, where input parameters are inevitably subject to revision, the results of Table~\ref{tab:simul_costs} suggest that the BFEP is best deployed in such a rolling-horizon scheme rather than used to define a one-shot commitment.
\end{revblock}
\newpage

\section{Notation} \label{ec:notation}
\begin{table}[htbp]
    \caption{Indices and sets, strategic variables, operational variables, parameters, value functions}
    \label{tab:notations}
    \centering
    \resizebox{1.0\textwidth}{!}{  
    \SingleSpacedXII
    \begin{tabular}{cl}
        \toprule
        Notation & Description \\
        \midrule
        \( p \in \mathcal{P} \) & Investment period (year); $\mathcal{P} = \{1,2,\dots,P\}$\\
        \( t \in \mathcal{T} \) & Operational time interval (hour); $\mathcal{T} = \mathbb{Z}/ T \mathbb{Z}$ the cyclic group of integers modulo $T$  \\
        \( r \in \mathcal{R} \) & Bus route \\
        \( b \in \mathcal{B} \) & Type of depot BEB \\
        \( i \in \mathcal{I} \) & Depot location \\
        \( j \in \mathcal{J} \) & On-route terminal \\
        \( s \in \mathcal{S}_b\) & Charging state for buses of type $b$; $\mathcal{S}_b = [0:s_b]$; $\mathcal{S}^w_b = [1:s_b]$ ; $\mathcal{S}^z_b=[0:s_b{-}1]$ \\
        \(\mathcal{R}(j) \) & Set of routes connected to terminal $j$ \\
        \(\mathcal{J}(r) \) & Set of terminals connected to route $r$ \\
        \midrule
        \( x_p \) & Strategic variables of period $p$; $x_p = \left(\chi_p,\{\eta_{pr}\}_{r \in \mathcal{R}}\right) \in \mathbb{Z}^n_+$  \\
        \( \chi_p \) & Chargers location decisions; $\chi_p = (\bar{\chi}^p, \widetilde{\chi}^p) \in \mathbb{Z}_+^{\mathcal{I}} \times \mathbb{Z}_+^{\mathcal{J}}$ \\
        \( \eta_{pr} \) & Fleet size decisions on route $r$; $\eta_{pr} = (\bar{\eta}^p_r, \widetilde{\eta}^p_r, \widehat{\eta}^p_r) \in \mathbb{Z}_+^{\mathcal{B}} \times \mathbb{Z}_+ \times \mathbb{Z}_+$ \\
        \( \bar{\chi}^p_{i} \) & \# of chargers at depot $i$ in period $p$ \\
        \( \widetilde{\chi}^p_{j} \) &  \# of fast on-route chargers at terminal $j$ in period $p$ \\
        \( \bar{\eta}^p_{rb} \) & \# of depot BEBs of type $b$ assigned to route $r$ in period $p$ \\
        \( \widetilde{\eta}^p_r \) &  \# of on-route BEBs assigned to route $r$ in period $p$ \\
        \( \widehat{\eta}^p_r \) &  \# of conventional buses assigned to route $r$ in period $p$ \\
        \midrule
        \multirow{2}{*}{\specialcell{$y_{pr}$}}  & Operational variables of period $p$ on route $r$;\\ & $y_{pr} = \left(w^p_r ,v^p_r ,z^p_r, \widetilde{w}^p_r ,\widehat{w}^p_r\right) \in \mathbb{Z}_+^{\mathcal{T}\times\mathcal{B}\times\mathcal{S}^w_b} \times
        \mathbb{Z}_+^{\mathcal{T}\times\mathcal{B}\times\mathcal{S}_b} \times
        \mathbb{Z}_+^{\mathcal{T}\times\mathcal{B}\times\mathcal{I}\times\mathcal{S}^z_b} \times
        \mathbb{Z}_+^{\mathcal{T}\times\mathcal{J}(r)} \times
        \mathbb{Z}_+^{\mathcal{T}}$  \\ 
        \( w_{rbs}^{p t} \) & \# of depot BEBs of type $b$ in service in state $s$ on route $r$ during interval $t$ of period $p$ \\
        \( v_{rbs}^{p t} \) & \# of depot BEBs of type $b$ idling in state $s$ on route $r$ during interval $t$ of period $p$ \\
        \multirow{2}{*}{\specialcell{$z_{rbis}^{p t}$}} & \# of depot BEBs of type $b$ assigned to route $r$ initiating a charging trip to \\ & \phantom{\#} depot $i$ from state $s$ in interval $t$ of period $p$ \\
        \multirow{2}{*}{\specialcell{$\widetilde{w}_{rj}^{p t}$}} & \# of on-route BEBs in service on route \begin{revblock}$r$ assigned to terminal $j$\end{revblock}\\ & \phantom{\#} \begin{revblock}for opportunity charging\end{revblock} in interval $t$ of period $p$ \\
        \( \widehat{w}_{r}^{p t} \) & \# of conventional buses in service on route \begin{revblock}$r$\end{revblock} in interval $t$ of period $p$ \\
        \midrule
        \( \gamma \) & Yearly discount factor \\
        \( c^{u}_{l} \) & Coefficient (coefficient vector) of component(s) $l$ of decision variables $u$ \\
        \( I_p^{\text{UB}} \) & Investment budget in period $p$ \\
        $\eta_p^{\text{LB}}$ & Minimum number of BEBs required in period $p$ (electrification target) \\
        $\widehat{\eta}_p^{\text{UB}}$ & Maximum number of conventional buses allowed in period $p$ (retirement target) \\
        $\bar{\chi}^{\text{UB}}_{i}$ & Maximum number of chargers that can be installed in depot $i$ \\
        $\widetilde{\chi}^{\text{UB}}_{j}$ & Maximum number of chargers that can be installed at terminal $j$ \\
        \( d_r^{p t} \) & Demand (minimum number of buses in service) on route $r$ in interval $t$ of period $p$ \\
        \( \rho \) & \# of on-route BEBs allowed to share the same charger during a service interval \\
        $\kappa_{rbis}$ &  \# of time intervals needed to perform a round-trip between route $r$ and depot $i$ \\ & \phantom{\#} and fully charge a BEB of type $b$ starting from state $s$ \\
        \midrule
        $I_p(x_p{-}x_{p{-}1})$ & Investment costs of period $p$ \\
        $O_p(x_p)$ & Total operational costs of period $p$; $O_p(x_p) = H_p(x_p) + \mathcal{Q}_p(x_p)$ \\
        $H_p(x_p)$ & Fixed maintenance costs of period $p$ given $x_p$ (do not depend on operational decisions) \\
        $\mathcal{Q}_p(x_p)$ & Optimal variable operational costs of period $p$ given $x_p$ (depend on operational decisions) \\
        \bottomrule
    \end{tabular}
    }  
\end{table}

\end{document}